\documentclass[12pt]{article}
\usepackage{pstricks}
\usepackage{amsfonts}
\usepackage{amsmath, graphicx, epsfig,psfrag, verbatim}
\usepackage{eucal,  amssymb,amsmath,amsthm,graphicx, amsfonts, latexsym,xcolor}
\usepackage{tikz}
\usepackage[utf8]{inputenc}

\newtheorem{theorem}{Theorem}[section]
\newtheorem{lem}{Lemma}[section]
\newtheorem{cor}{Corollary}[section]

\def\E{{\rm E}}
\def\P{{\rm P}}

\theoremstyle{remark}
\newtheorem{remark}[theorem]{Remark}

\newcommand{\eqd}{\stackrel{D}{=}}
\newcommand{\sge}{\stackrel{{\rm st}}{\ge}}
\newcommand{\sle}{\stackrel{{\rm st}}{\le}}

\newcommand{\convp}{\stackrel{p}{\longrightarrow}}

\def\mcE{\mathcal{E}}
\def\mcEt{\tilde{\mathcal{E}}}

\def\re{{\rm e}}
\def\var{{\rm var}}
\def\cov{{\rm cov}}

\def\convas{\stackrel{{\rm a.s.}}{\longrightarrow}}
\def\convD{\stackrel{{\rm D}}{\longrightarrow}}
\def\convw{\stackrel{{\rm w}}{\longrightarrow}}

\numberwithin{equation}{section}

\newcommand{\bmalpha}{\mbox{\boldmath{$\alpha$}}}

\begin{document} \parskip=5pt plus1pt minus1pt \parindent=0pt
\title{Household epidemic models revisited}
\author
{Frank Ball, Tom Britton and Peter Neal }
\date{\today}
\maketitle

\begin{abstract}
We analyse a generalized stochastic household epidemic model defined by a bivariate random variable $(X_G, X_L)$, representing the number of global and local infectious contacts that an infectious individual makes during their infectious period. Each global contact is selected uniformly among all individuals and each local contact is selected uniformly among all other household members. The main focus is when all households have the same size $h \geq 2$, and the number of households is large. Large population properties of the model are derived
including a central limit theorem for the final size of a major epidemic, the proof of  which utilises an enhanced embedding argument. A modification of the epidemic model is considered where local contacts are replaced by global contacts independently with probability $p$. We then prove monotonicity results for the probability of the major outbreak and the limiting final fraction infected $z$ (conditioned on a major outbreak). a) The probability of a major outbreak is shown to be increasing in both $h$ and $p$ for any distribution of $X_L$.  b)  The final size $z$ increases monotonically with both $h$ and $p$ if the probability generating function (pgf) of $X_L$ is log-convex, which is satisfied by traditional household epidemic models where $X_L$ has a mixed-Poisson distribution. Additionally, we provide counter examples to b) when the pgf of $X_L$ is not log-convex.

\end{abstract}


\section{Introduction}
The spreading of an infectious disease in a community is highly affected by how individuals mix. The earliest epidemic models assumed homogeneous mixing between all individuals, but this has later been generalized to allow for e.g.\ a multitype community (where mixing rates depend on the types of the two involved individuals involved), social networks (where individuals are connected through some underlying social structure, often characterized by a degree distribution for how many acquaintances individuals have), household epidemic models (where the population is divided into small units with higher mixing within the units) and spatial models (where the rate at which people mix depend on their spatial distance from each other), or a combination of these heterogeneous mixing features.

Stochastic household epidemic models were first considered by McKendrick \cite{McKen26}, with increased interest in the past 30 years starting with Becker and Dietz \cite{BeckerD95} and Ball {\it et al.}\ \cite{BMST}. A common definition of such a model uses the SIR or SEIR model concept to define an infectious period $I$, possibly having a random duration, during which an infectious individual has infectious contacts on the global scale and within the household according to independent Poisson processes having rates $\beta_G$ and $\beta_L$ (global and local). Each global contact is with a uniformly selected individual from the entire population (including also household members for convenience) and each local contact is with a uniformly selected individual of the same household. This implies that, conditional on the duration of the infectious period of an infective $I=x$, the random number of (not necessarily unique) global and local contacts are Poisson distributed with means $\beta_G x$ and $\beta_L x$, respectively. The corresponding unconditional numbers of global and local contacts $(X_G, X_L)$ are hence mixed-Poisson distributed: $X_G\sim {\rm MixPo}(\beta_G I)$ and $X_L\sim {\rm MixPo}(\beta_L I)$, where both random variables depend on the same random variable $I$ (the infectious period). If the population size is large it is unlikely that an individual makes multiple global contacts with the same individual. However, within small households multiple infectious contacts with the same individual are common.

The time dynamics of this household epidemic depends on the infectious period $I$ and its possible preceding latent period $L$, but the final size describing who eventually gets infected is independent of $L$ and only depends on $I$ through $(X_G, X_L)$. In the present paper we study a more general model where $(X_G, X_L)$ may follow an arbitrary but specified distribution on $\mathbb{Z}_+^2$. Hence, $X_G$ and $X_L$ need not be mixed-Poisson, nor do they have to be positively correlated as in the traditional model. In fact, it is quite possible that someone who becomes ill soon after infection makes fewer global infectious contacts but on the other hand more local (household) infectious contacts, thus making $X_G$ and $X_L$ negatively correlated. We consider also a modification of the model where each local contact is replaced by a global contact with probability $p$, in order to analyse what happens with the epidemic as more of the contacts become global,  i.e.\ $p$ increases.

The focus of this paper is the distribution of the final size of the epidemic as the population size, $N$, tends to infinity in the case where all households have the same size $h \ge 2$. We extend limiting results for the traditional epidemic models to our more general household epidemic model: a branching process approximation of the initial stages of the epidemic, an expression for the basic reproduction number $R_\ast$, and a law of large number and central limit theorem for the final size, conditional on the epidemic taking off. The central limit theorem (Theorem \ref{thm:clt}) uses an embedding argument based on the approach introduced in Scalia-Tomba \cite{ST85} and successfully applied to household epidemics in \cite{BMST}. However, given that the number of global infectious contacts made by infectives is not necessarily mixed-Poisson an additional layer of embedding is required to allow for a general distribution for global contacts, necessitating a novel proof. An explicit, and relatively easy to compute, expression for the variance of the central limit theorem is given in \eqref{eq:sigma} with details on numerical computation given in Appendix \ref{app:Calc}. 
In the absence of local infection, $X_L \equiv 0$, the model behaves as a homogeneously mixing epidemic and Theorem \ref{thm:clt} holds for extensions of the Reed-Frost model considered in Martin-L\"of \cite{ML86} and Picard and Lef\`evre \cite{PL90}, with Theorem \ref{thm:clt} corresponding to \cite{ML86}, Theorem 1. 
Provided that $\mu_G = \E [X_G] < \infty$, in the limit as $N \to \infty$, we obtain the same asymptotic final size distribution whether the global contacts made by an individual are {\it with} or {\it without} replacement. The standard household model assumes the local contacts are made with replacement and such a model is the main focus of this paper. However, the central limit theorem given in Theorem \ref{thm:clt} holds if instead we assume that local contacts are made without replacement through minor modifications to the arguments.

We provide novel insight into how the household size, $h$, and the probability, $p$, that a local contact is replaced by a global contact, affect the probability that the epidemic takes off (a major epidemic occurs) and the (asymptotic) final size of a major epidemic. For larger households (increasing $h$) and a greater proportion of global contacts (increasing $p$), the epidemic more closely resembles a homogeneously mixing epidemic with fewer multiple contacts made by an infective with the same individual. Hence, intuitively the probability that the epidemic takes off and the final size of a major epidemic are increasing in both $h$ and $p$. In Theorem  \ref{thm:prob:major}, we show that this is the case for the probability a major epidemic outbreak regardless of the choice of $(X_G, X_L)$. The effect of $h$ and $p$ on the final size of the epidemic depends on the distribution of $X_L$, with $X_G$ only affecting the final size through its mean $\mu_G$. Specifically, in Theorem \ref{thm:finalsize}, we show that the final size of the epidemic is increasing in both $h$ and $p$ if the logarithm of the probability generating function (pgf) of $X_L$ is convex. This is the case if $X_L$ follows a mixed-Poisson distribution, so the monotonicity results hold for the standard construction of the SIR household model. However, for more general $X_L$ the situation is more complex, with scenarios where the counter-intuitive result holds of smaller household sizes and increased local infectious contacts (with repeated contacts) leading to a larger final size.

The remainder of the paper is structured as follows. In Section \ref{sec:main:result}, we present the general stochastic household epidemic model and state the main results of the paper, a central limit theorem for the final size of the epidemic (Theorem \ref{thm:clt}) and sufficient conditions for the probability of a major outbreak (Theorem \ref{thm:prob:major}) and the final size of the epidemic (Theorem \ref{thm:finalsize}) to be increasing in $h$ and $p$. In Section \ref{sec:illus}, we present numerical illustrations of the main results, demonstrating the usefulness of the central limit theorem for finite $N$ and providing  examples where the final size of the epidemic is not increasing in $h$ and/or $p$. The proofs of the central limit theorem and of the effects of $h$ and $p$ on the probability and final size of a major outbreak are given in Sections \ref{sec:clt} and \ref{Sec-proof}, respectively. In Section \ref{Sec-discussion}, we discuss the findings of the paper and possible extensions. Finally, in the appendices we present details of how to compute key quantities such as the probability of a major outbreak, the final proportion infected and the variance of the final size (Appendix \ref{app:Calc}), along with  Appendices \ref{app:pnear1} and \ref{app:hinf}, which provide the proofs of Theorems \ref{thm:pnear1} and \ref{lem:hinf} concerning how the final size of a major epidemic behaves near $p=1$ (almost all local contacts replaced by global contacts) and as $h \to \infty$, respectively.

\section{Model and main results}\label{sec:main:result}

\subsection{The general household epidemic model}

The main ingredient for our epidemic model is the bivariate random variable $(X_G, X_L)$ with distribution on $\mathbb{Z}_+^2$.
$X_G$ and $X_L$, respectively, denote the number of global and local contacts that a randomly selected individual makes.

Consider a population consisting of $n$ households, all having size $h$. We investigate the limiting situation where the population size $N=nh$ tends to infinity in such a way that the household size $h$ remains fixed and the number of household $n\to\infty$.

An individual who gets infected draws their random pair $(X_G, X_L)$. Each of the $X_G$ global infectious contacts is with a uniformly selected individual from the entire population. Each of the $X_L$ local infectious contact is selected uniformly among the other $h-1$ household members. All contact selections are made independently, and a susceptible individual who is contacted gets infected (and repeats the procedure), whereas contacts with previously infected individuals have no effect. It is worth pointing out that the contacts of an individual may not all be to unique individuals. In particular, the $X_L$ local contacts may very well include multiple contacts to some individual(s). Such multiple contacts have no effect on the epidemic - it is the number of unique contacts that determines the propagation of the epidemic.

We consider also a modification of the model containing an additional parameter $p$. In this model, each local contact is, independently of everything else, replaced by a global contact with probability $p$. Thus, as $p$ increases, there are fewer local and more global contacts.

The epidemic is initiated by a number of individuals, chosen uniformly at random from the population, being infected and all other individuals being uninfected and susceptible. The epidemic continues until it eventually stops by no new individuals getting infected. The final number infected is denoted $Z$, or $Z_{n,h}$ if we want to emphasize its dependence on the number and size of households. Clearly $1\le Z\le N (=nh)$.

We denote the original model by $\mathcal{E}_{n,h}  (X_G, X_L)$ and the model with swapping of local contacts to global contacts by $\mathcal{E}_{n,h}  (X_G, X_L,p)$, so $\mathcal{E}_{n,h}  (X_G, X_L)$ is identical to $\mathcal{E}_{n,h}  (X_G, X_L,0)$. Note that $\mathcal{E}_{n,h} (X_G, X_L,1)$ is a homogeneously mixing epidemic.

\subsection{Relation to traditional household epidemic models} \label{sec:main:trad}

As described in the introduction, traditional household epidemic models are often defined by infectious individuals having a random  infectious period $I$, during which the infective has global contacts at rate $\beta_G$ and household contacts at rate $\beta_L$ (or $(h-1)\beta_L$, so $\beta_L$ to each household member, but we choose the former parametrisation). In that case, the numbers of global and local contacts have distribution $(X_G, X_L)=({\rm MixPo} (\beta_G I), {\rm MixPo}(\beta_L I))$, where we note that the two random variables are dependent having parameter containing the same random variable $I$. The final size of the epidemic depends only on the distribution of $(X_G, X_L)$, so the traditional model can be viewed as a subclass of $\mathcal{E}_{n,h}  (X_G, X_L)$.

\subsection{The $\mathcal{E}_{n,h}  (X_G, X_L,p)$ model described as an $\mathcal{E}_{n,h}  (X_G^\prime, X_L^\prime)$ model}\label{sec:main:p}

It is worth mentioning that $\mathcal{E}_{n,h}  (X_G, X_L,p)$ can, for a fixed value $p$, be described by $\mathcal{E}_{n,h}  (X_G^\prime, X_L^\prime)$, i.e.~the model without swapping, where the new random vector $(X_G^\prime, X_L^\prime)$ is different from the original vector $(X_G, X_L)$. More precisely, the new vector is simply the (random) number of global and local contacts that occur \emph{after} the swapping has happened. Suppose that  $X_L=k$ and let $Y_L\sim {\rm Bin}(k,p)$ denote how many contacts are swapped, then $X^\prime_G=X_G+Y_L$ and $X^\prime_L=X_L-Y_L$. Unconditionally, and showing the dependence on $p$, we hence have
$$
(X^{(p)}_G, X^{(p)}_L)= (X_G + Y^{(p)}_L,\ X_L -Y^{(p)}_L ) \text{, where }Y^{(p)}_L \sim {\rm MixBin} (X_L, p).
$$
Note that, in the expressions above, $Y^{(p)}_L$ depends on $X_L$ which is evident from the mixed-binomial distribution but hidden when writing the random vector $(X_G + Y^{(p)}_L,\ X_L -Y^{(p)}_L )$.

\subsection{Main results for the general household epidemic model}
\label{sec:SSHEmain:result}
We  now state our main results, firstly for the $\mathcal{E}_{n,h}  (X_G,X_L)$ model and then for the $\mathcal{E}_{n,h}  (X_G,X_L,p)$ model. These results are asymptotic results as $n \to \infty$ and for fixed $h$, we consider a sequence of epidemics, indexed by the number of households $n$.  The epidemic $\mathcal{E}_{n,h}  (X_G,X_L)$ is initiated by $m_n$ individuals, chosen uniformly at random from the population, being infected, with the remaining $nh-m_n$ individuals being susceptible.
Let $\bar{Z}_{n,h} = (nh)^{-1} Z_{n,h}$ denote the proportion of the population infected in  $\mathcal{E}_{n,h}  (X_G,X_L)$ and let $V_{n,h}$ denote the number of households where at least one individual is infected. 
Let $\mathcal{G}^{n,h} = \{ V_{n,h} \geq \lfloor \log n \rfloor\}$, the event that the epidemic infects at least $k_n = \lfloor \log n \rfloor$ households. We say that a major epidemic has occurred if $\mathcal{G}^{n,h}$ occurs. The choice of $k_n = \lfloor \log n \rfloor$ households being infected to define a major epidemic is somewhat arbitrary and the results in this paper hold for any sequence $k_n$ such that $k_n \to \infty$ and $k_n/\sqrt{n} \to 0 $ as $n \to \infty$. 

Before stating Theorem \ref{thm:clt}, which extends known results for the traditional household epidemic models (where $(X_G, X_L)=({\rm MixPo} (\beta_G I), {\rm MixPo}(\beta_L I))$) to a general random vector $(X_G, X_L)$,
we require some extra notation.

Consider a household of size $h$, with initially 1 infective and $h-1$ susceptibles. Let $\mcE_h^H(X_G, X_L)$ denote the ensuing {\it within-household} epidemic in which infected individuals make global and local infections according to the random pair $(X_G, X_L)$.  Let $C$ denote the number of global contacts that emanate from $\mcE_h^H(X_G, X_L)$.  Let $S$ denote the size of the susceptibility set of a typical individual in the household, where the susceptibility set of a given individual is the set of individuals, including themselves, who if infected globally will lead to the chosen individual being infected locally.  A formal definition is given in Section \ref{sec:clt:Sellke}.  Note that $S$ has support $\{1,2,\ldots,h\}$.  Let
\begin{align}
\label{eq:sus:pgf}
f_S (s) = \sum_{k=1}^h \P (S =k) s^k \qquad (0 \le s \le 1)
\end{align}
and
\begin{align*}
f_C (s) = \sum_{k=0}^\infty \P (C =k) s^k  \qquad (0 \le s \le 1)
\end{align*}
denote the pgfs of $S$ and $C$, respectively. Note that the distributions of $S$ and $C$ depend on $h$ but for notational convenience we suppress explicitly mentioning the dependence on $h$ unless it is the focus of our study. Let $R_\ast = \E [C]$, the mean number of global contacts emanating from a household epidemic. Then, letting $\mu_G = \E [X_G]$, it is straightforward (see the appendix of Ball {\it et al.}\ \cite{BMST}) to show that
\[
R_\ast = \E[C] = \mu_G \E [S].
\]

We now consider the household exposed to global  infection.  For $\pi \in [0,1]$, let $\mcEt^H_h(X_G, X_L, \pi)$ denote the following epidemic.  Initially the whole household is susceptible.  During the course of the epidemic, individuals avoid external infection independently with probability $\pi$.  Infected individuals make global and local infections according to the random pair $(X_G, X_L)$. For $t \ge 0$, let $R(t)$ and $G(t)$ be respectively the total number  infected in the household and the total number of global contacts emanating from the household in $\mcEt^H_h(X_G, X_L, \re^{-t})$.

For $t \geq 0$, let
\begin{align}
\label{eq:nuR}
\nu_R (t) \left( = \frac{1}{h}\E [R(t)]\right) = 1- f_S (e^{-t}).
\end{align}
Suppose that $R_\ast>1$ and define $z$ to be the solution in $(0,1]$ of
\begin{align}
\label{eq:z}
z = 1- f_S (e^{-z \mu_G}) = \nu_R (\mu_G z).
\end{align}
(It is seen easily that $z$ exists and is unique, since $\nu_R (\cdot)$ is concave, $\mu_G\nu_R^\prime (0)=R_\ast$ and $\nu_R (\infty) =1$.)  Let
\begin{align}
\label{eq:sigma}
\sigma^2 &= \frac{1}{h}\left[(1+b(\tau)\mu_G)^2\var(R(\tau)) + b(\tau)^2h\nu_R(\tau)(\sigma_G^2-\mu_G)\right.\\
&\qquad+\left. 2b(\tau)(1+b(\tau)\mu_G)(\cov(R(\tau), G(\tau))-\mu_G \var(R(\tau)))\right],\nonumber
\end{align}
where $\sigma_G^2=\var(X_G)$, $\tau = \mu_G z$ and $b (t) = \nu_R^\prime (t)/[1 - \mu_G \nu_R^\prime (t)]$.
\begin{theorem} \label{thm:clt}
Suppose that $R_\ast >1$, and that there exists $m \geq 1$ such that $m_n=m$ for all sufficiently large $n$, and $a > 0$ such that $\E [X_G^{2+a}] <\infty$. Let $z>0$ be given by \eqref{eq:z}
and $\rho$ be the unique solution in $[0,1)$ of
\begin{align}
\label{eq:rho}
\rho =  f_C (\rho).
\end{align}
Then
\[ \bar{Z}_{n,h} \convD Z \qquad \mbox{as } n \to \infty, \]
where the random variable $Z$ has probability mass function
\begin{align}
\label{eq:probZ}
\P (Z=0) = 1 -\P (Z=z) = \rho^m.
\end{align}
Furthermore, there exists $0 < \sigma^2 < \infty$ given by  \eqref{eq:sigma}, such that
\begin{align}
\label{eq:thmclt1}
\left. \sqrt{n h} \left(\bar{Z}_{n,h} - z\right) \right| \mathcal{G}^{n,h} \convD {\rm N} (0, \sigma^2) \qquad \mbox{as } n \to \infty.
\end{align}
\end{theorem}

Theorem \ref{thm:clt} holds if instead $(X_G,X_L)$ are the numbers of unique individuals contacted by an infective in the population and their household, respectively. In this case $X_L$ has support $\{0,1,\ldots,h-1\}$ and corresponds to sampling local infectious contacts {\it without} replacement from the other members of the household. Sampling without replacement affects the distributions of $C$ and $S$ but does not otherwise affect the derivation of the central limit theorem. We discuss this in more detail in Section \ref{sec:clt:without}.

In Section \ref{sec:clt:variance} we give two alternative but equivalent expressions for $\sigma^2$.  Note that $z$ depends on the distribution of $X_G$ only through its mean $\mu_G$.  In Appendix~\ref{app:Calc}, we give expressions for $\E [R(t)]$, $\var(R(t))$, $\cov(R(t), G(t))$ and $f_C(s)$ in terms of Gontcharoff polynomials, which enable $\rho$, $z$ and $\sigma^2$ to be computed.

We now turn our attention to the $\mathcal{E}_{n,h}  (X_G,X_L,p)$ model.  
Theorems \ref{thm:prob:major} and \ref{thm:finalsize} analyse  $\pi^{(h,p)}$, the limiting probability of a major outbreak assuming a single initial infective, and $z^{(h,p)}$, the limiting final fraction getting infected in the event of a major outbreak,  in particular their dependence on $h$ and $p$ for a given vector $(X_G, X_L)$.  (To connect with Theorem~\ref{thm:clt}, note that in an obvious notation, $\pi^{(h,p)}=1-\rho^{(h,p)}$.)

\begin{theorem} \label{thm:prob:major}
The limiting probability of a major outbreak $\pi^{(h,p)}$ is monotonically  increasing in $h$ and $p$ for any random vector $(X_G, X_L)$.
\end{theorem}
This hence means that the probability of a major outbreak increases if households are larger and/or local contacts are replaced by global contacts, both features making the epidemic model becoming closer to homogeneously mixing.

The second theorem concerns the final outbreak size $z^{(h,p)}$ assuming a major outbreak has occurred. Here the result depends on the \emph{distribution} of $X_L$ and in particular how much randomness there is. To this end we define the pgf of $X_L$: $f_{X_L}(s)=\sum_{k=1}^\infty s^k \P(X_L=k)$.

\begin{theorem}  \label{thm:finalsize}
Assume that $\log (f_{X_L}(s))$ is convex on $0\le s\le 1$. Then the limiting final size $z^{(h,p)}$ is monotonically  increasing in $h$ and $p$ for any $X_G$ (dependent or independent of $X_L$).
\end{theorem}
The mixed-Poisson distribution has a log-convex pgf, so Theorem \ref{thm:finalsize} holds for the traditional household epidemic model. Log-convexity of the pgf of $X_L$ implies that $\sigma_L^2 \geq \mu_L$ where $\sigma_L^2 = \var(X_L)$ and $\mu_L = \E[X_L]$.
In Section \ref{Sec-counter}, we present counter examples to Theorem \ref{thm:finalsize} in the case where  $\sigma_L^2 < \mu_L$.   (See Theorem \ref{thm:pnear1} (a) below.)

The following theorem is proved in Appendix~\ref{app:pnear1}.  We define $z^{(h,p)}$ to be strictly increasing (decreasing) in $p$ near 1 if there exists $p_*^{(h)} \in [0,1)$ such that $z^{(h,p)}$ is strictly increasing (decreasing) in $p$ for $p \in [p_*^{(h)},1]$. 
For $\sigma_L^2<\mu_L$, let $z^*(\mu_L, \sigma_L^2)=1-\frac{\mu_L-\sigma_L^2}{3\mu_L^2}$ and note that $z^*(\mu_l, \sigma_L^2) \in (0,1)$ since $X_L$ takes values in $\mathbb{Z}_+$.  Let $\alpha=\mu_L+\mu_G$ and, for $\alpha>1$, let $z_{\rm hom}(\alpha)$ be the unique solution of $1-z=\re^{-\alpha z}$ in $(0,1)$.  Note that $z_{\rm hom}(\alpha)$ is the proportion infected by a major outbreak in a homogeneously mixing epidemic, where each individual makes on average $\alpha$ infectious contacts.
\begin{theorem}  \label{thm:pnear1}
\begin{itemize}
Suppose that $h \ge 2$ and $\alpha=\mu_G+\mu_L>1$, so $\mathcal{E}_{n,h} (X_G, X_L,1)$ is supercritical.
\item[(a)]
If $\sigma_L^2 \ge \mu_L$, then $z^{(h,p)}$ is strictly increasing in $p$ near 1.
\item[(b)]
Suppose that $\sigma_L^2 < \mu_L$.  Then $z^{(h,p)}$ is strictly increasing in $p$ near 1 if $z_{\rm hom}(\alpha)<z^*(\mu_L, \sigma_L^2)$
and strictly decreasing in $p$ near 1 if $z_{\rm hom}(\alpha)>z^*(\mu_L, \sigma_L^2)$.
\end{itemize}
\end{theorem}

Finally, we consider the final size $z^{(h,p)}$ in the limit as household size $h \to \infty$, with the proof given in  Appendix~\ref{app:hinf}.
\begin{theorem}  \label{lem:hinf} Suppose that $\alpha=\mu_G + \mu_L >1$.
Then for any $0 \leq p \leq 1$, $z^{(h,p)} \to z_{\rm hom}(\alpha)$ as $h \to \infty$.
\end{theorem}

\subsection{Counter examples to Theorem \ref{thm:finalsize} when  $\sigma_L^2 < \mu_L$}\label{Sec-counter}

In this section we provide simple counter examples showing that our main results are not necessarily true when $X_L$ has too little randomness.

\subsubsection{An example where final size decreases with household size}

Consider the simple case where $X_L\equiv 1$, meaning that all infected individuals have exactly one household contact, uniformly selected among all household neighbours, and some fixed $\mu_G$. Note that $\log (f_{X_L} (s)) = \log (s)$, so the pgf of $X_L$ is a concave function.
From \eqref{eq:z}, we know that the final size $z$ is given by the solution in $(0,1)$ of the equation $1-z=f_S ( \re^{- \mu_G z})$. 

We start with the case $h=2$. The susceptibility set is then identical to 2, since the other household member must contact the index locally. So $S\equiv 2$, and the right-hand side of the final size equation equals $e^{-2\mu_G z}$.

When $h=3$ the susceptibility set of an individual can in fact take only the values 1 or 3. The former if both housemates contact each other locally, and the latter otherwise.  Consequently, we have $\P(S_3=1)=0.25$ and $\P(S_3=3)=0.75$. The right-hand side of the final size equation then equals $0.25\re^{-\mu_Gz} + 0.75 \re^{-3\mu_Gz}$.

If we choose $\mu_G=2$ the final size equation for $h=2$ becomes $1-z=e^{-4z}$, with solution $z_2=0.980$. When $h=3$ the final size equation is $1-z=0.25e^{-2z}+0.75e^{-6z}$ with solution $z_3=0.961$, thus showing that $h=2$ gives a larger major outbreak than $h=3$.

\subsubsection{An example where moving local to global contacts lead to smaller final size}

For an example such that the final size decreases as local contacts are swapped to global contacts we continue the example from the previous subsection with $X_L\equiv 1$, $\mu_G=2$ and $h=2$. When $p=0$ we have the final size equation considered above, leading to final size $z_2=0.980$. If we swap \emph{all} local contacts to global contacts (so $p=1$) we simply have a homogeneous community where all individuals have $\mu_G=3$ global contacts. The final size equation is then $1-z=\re^{-3z}$, with solution $z=0.941$. So, if \emph{all} local  contacts are swapped to global contacts we get a \emph{smaller} outbreak, implying that the final size cannot increase monotonically with $p$ (in fact it decreases monotonically).

\section{Numerical illustrations}
\label{sec:illus}

\subsection{Accuracy of asymptotic approximations}
\label{subsec:accuracy}

\begin{figure}
\begin{center}
\includegraphics[width=7cm]{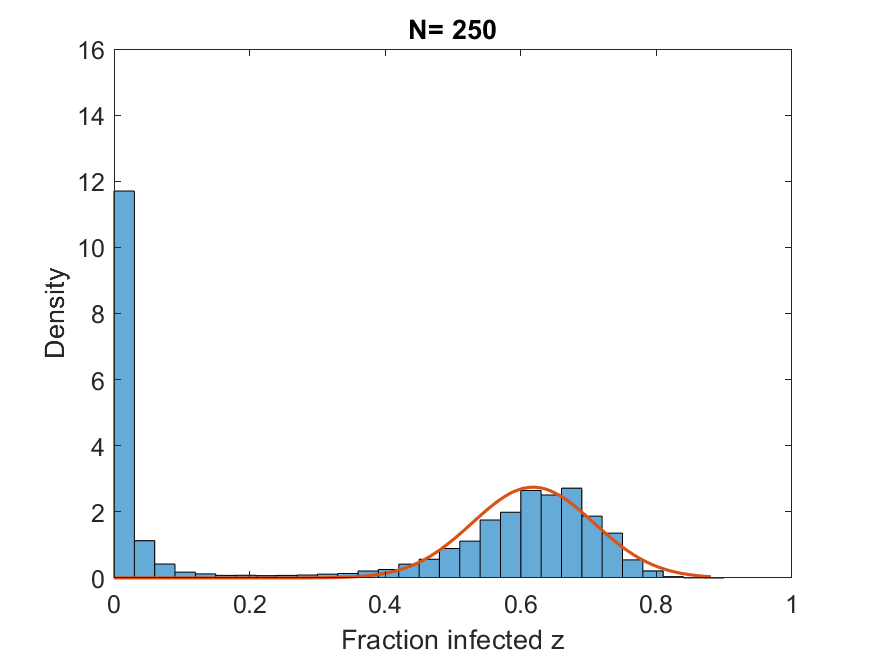}  \includegraphics[width=7cm]{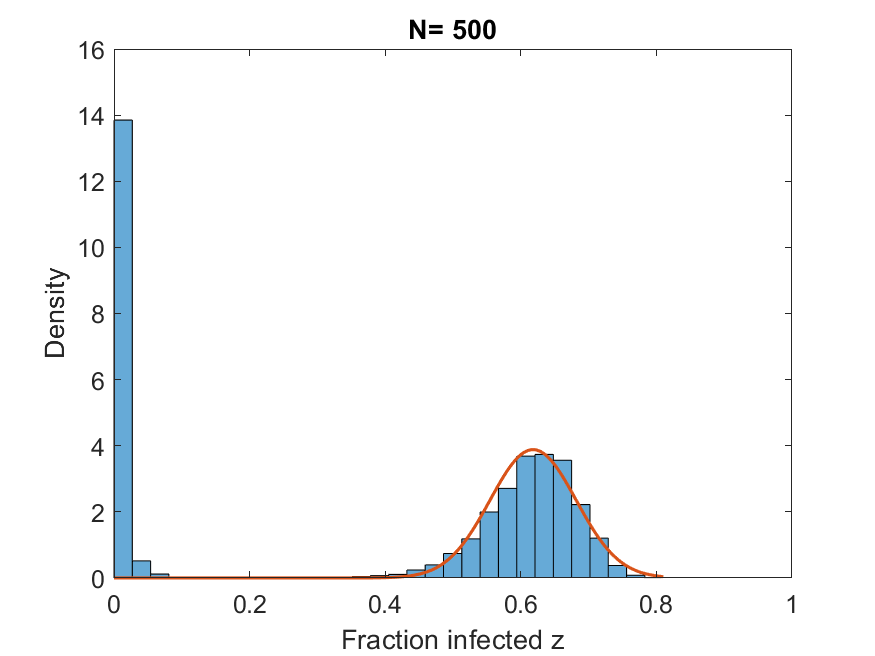}
\includegraphics[width=7cm]{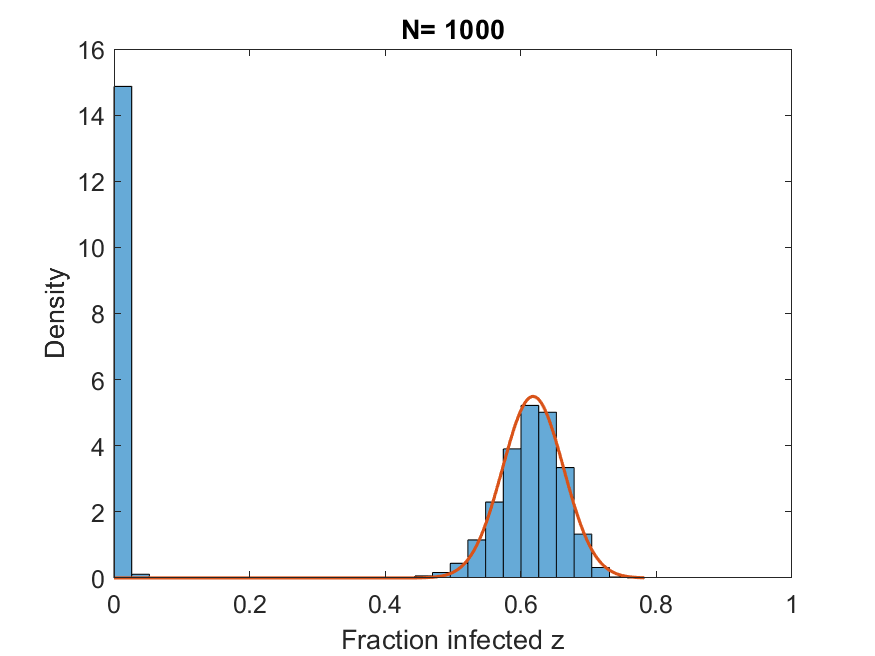}  \includegraphics[width=7cm]{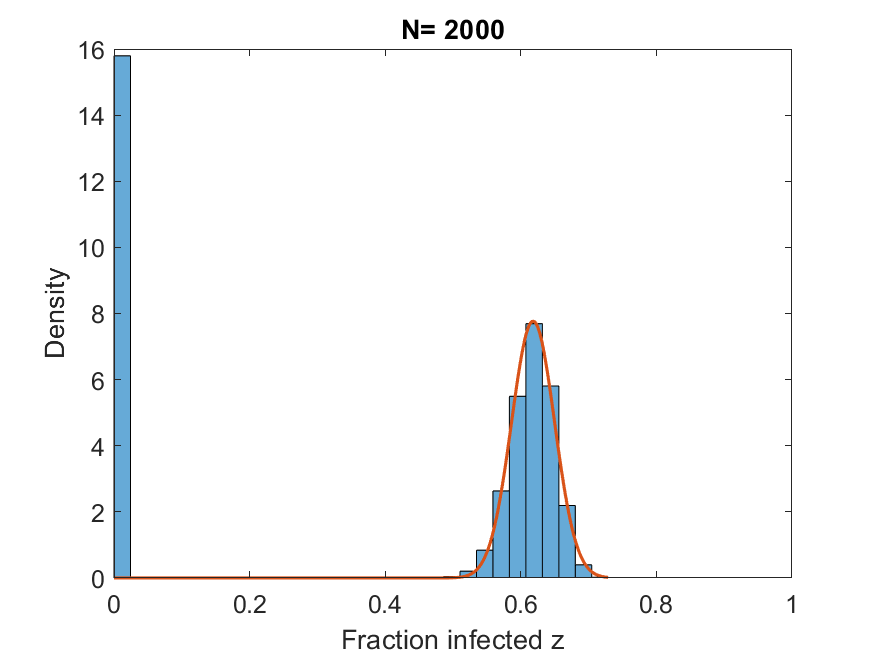}
\end{center}
\caption{Histograms of 100,000 simulations of the fraction of the population infected in $\mathcal{E}_{n,2} (X_G, X_L)$ when $X_G \sim {\rm Po}(1)$ and $X_L \sim {\rm Po}(1)$ independently, for population sizes $N=nh=250, 500, 1,000$ and $2,000$, with a normal approximation superimposed; see text for details.}
\label{fig:poisson}
\end{figure}

\begin{figure}
\begin{center}
\includegraphics[width=7cm]{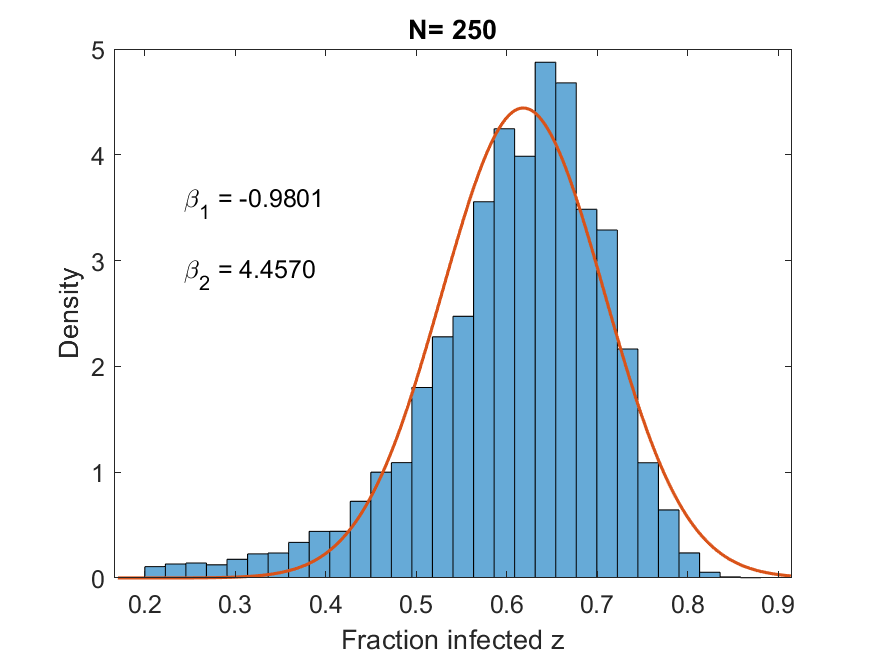}  \includegraphics[width=7cm]{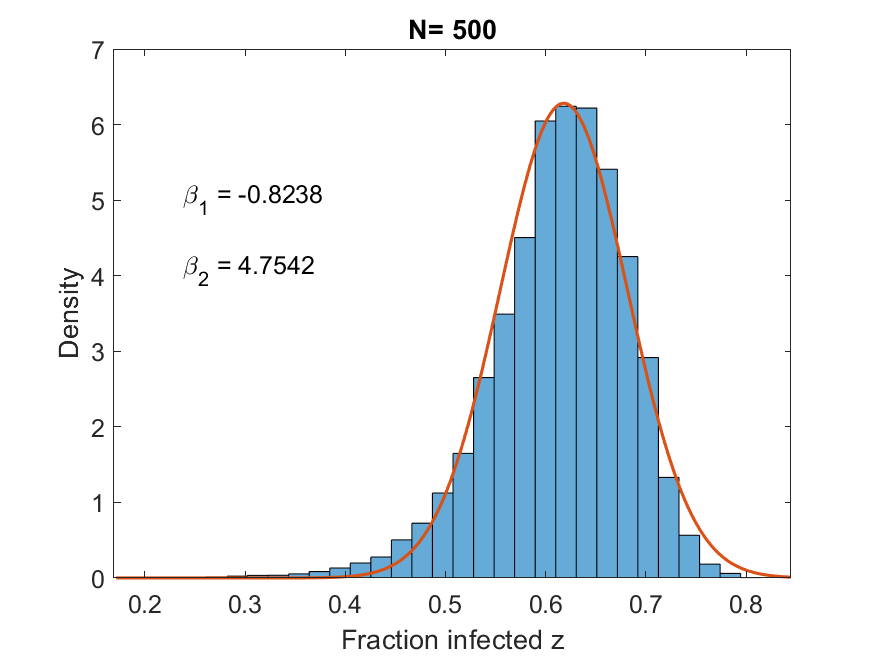}
\includegraphics[width=7cm]{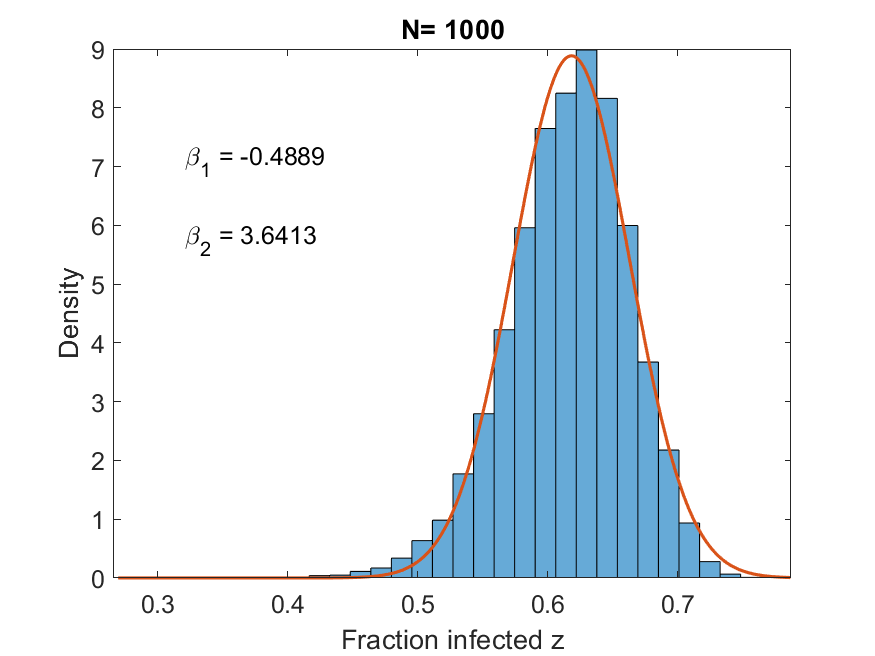}  \includegraphics[width=7cm]{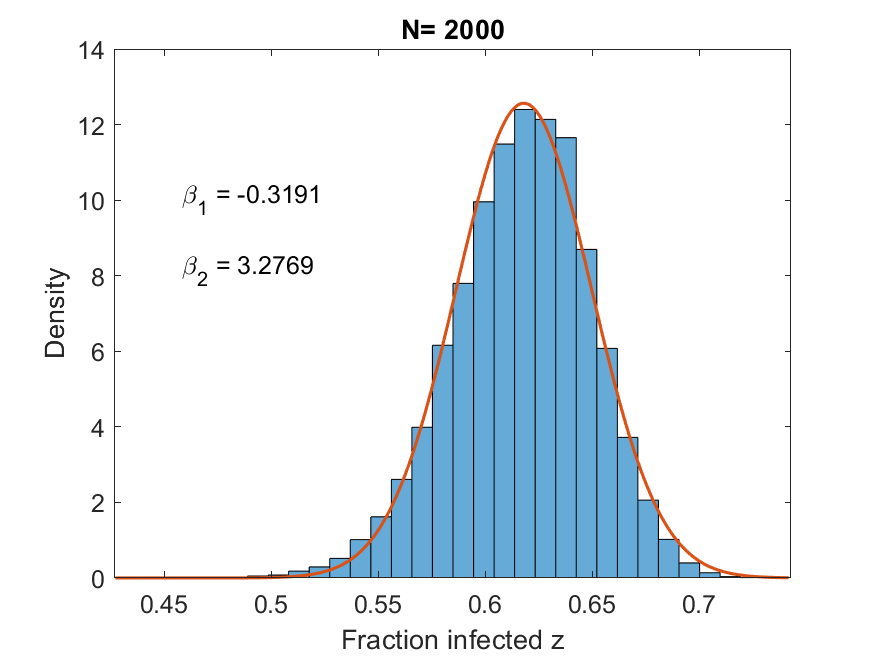}
\end{center}
\caption{Histograms of 100,000 simulations of the fraction of the population infected in a major outbreak in $\mathcal{E}_{n,2} (X_G, X_L)$ when $X_G \sim {\rm Po}(1)$ and $X_L \sim {\rm Po}(1)$ independently, for population sizes $N=nh=250, 500, 1,000$ and $2,000$, with a normal approximation superimposed; see text for details.}
\label{fig:poissonmaj}
\end{figure}

Figure~\ref{fig:poisson} shows histograms of the fraction of the population infected, $\bar{Z}_{n,h}$, in the epidemic $\mathcal{E}_{n,h}  (X_G,X_L)$ when $h=2$, $X_G \sim {\rm Po}(1)$ and $X_L \sim {\rm Po}(1)$ independently,
and $n=125, 250, 500$ and $1,000$ (so the total population size $N=250, 500, 1,000$ and $2,000$).  Each epidemic is initiated by a single infective and each histogram is based on $100,000$ simulations.  Superimposed on each histogram is the density $\pi^{(h,0)}f_N(x)$, where $f_N(x)$ is the probability density function of the normal distribution ${\rm N}(z, \frac{\sigma^2}{N})$, which approximates the distribution of $\bar{Z}_{n,h}$ for a major outbreak by Theorem~\ref{thm:clt}. For $N=1,000$ and $2,000$, there is a clear distinction between major and minor outbreaks.  The distinction is fairly clear for $N=500$ but not when $N=250$, where the choice of a cutoff to separate minor and major outbreaks is far from clear.  Figure~\ref{fig:poissonmaj} shows histograms of $100,000$ simulated major epidemics, using the same parameters as in Figure~\ref{fig:poisson} and a cutoff of $z=0.2$, with the ${\rm N}(z, \frac{\sigma^2}{N})$ probability density function superimposed.  Also shown are estimates of the skewness $\beta_1$ and kurtosis $\beta_2$ of the distribution of $\bar{Z}_{n,h}$ conditional upon a major outbreak.  (Note that $\beta_1=0$ and $\beta_2=3$ for a normal distribution.) The asymptotic normal distribution gives a good approximation for $N \ge 500$.  The true distribution of $\bar{Z}_{n,h}$ is skewed slightly to the left, with the degree of skewness decreasing as $N$ increases, and slightly more peaked than the asymptotic normal distribution. Note that Theorem~\ref{thm:clt} implies that, for any $z_* \in (0,z)$, the probability a major outbreak infects at least a fraction $z_*$ of the population tends to one as $n \to \infty$.  In the numerical study below, following inspection of histograms, we define a major outbreak to be one with $\bar{Z}_{n,h} \ge 0.2$.  Of course, the choice of cutoff depends on the parameters of an epidemic.

For a population of size $N$ consisting of households of size $h=2$, let $\pi_N$ be the major outbreak probability, and $z_N$ and $\sigma_N$ be the mean and scaled standard deviation of the fraction infected by a major outbreak.  (Thus $\sigma_N^2=N\var(\bar{Z}_{n,h}|\bar{Z}_{n,h} \ge 0.2)$, cf.~Theorem~\ref{thm:clt}.)
Table~\ref{table:hequals2} shows estimates of $\pi_N$, $z_N$ and $\sigma_N$ for the epidemic $\mathcal{E}_{n,h}  (X_G,X_L)$ with household size $h=2$ and various choices for the population size $N=nh$ and distribution for $(X_G, X_L)$.  For each choice of $N$ and distribution for $(X_G, X_L)$, $n_{\rm sim}=100,000$ epidemics were simulated and $\pi_N$ was estimated by $\hat{\pi}_N$, the fraction of simulations with $\bar{Z}_{n,2} > 0.2$, with an approximate 95\% confidence interval for $\pi_N$ given by $\hat{\pi}_N \pm 1.96 \sqrt{\hat{\pi}_N(1-\hat{\pi}_N)/n_{\rm sim}}$.  The simulations with $\bar{Z}_{n,h}\le 0.2$ were then discarded and further simulations made until there were $n_{\rm sim}$ simulations with $\bar{Z}_{n,h} > 0.2$, which were used to estimate
$z_N$ and the scaled standard deviation $\sigma_N$.  Let $\hat{z}_N$ and $\tilde{\sigma}_N^2$ be the sample mean and variance of these $n_{\rm sim}$ simulations of $\bar{Z}_{n,h}$.  Then $z_N$ was estimated by $\hat{z}_N$, with an approximate 95\% confidence interval given by $\hat{z}_N \pm  1.96 \tilde{\sigma}_N/\sqrt{n_{\rm sim}}$ and $\tilde{\sigma}_N$ was estimated by $\hat{\sigma}_N=\sqrt{N}\tilde{\sigma}_N$, with an approximate 95\% confidence interval given by $\left[\hat{\sigma}_N\sqrt{(n_{\rm sim}-1)/q_2},
\hat{\sigma}_N\sqrt{(n_{\rm sim}-1)/q_1}\right]$, where $q_1$ and $q_2$ are respectively the 2.5\% and 97.5\% quantiles of the $\chi^2_{n_{\rm sim}-1}$ distribution. The $N=\infty$ entries in Table~\ref{table:hequals2} give the asymptotic values $\pi, z$ and $\sigma$ given by Theorem~\ref{thm:clt}.

The distributions of $(X_G, X_L)$ in Table~\ref{table:hequals2} all have $\E[X_G]=\E[X_L]=1$ and are defined as follows.  Constant: $(X_G, X_L) \equiv (1,1)$.  Binomial: $X_G \sim {\rm Bin}(2, \frac{1}{2})$ and $X_L \sim {\rm Bin}(2, \frac{1}{2})$ independently.
Poisson: $X_G \sim {\rm Po}(1)$ and $X_L \sim {\rm Po}(1)$ independently.  Mixed-Poisson $I$: $X_G|I \sim {\rm Po}(I)$ and $X_L|I \sim {\rm Po}(I)$ independently, where $I$ is a single realisation of the given distribution.  Note that in Table~\ref{table:hequals2}, the distributions are listed in increasing order of $\var(X_G)$ and $\var(X_L)$.  There are no entries under $\hat{\pi}_N$ when $(X_G, X_L)$ has the Constant distribution since, then $\pi=1$ and for the values of $N$ considered, $\pi_N$ is extremely close to one.

It can be seen from Table~\ref{table:hequals2} that $\hat{\pi}_N$ generally increases with $N$ and $\pi$ is an overestimate of $\pi_N$  for finite $N$, as one would expect on intuitive grounds.  Further, the convergence of $\pi_N$ to its asymptotic value $\pi$ is faster when $X_G$ and $X_L$ have a smaller variance.  A similar comment holds for the fraction infected by a major outbreak $z$, though convergence of $z_N$ to $z$ is generally faster than that of $\pi_N$ to $\pi$.  Note that the confidence intervals for $z_N$ are smaller than those for $\pi_N$.  The simulations suggest that $\sigma$ is an underestimate of $\sigma_N$ and that the scaled
standard deviation of the size of a major outbreak converges to its asymptotic value more slowly than the mean.  Caution is required when interpreting results for small $N$, since then the distinction between major and minor outbreaks is less clear, particularly for distributions with larger $\var(X_G)$ and $\var(X_L)$.

The accuracy of the asymptotic normal distribution as an approximation for the size of a major epidemic in a finite population is explored further in Table~\ref{table:ks}, which is based on $n_{\rm sim}=100,000$ simulations for each choice of distribution for $(X_G, X_L)$, population size $n$ and household size $h$.  For each such choice, the table shows the value of the Kolmogorov-Smirnov one-sample test statistic $\displaystyle D_{n_{\rm sim}}=\sup_x|F_{n_{\rm sim}}(x)-F(x)|$, where $F_{n_{\rm sim}}$ is the empirical distribution function of the $n_{\rm sim}$ simulated fractions infected by a major outbreak and $F$ is the distribution function of the approximating ${\rm N}(z, \frac{\sigma^2}{N})$ distribution obtained using Theorem~\ref{thm:clt}.  Note that the corresponding tests all reject the null hypothesis that the fraction infected by a major outbreak follows a ${\rm N}(z, \frac{\sigma^2}{N})$ distribution, with a very low $p-$value, since the true distribution is not ${\rm N}(z, \frac{\sigma^2}{N})$ and the sample size $n_{\rm sim}$ is very large.  Nevertheless, the values of $D_{n_{\rm sim}}$ give a measure of the accuracy of the normal approximation. The values of $D_{n_{\rm sim}}$ clearly decrease with $N$, consistent with the convergence in Theorem~\ref{thm:clt}.  They also generally decrease with increasing household size $h$, though that is less clear for the Constant and Binomial cases.  Among the Poisson and mixed-Poisson choices for the distribution of $(X_G, X_L)$, the accuracy of the approximation generally decreases with increasing variance.  Overall, Table~\ref{table:ks} confirms the usefulness of the asymptotic normal approximation for finite population sizes.

\begin{table}
\small
\begin{center}
\begin{tabular}{|c|c|l|l|l|}
\hline
$(X_G, X_L)$ & $N$ & \qquad \qquad$ \quad\hat{\pi}_N$  &\qquad \qquad \quad $\hat{z}_N$ &\qquad \qquad \quad $\hat{\sigma}_N$ \\
\hline
& 250& &0.7960 \; (0.7957, 0.7962) &0.7471 \; (0.7438, 0.7504)\\
& 500& &0.7966 \;  (0.7964,    0.7968) &  0.7402 \;   (0.7369,    0.7434)\\
& 1,000& &0.7967 \;   (0.7965,    0.7968) &    0.7421 \;   (0.7389,    0.7454)\\
Constant  & 2,000& &0.7967 \;   (0.7966,    0.7968) &   0.7394 \;  (0.7362,    0.7427)\\
& 5,000 &&0.7967 \;   (0.7967,    0.7968) &    0.7384 \;     (0.7352,   0.7417)\\
& 10,000 &&0.7968 \;   (0.7967,    0.7968) &    0.7367 \;    (0.7335,    0.7399)\\
& $\infty$ &    &0.7968 & 0.7386  \\
\hline
& 250 &0.8103 \;   (0.8078,    0.8127) &    0.6762 \;   (0.6757,    0.6766) &    1.1642 \;    (1.1591,    1.1693) \\
& 500 &0.8181 \;   (0.8157,    0.8205) &   0.6794 \;   (0.6791,    0.6797) &    1.1138  \;  (1.1089, 1.1187)  \\
& 1,000 &0.8199 \;   (0.8175,    0.8223) &    0.6805 \;    (0.6803,    0.6808) &    1.0963 \;   (1.0915,    1.1011) \\
Binomial  &2,000 & 0.8230 \;    (0.8207,    0.8254) &   0.6812 \;   (0.6811,    0.6814) &    1.0902 \;    (1.0854,    1.0950)\\
& 5,000 & 0.8233 \;    (0.8210,   0.8257) &    0.6814 \;   (0.6813,    0.6815) &    1.0928 \;    (1.0881,    1.0976) \\
& 10,000 & 0.8238 \;   (0.8215,    0.8262) &    0.6816 \;    (0.6815,    0.6817) &    1.0852 \;   (1.0805,    1.0900) \\
 & $\infty$ &0.8238    &0.6817 & 1.0854  \\
 \hline
 & 250 & 0.5916 \;  (0.5885,    0.5946) &   0.6084 \;   (0.6078,    0.6091) &   1.5814 \;    (1.5745,    1.5884) \\
& 500 & 0.6053  \;  (0.6023,    0.6083)   & 0.6135 \;    (0.6131,    0.6139)   & 1.5249 \;   (1.5182,    1.5316) \\
& 1,000 &0.6126 \;   (0.6096,    0.6157   & 0.6159 \;    (0.6156,    0.6162)   & 1.4670 \;   (1.4606,    1.4735) \\
Poisson  &2,000 &0.6169 \;   (0.6139,  0.6199) &    0.6170 \;   (0.6168,    0.6172) &    1.4359 \;    (1.4297,    1.4423) \\
& 5,000 &0.6153 \;   (0.6123,    0.6183)  &  0.6178 \;   (0.6177,    0.6179)  &  1.4270 \;   (1.4208,    1.4333) \\
& 10,000 &0.6179 \;   (0.6149,    0.6209) &  0.6179 \;    (0.6178,    0.6180)  &  1.4196 \;   (1.4134,    1.4259)\\
 & $\infty$ &0.6181    &0.6181 & 1.4201  \\
 \hline
  & 250 & 0.3992 \;    (0.3962,    0.4023) &    0.5640 \;   (0.5633,    0.5648) &    1.9193 \;   (1.9110,    1.9278)\\
& 500 & 0.4127 \;   (0.4097,    0.4158) &  0.5661 \;   (0.5655,    0.5666) &    2.0076 \;   (1.9989,    2.0165) \\
Mixed-Poisson& 1,000 & 0.4252  \;  (0.4221,    0.4283) &    0.5687 \;   (0.5683,    0.5691) &    1.9666  \;  (1.9580,    1.9752)\\
$I \sim {\rm Gamma}(2,2)$  &2,000 &  0.4284 \;    (0.4253,    0.4314) &    0.5708 \;   (0.5705,    0.5710) &    1.8831 \;  (1.8749,    1.8914)\\
& 5,000 & 0.4316  \;  (0.4285,   0.4346) &    0.5718 \;   (0.5716,    0.5719) &    1.8508 \;   (1.8428,    1.8590)\\
& 10,000 & 0.4350 \;   (0.4319,    0.4381) &    0.5722 \;   (0.5721,    0.5723) &    1.8461 \;   (1.8381,    1.8542) \\
 & $\infty$ &0.4391    &0.5725 & 1.8378  \\
 \hline
   & 250 & 0.2933 \;   (0.2905,    0.2961) &    0.5357 \;   (0.5348,    0.5365) &    2.0870 \;   (2.0779,    2.0962)\\
& 500 & 0.3024 \;   (0.2995,    0.3052) &    0.5320 \;   (0.5313,    0.5326) &    2.3291 \;   (2.3190,    2.3394)\\
Mixed-Poisson& 1,000 & 0.3150  \;  (0.3122,    0.3179) &    0.5326 \;   (0.5322,    0.5331) &    2.4141 \;   (2.4035,    2.4247) \\
$I \sim {\rm Exp}(1)$  &2,000 & 0.3224 \;   (0.3195,    0.3253) &    0.5346 \;   (0.5343,    0.5349) &     2.3315 \;   (2.3213,    2.3417)\\
& 5,000 & 0.3254  \;  (0.3225,    0.3283) &     0.5359 \;   (0.5357,    0.5361) &    2.2697 \;   (2.2598,    2.2797) \\
& 10,000 & 0.3274  \;  (0.3245,    0.3303) &    0.5363 \;  (0.5362,   0.5365) &    2.2453 \;   (2.2355,    2.2552) \\
 & $\infty$ &0.3247    &0.5368 & 2.2347  \\
 \hline
   & 250 &0.1892 \;   (0.1868,    0.1917) &    0.5013 \;   (0.5004,    0.5021) &    2.2397 \;   (2.2299,    2.2496)\\
& 500 &  0.1838 \;    (0.1814,    0.1862) &    0.4900 \;   (0.4892,    0.4907) &    2.6346 \;   (2.6231,    2.6462)\\
Mixed-Poisson& 1,000 &0.1900 \;   (0.1876,    0.1924) &    0.4831 \;   (0.4825,    0.4837) &    2.9566 \;   (2.9437,    2.9696)\\
$I \sim {\rm Gamma}(\frac{1}{2},\frac{1}{2})$  &2,000 & 0.1984 \;   (0.1959,    0.2009) &    0.4810 \;   (0.4806,    0.4815) &    3.1439 \;   (3.1302,    3.1578)\\
& 5,000 & 0.2011 \;   (0.1986,    0.2036) &    0.4816 \;   (0.4813,    0.4819) &    3.0990 \;    (3.0854,   3.1126)\\
& 10,000 & 0.2044 \;   (0.2019,    0.2069) &    0.4819 \;   (0.4818,    0.4821) &   3.0495 \;   (3.0362,    3.0629)\\
 & $\infty$ &0.2060    &0.4829 & 2.9959  \\
 \hline
\end{tabular}
\end{center}
\caption{ Simulation results against
theoretical (asymptotic) calculations for epidemics with $h=2$.  See text for details.}
\label{table:hequals2}
\end{table}

\begin{table}
\small
\begin{center}
\begin{tabular}{|c|c|c|c|c|c|}
\hline
$(X_G, X_L)$ & $N$ & $h=2$ & $h=3$ & $h=4$ & $h=5$ \\
\hline
& 250&0.0477 &   0.0391 &   0.0455  &  0.0469\\
& 500& 0.0350 &   0.0293  &  0.0314  &  0.0332  \\
& 1,000& 0.0228  &  0.0196 &   0.0266 &   0.0226\\
Constant  & 2,000&  0.0168  &  0.0135 &   0.0174  &  0.0185\\
& 5,000 & 0.0092 &   0.0098 &   0.0149 &   0.0099\\
& 10,000 &0.0083 &   0.0081 &   0.0085 &   0.0070\\
\hline
& 250 &0.0303   & 0.0272   & 0.0320   & 0.0349 \\
& 500 &  0.0215  &  0.0204  &  0.0224  &  0.0233 \\
& 1,000 & 0.0154 &   0.0127 &   0.0200 &   0.0197\\
Binomial  &2,000 &0.0107 &   0.0100 &   0.0107 &   0.0116 \\
& 5,000 & 0.0085 &   0.0072 &   0.0071  &  0.0100\\
& 10,000 & 0.0048  &  0.0072 &   0.0057 &   0.0090\\
 \hline
 & 250 & 0.0414  &  0.0357  &  0.0314 &   0.0342 \\
& 500 & 0.0285   & 0.0244   & 0.0225  &  0.0224 \\
& 1,000 & 0.0193  &  0.0173 &   0.0169 &   0.0162 \\
Poisson & 2,000  & 0.0154 &   0.0115 &   0.0123 &   0.0123\\
& 5,000 & 0.0100  &  0.0081   & 0.0077   & 0.0110 \\
& 10,000 & 0.0076 &   0.0081   & 0.0062   & 0.0072\\
\hline
  & 250 & 0.0469 &   0.0479 &   0.0431 &   0.0371\\
& 500 & 0.0363   & 0.0322   & 0.0270   & 0.0267 \\
Mixed-Poisson& 1,000 &0.0276 &   0.0212 &   0.0196 &   0.0190\\
$I \sim {\rm Gamma}(2,2)$  &2,000 &0.0176  &  0.0149 &   0.0152 &   0.0149 \\
& 5,000 &  0.0115 &   0.0103  &  0.0098  &  0.0100\\
& 10,000 & 0.0102  &  0.0078 &   0.0067  &  0.0070 \\
\hline
   & 250 &0.0643  &  0.0542  &  0.0517  &  0.0493 \\
& 500 &0.0481   & 0.0387 &   0.0342  &  0.0324 \\
Mixed-Poisson& 1,000 & 0.0371  &  0.0267 &   0.0236  &  0.0224 \\
$I \sim {\rm Exp}(1)$  &2,000 &  0.0240  &  0.0187  &  0.0174  &  0.0162 \\
& 5,000 &   0.0152 &   0.0121 &   0.0118  &  0.0098\\
& 10,000 & 0.0110  &  0.0082  &  0.0079 &   0.0074 \\
\hline
   & 250 & 0.0923 &   0.0589  &  0.0597 &   0.0623\\
& 500 & 0.0797  &  0.0509   & 0.0461  &  0.0438  \\
Mixed-Poisson& 1,000 & 0.0616 &   0.0357  &  0.0297  &  0.0283\\
$I \sim {\rm Gamma}(\frac{1}{2},\frac{1}{2})$  &2,000 &0.0434 &   0.0256 &   0.0202 &   0.0217\\
& 5,000 & 0.0213  &  0.0164 &   0.0148  &  0.0122\\
& 10,000 & 0.0134  &  0.0118 &   0.0100  &  0.0094\\
\hline
\end{tabular}
\end{center}
\caption{Kolmogorov-Smirnov one-sample test statistics $D_{n_{\rm sim}}$ for testing the goodness-of-fit of the
approximating ${\rm N}(z, \frac{\sigma^2}{N})$ distribution, obtained using Theorem~\ref{thm:clt}, to a random sample of $n_{\rm sim}=100,000$ simulated major outbreaks for each parameter combination.  See text for details.}
\label{table:ks}
\end{table}

\subsection{Exploring model behaviour}
\label{subsec:modelbehaviour}
\begin{figure}
\begin{center}
\begin{tabular}{ccc}
  (a) Constant &
  (b) Binomial \\
 \includegraphics[width=7cm]{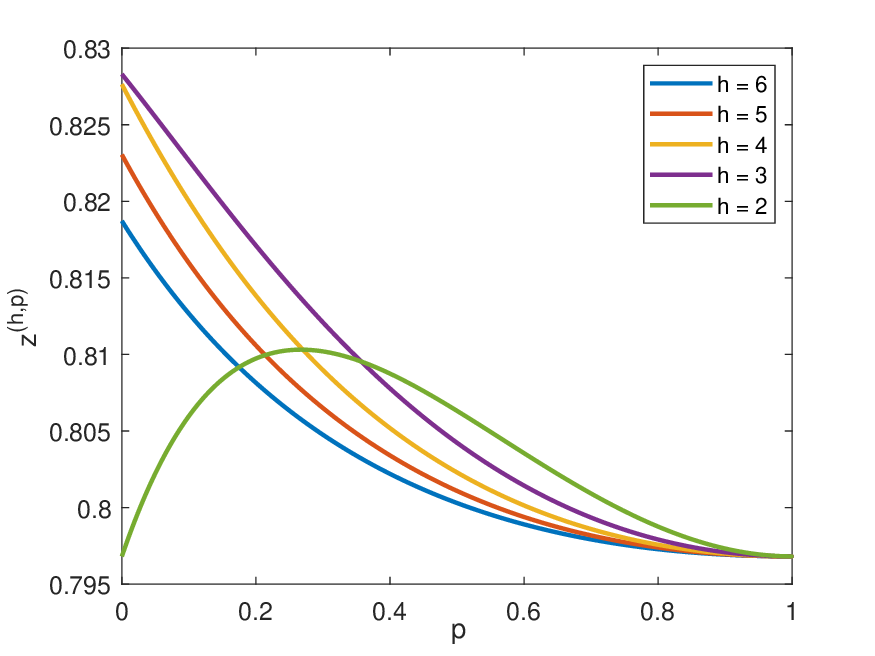} &
 \includegraphics[width=7cm]{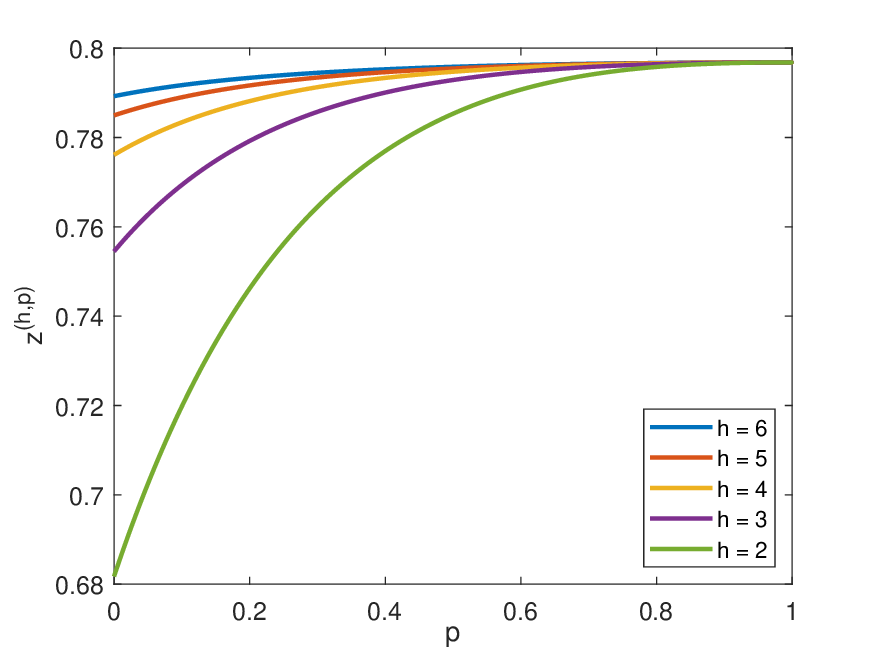} \\
  (c) Poisson &
  (d) Mixed-Poisson, $I \sim {\rm Exp}(1)$\\
 \includegraphics[width=7cm]{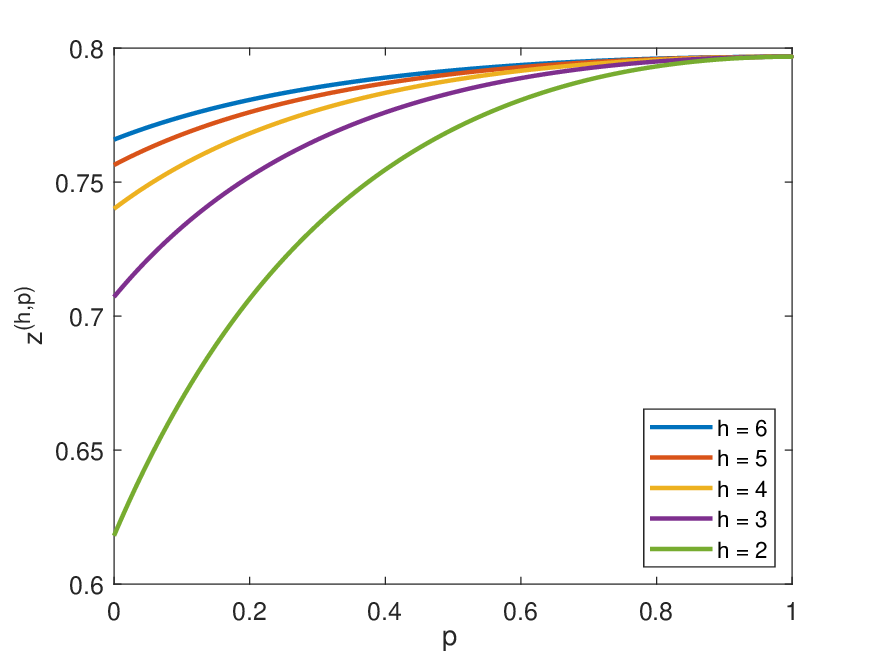} &
 \includegraphics[width=7cm]{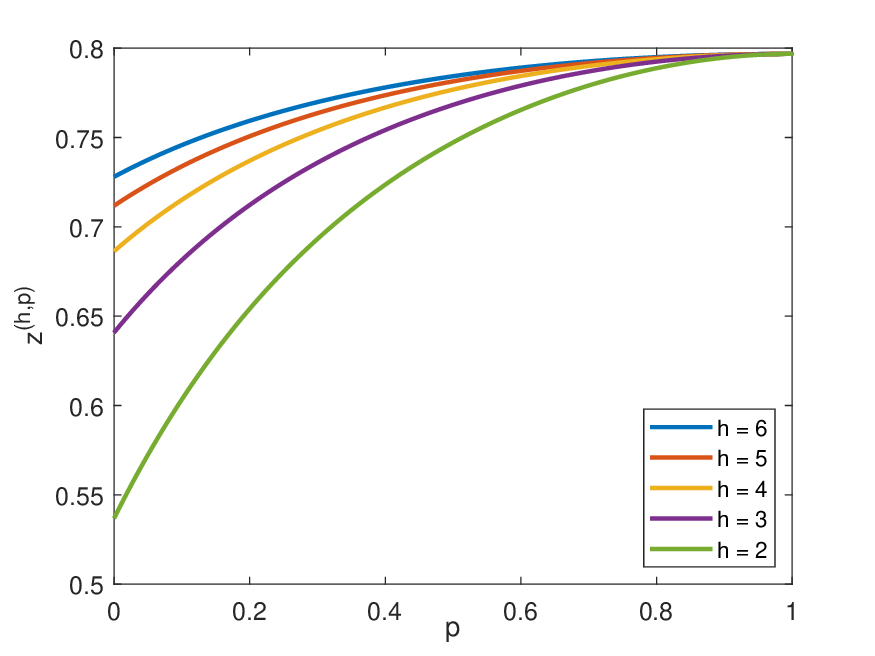}
\end{tabular}
\end{center}
\caption{Graphs of the fraction of the population infected by a major outbreak, $z^{(h,p)}$, against $p$ for different choices of household size $h$ and distribution of $(X_G, X_L)$.
}
\label{fig:z}
\end{figure}

\begin{figure}
\begin{center}
\begin{tabular}{ccc}
 (a) Constant &
  (b) Binomial \\
 \includegraphics[width=7cm]{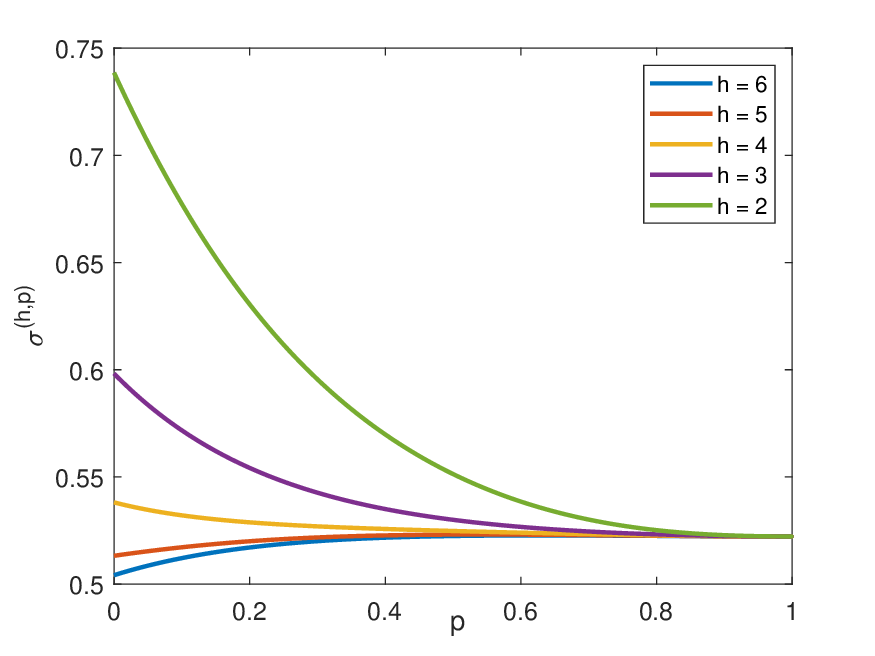} &
 \includegraphics[width=7cm]{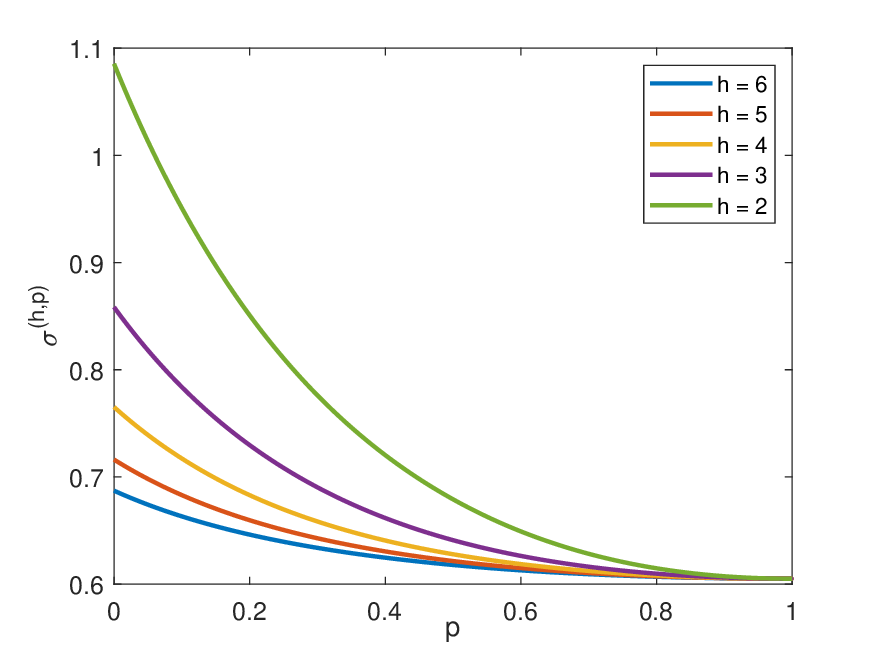} \\
 (c) Poisson &
 (d) Mixed-Poisson, $I \sim {\rm Exp}(1)$\\
 \includegraphics[width=7cm]{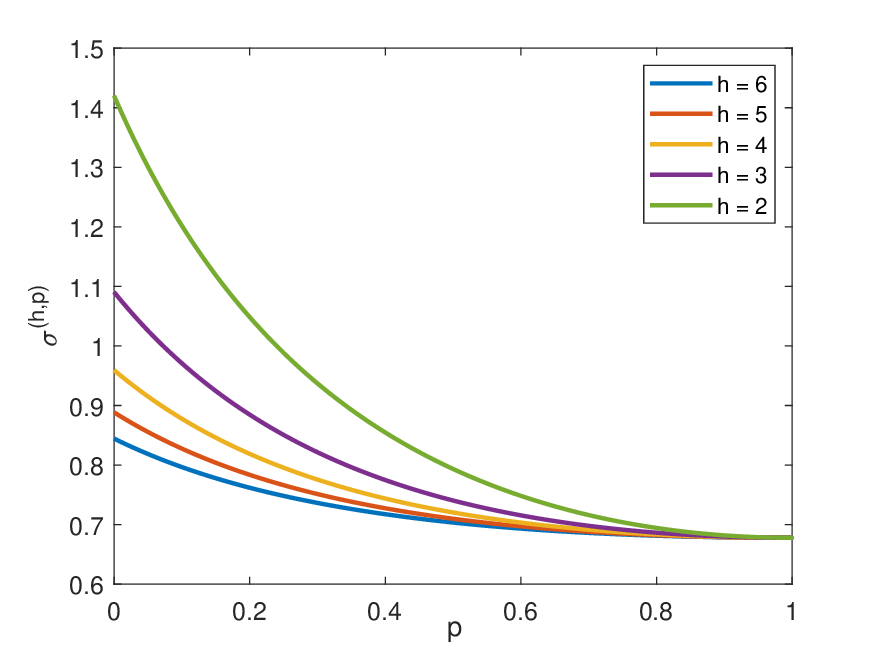} &
 \includegraphics[width=7cm]{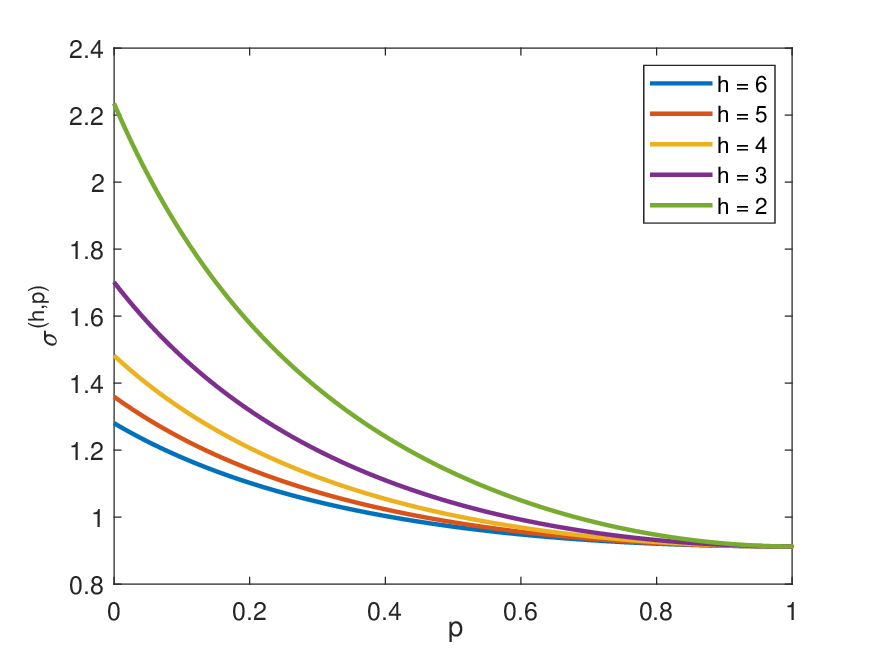}
\end{tabular}
\end{center}
\caption{Graphs of the scaled variance, $\sigma^{(h,p)}$, of the fraction of the population infected by a major outbreak against $p$ for different choices of household size $h$ and distribution of $(X_G, X_L)$.
}
\label{fig:sigma}
\end{figure}

In this section, we illustrate numerically the dependence of $\pi^{(h,p)}$, $z^{(h,p)}$ and $\sigma^{(h,p)}$ on $h$, $p$ and the distribution of $(X_G, X_L)$.  (Recall that  $h$ is the household size, $p$ is the probability that a local contact is replaced by a global contact,
$\pi^{(h,p)}$ is the asymptotic probability of a major outbreak, given one initial infective, and $z^{(h,p)}$ and $\sigma^{(h,p)}$ are the asymptotic mean and scaled standard deviation of the fraction of the population infected by a major outbreak.)   Unless specified otherwise, the naming of the distributions follows exactly that used in Table~\ref{table:hequals2}.  Figures~\ref{fig:z} and~\ref{fig:sigma} show the dependence of $z^{(h,p)}$ and $\sigma^{(h,p)}$ on $h$ and $p$ when $(X_G,X_L)$ is (a) Constant, (b) Binomial, (c) Poisson and (d) Mixed-Poisson with $I \sim {\rm Exp}(1)$.  Note that in both the Poisson and Mixed-Poisson cases, $z^{(h,p)}$ is increasing in both $h$ and $p$, as predicted by Theorem~\ref{thm:finalsize} since for both of these distributions $\log(f_{X_L}(s))$ is convex.  The same holds for this Binomial case, even though then $\log(f_{X_L}(s))$ is not convex, so the condition that $\log(f_{X_L}(s))$ is convex is not necessary for the conclusions of Theorem~\ref{thm:finalsize} to hold.  Observe that in this Constant case, where $(X_G, X_L) \equiv (1,1)$, $z^{(h,p)}$ is decreasing with $p$ when $h=3,4,5,6$, while $z^{(2,p)}$ first increases and then decreases with $p$, and $z^{(2,0)}=z^{(2,1)}$.  The final observation has a simple explanation.  When $p=0$, an infected individual necessarily contacts their housemate, so the epidemic can be viewed as a homogeneously mixing one of fully infected households in which each infected household makes precisely two global contacts.  When $p=1$, the epidemic is homogeneously mixing with each individual making two (global) contacts. Hence, $z^{(2,0)}=z^{(2,1)}$.  For $h=3,4,5,6$, $z^{(h,p)}$ is decreasing with $h$ but the comparison with $h=2$ depends on the value of $p$.

Turning to the scaled standard deviation, note that in the Poisson and Mixed-Poisson cases, $\sigma^{(h,p)}$ is decreasing in both $h$ and $p$.  The same observation holds for all of the cases we have considered in which $\log(f_{X_L}(s))$ is convex.  A possible intuitive explanation is that increasing $h$ and increasing $p$ both have the effect of making the epidemic more homogeneous.  The observation also holds for this Binomial case but, as we illustrate below, it and the above observation concerning $z^{(h,p)}$, do not hold generally when $X_G$ and $X_L$ follow independent Binomial distributions.  In this Constant case, $\sigma^{(h,p)}$ is decreasing with $h$, however it is decreasing with $p$ for $h=2,3,4$ and increasing with $p$ for $h=5,6$.  Note that for the distributions considered, $z^{(h,p)}$ decreases and $\sigma^{(h,p)}$ increases as the variances of $X_L$ and $X_G$ increase.

\begin{figure}
\begin{center}
\includegraphics[width=7cm]{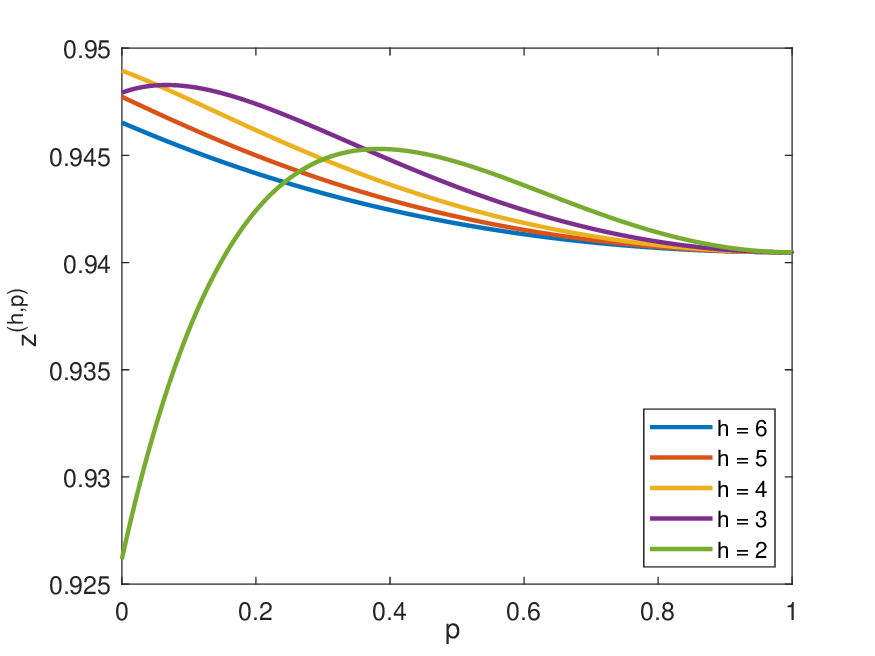}  \includegraphics[width=7cm]{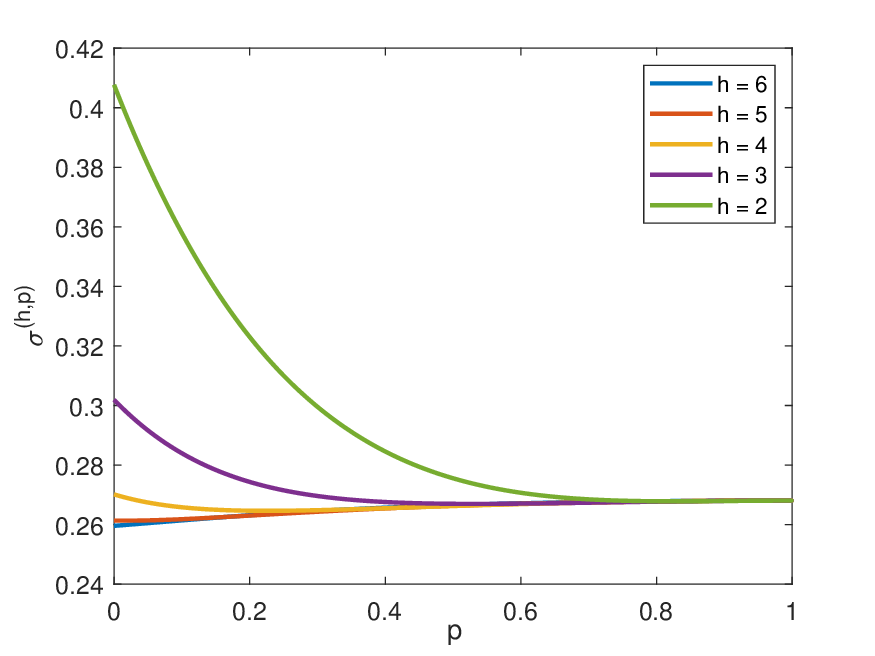}
\end{center}
\caption{Graphs of $z^{(h,p)}$ (left panel) and  $\sigma^{(h,p)}$ (right panel) when $X_G \sim {\rm Bin}(2, \frac{3}{4})$ and $X_L \sim {\rm Bin}(2, \frac{3}{4})$ independently.}
\label{fig:binomialzsigma}
\end{figure}

Figure~\ref{fig:binomialzsigma} shows plots of $z^{(h,p)}$ and $\sigma^{(h,p)}$ when $X_G \sim {\rm Bin}(2, \frac{3}{4})$ and $X_L \sim {\rm Bin}(2, \frac{3}{4})$ independently.  Note that these plots are broadly similar to the corresponding plots in the above Constant case, except here $z^{(h,p)}$ is also non-monotonic with $p$ when $h=3$.

\begin{figure}
\begin{center}
\begin{tabular}{ccc}
  (a) Mixed-Poisson, $I \sim {\rm Gamma}(\frac{1}{2}, \frac{1}{2})$  &
  (b) Binomial \\
 \includegraphics[width=7cm]{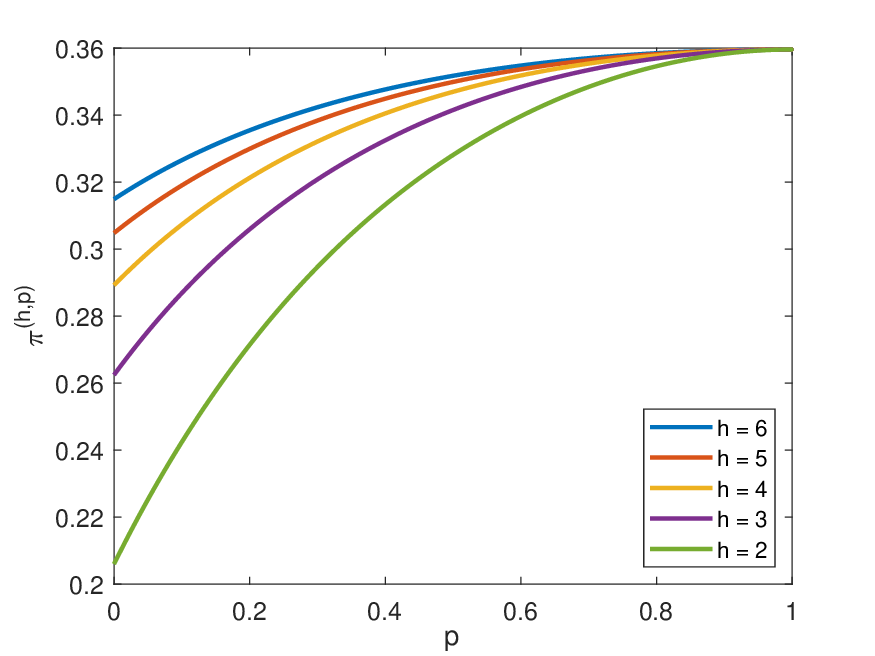} &
 \includegraphics[width=7cm]{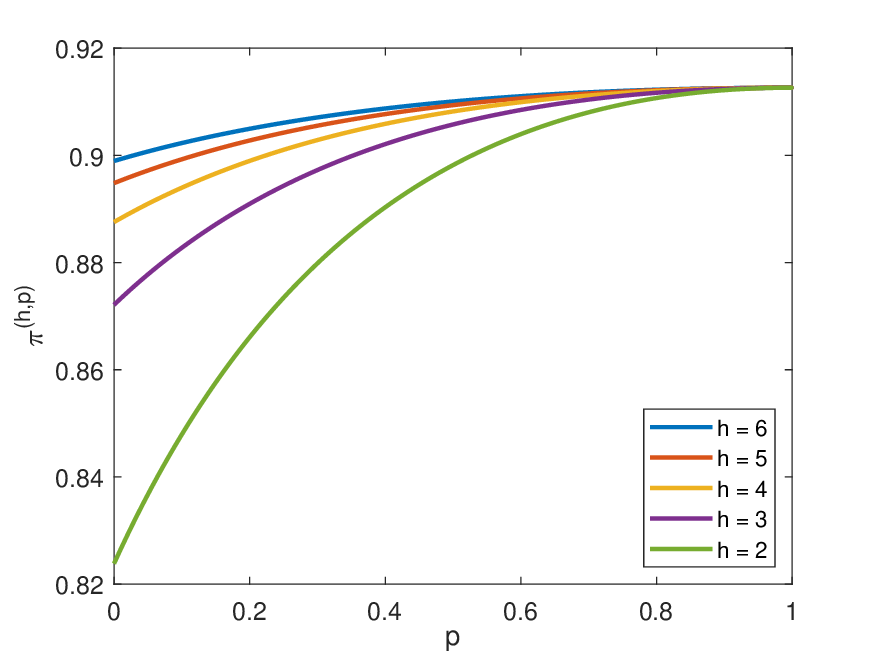} \\
  (c) Poisson &
  (d) Mixed-Poisson, $I \sim {\rm Exp}(1)$\\
 \includegraphics[width=7cm]{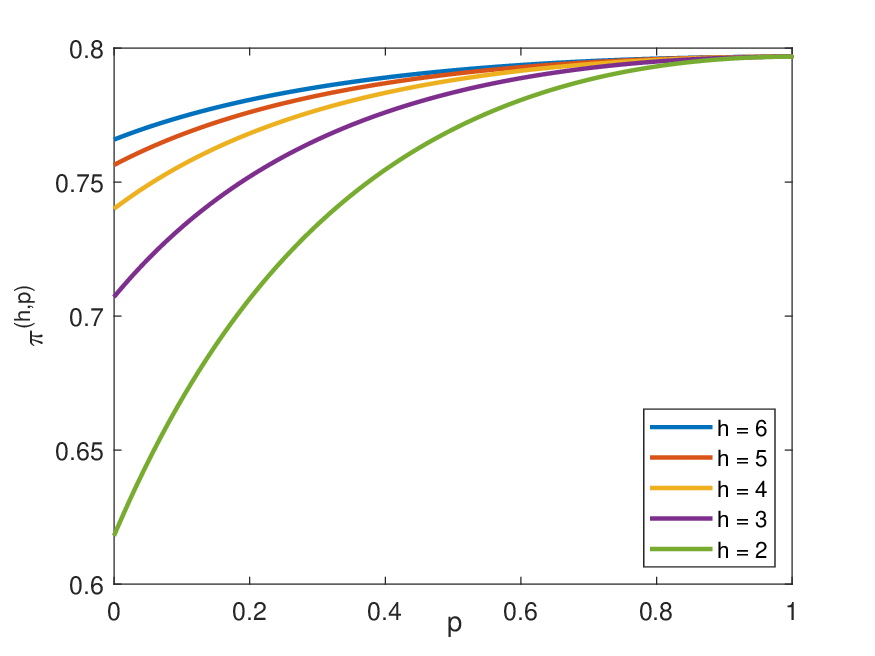} &
 \includegraphics[width=7cm]{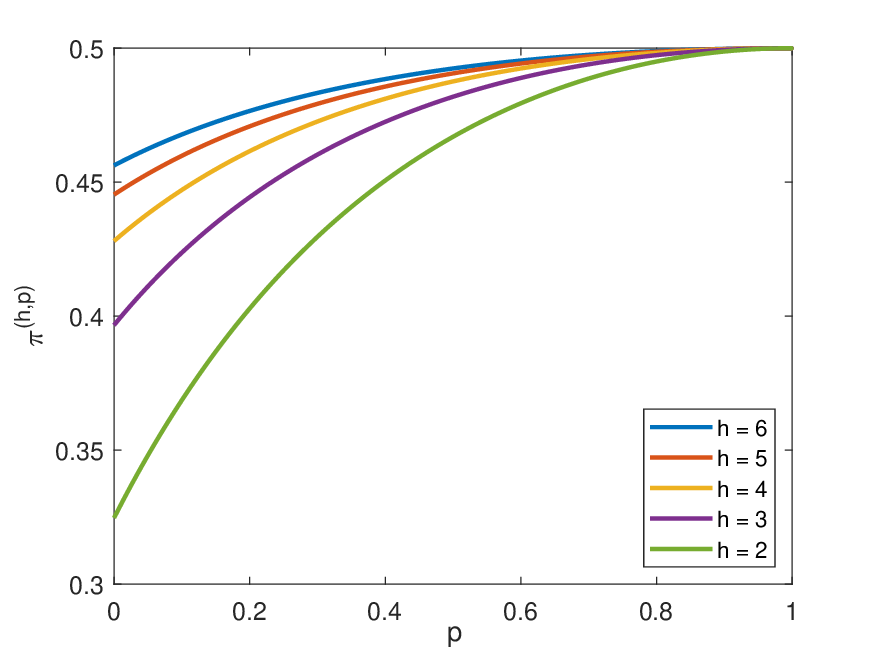}
\end{tabular}
\end{center}
\caption{Graphs of the probability of a major outbreak, $\pi^{(h,p)}$, against $p$ for different choices of household size $h$ and distribution of $(X_G, X_L)$.
}
\label{fig:pmaj}
\end{figure}

Finally, Figure~\ref{fig:pmaj} shows plots of the probability of a major outbreak, $\pi^{(h,p)}$, for various choices of distribution for $(X_G, X_L)$.  Note that in all cases, $\pi^{(h,p)}$ is increasing in both $h$ and $p$, as predicted by Theorem~\ref{thm:prob:major}.  For fixed $(h,p)$, $\pi^{(h,p)}$  decreases as the variances of $X_G$ and $X_L$ increase.  Note that in the Poisson case, $\pi^{(h,p)}=z^{(h,p)}$, while in the other cases in which $\log(f_{X_L}(s))$ is convex, $\pi^{(h,p)}<z^{(h,p)}$ (see also Table~\ref{table:hequals2} when $(h,p)=(2,0)$).  This is usually the case for epidemic models.  However, in the Binomial case, $\pi^{(h,p)}>z^{(h,p)}$.

\section{Central limit theorem proof} \label{sec:clt}

\subsection{Introduction} \label{sec:clt:intro}

In this section we prove Theorem \ref{thm:clt}. We begin in Section \ref{sec:clt:model} by defining  a sequence of $\mathcal{E}_{n,h}  (X_G,X_L)$ epidemics, $\tilde{\mathcal{E}}_n$, indexed by $n$ the number of households. 
In Section \ref{sec:clt:BP}, we give a branching process approximation for the early stages of the epidemic and show that the probability of a minor outbreak (which infects at most $\lfloor \log n \rfloor$ households) converges to $\rho^m$ as $n \to \infty$, where $\rho$ satisfies  \eqref{eq:rho}.
 In Section \ref{sec:clt:Sellke} we define the embedding process which is utilised for the central limit theorem. The embedding process is based on a Sellke construction, see Sellke \cite{Sellke83}, of the epidemic with an extra level of embedding. We define a sequence of epidemics $\mathcal{E}_n$ based on the embedded construction and show that $\tilde{\mathcal{E}}_n$ and $\mathcal{E}_n$ can be coupled to give the same epidemic final size, albeit with potentially different global infectors of individuals. This enables us to focus on the embedded construction in the remainder of the section. In Section \ref{sec:clt:lln}, we prove a law of large numbers result and show that $\bar{Z}_{n,h} \convD Z$ as $n \to \infty$, where the probability mass function of $Z$ satisfies \eqref{eq:probZ}.  In Section \ref{sec:clt:proof} we prove Theorem \ref{thm:clt} by exploiting an upper and lower bound for the proportion infected in the event of a major epidemic and showing that both these bounds have the same limit. A key component in the proof is Theorem \ref{thm:cltbound} whose proof is postponed to Section \ref{sec:clt:proofhat}.
In Section \ref{sec:clt:variance}, we discuss $\sigma^2$ and give two equivalent expressions for $\sigma^2$ in \eqref{eq:sigma1} and \eqref{eq:sigma2}.  The first expression, \eqref{eq:sigma1}, arises naturally in the proof of Theorem \ref{thm:clt}, whilst the second expression,  \eqref{eq:sigma2}, is often simpler to work with in terms of computing $\sigma^2$ numerically.
The proof that the expressions in \eqref{eq:sigma1} and \eqref{eq:sigma2} are equal, and equivalent to that given by~\eqref{eq:sigma} in Section~\ref{sec:SSHEmain:result} are deferred to Appendix \ref{app:sigma2}. Finally, in Section \ref{sec:clt:without} we discuss the minor modifications to the central limit theorem for the case where the contacts $(X_G,X_L)$ are sampled {\it without} replacement from the population and household, respectively.

\subsection{Model description}  \label{sec:clt:model}

 For $i=1,2, \ldots$ and $j=1,2,\ldots,h$, let $\mathbf{X}_{ij}$ be i.i.d.~copies of $\mathbf{X} = (X_G, X_L)$ with $\mathbf{X}_{ij}$ determining the number of global and local infectious contacts made by the $j^{th}$ individual in household $i$.
We construct the epidemic $\tilde{\mathcal{E}}_n$ using $\{ \mathbf{X}_{ij} = (X_{G,(i,j)}, X_{L,(i,j)}); i=1,2,\ldots, n, j=1,2,\ldots,h\}$ as follows.  We assign to each individual a list of household contacts $\mathbf{H}_{ij} = (H_{ij1}, H_{ij2}, \ldots)$, where $H_{ijk}$ is the individual within the household contacted by the $k^{th}$ household infectious contact made by individual $j$ in household $i$. (Note that the $\{H_{ijk}\}s$ are independent and uniformly distributed on $\{1,2,\ldots h\} \backslash j$.)  The individual $(i,j)$ makes a household infectious contact with individual $(i,l)$ if $l \in \{ H_{ij1}, H_{ij2}, \ldots, H_{ijX_{L,(i,j)}}\}$.
 In addition, for each $n$, we let $U_1^n, U_2^n, \ldots$ be i.i.d.~copies of $U^n$, where
\[ \P (U^n = (i,j)) = \frac{1}{nh} \qquad  ( i=1,2,\ldots, n; j=1,2,\ldots,h). \] Therefore, $U^n$ can be used to choose an individual uniformly at random from the population underlying the epidemic $\tilde{\mathcal{E}}_n$.

The epidemic $\tilde{\mathcal{E}}_n$ starts with $m_n$ initial infectives, and we assume that there exists $m \geq 1$ such that $m_n=m$ for all sufficiently large $n$. The $m_n$ initial infectives are given by the first $m_n$ unique $U^n$.
For $n=1,2,\ldots$ and $k=1,2,\ldots$, let $\mathcal{I}_k^n =  \cup_{i=1}^k \{U_i^n\}$. Then $\mathcal{I}_{K_n}^n$ denotes the set of initial infectives where $K_n$ satisfies
\[ K_n = \min \left\{ k: \left|\mathcal{I}_k^n  \right| = m_n\right\}. \]
The epidemic is then constructed by  considering infectives one at a time. Suppose that prior to considering individual $(i_0,j_0)$, say, there has been a total of $M$ global infectious contacts.  The local infectious contacts made by individual $(i_0,j_0)$ are governed by $X_{L,(i_0,j_0)}$ and  $\mathbf{H}_{i_0j_0}$. The global infectious contacts made by individual $(i_0,j_0)$ are with, if $X_{G,(i_0,j_0)} >0$, individuals $U^n_{M+1}, U^n_{M+2},  \ldots, U^n_{M+X_{G,(i_0,j_0)}}$. The process continues until there are no more infectives in the population.

\subsection{Branching process approximation} \label{sec:clt:BP}

For the epidemic $\tilde{\mathcal{E}}_n$ we have defined a major epidemic as one that infects at least $k_n = \lfloor \log n \rfloor$ households. Therefore we define a minor epidemic as one that infects fewer than $\lfloor \log n\rfloor$ households, that is, if $V_{n,h} < \lfloor \log n \rfloor$ and in this section we show that 
\begin{align} \label{eq:V:BP}
\P (V_{n,h} < \lfloor \log n\rfloor) \to \rho^m \qquad \mbox{as } n \to \infty,
\end{align}
where $\rho$ satisfies \eqref{eq:rho}.

In order to prove \eqref{eq:V:BP}, we couple the sequence of epidemics  $\tilde{\mathcal{E}}_n$ to a Galton-Watson branching process $\mathcal{B}$. Specifically, the branching process $\mathcal{B}$ has $m$ ancestors and the number of offspring from individuals are i.i.d.~copies of $C$, defined just before~\eqref{eq:sus:pgf} in Section~\ref{sec:SSHEmain:result}. Hence, $\rho$ denotes the extinction probability of the branching process $\mathcal{B}$. Let $V$ denote the total size, including initial ancestors, of the branching process $\mathcal{B}$.

\begin{lem} \label{lem:BP}
For any $k =1,2,\ldots$,
\[ \P (V_{n,h} \leq k) \to \P (V \leq k) \qquad \mbox{as } n \to \infty. \]
\end{lem}
\begin{proof} We prove the lemma by constructing $\tilde{\mathcal{E}}_n$ and $\mathcal{B}$ on a common probability space. For $i=1,2, \ldots$ and $j=1,2,\ldots, h$, let $\bar{\mathbf{X}}_{ij}$ be i.i.d.~copies of $\mathbf{X}$ and let $\bar{\mathbf{H}}_{ij}$ be independent with $\bar{\mathbf{H}}_{ij} \eqd \mathbf{H}_{1j}$.
For $i=1,2,\ldots$, let $C_i$ denote the number of global contacts emanating from the $i^{th}$ household epidemic constructed using $\{ \bar{\mathbf{X}}_{ij}, \bar{\mathbf{H}}_{ij}; j=1,2,\ldots,h\}$, where the individual $(i,1)$ is the initial infective in the household. Let $C_i$ denote the number of offspring of the $i^{th}$ individual in the branching process $\mathcal{B}$ with $C_1, C_2, \ldots, C_m$ denoting the offspring of the $m$ ancestors.

Let $\tilde{U}_1^n, \tilde{U}_2^n, \ldots$ be i.i.d.~copies of $\tilde{U}^n$, where $\tilde{U}^n$ is a discrete uniform distribution on $\{1,2,\ldots,n\}$. We construct a realisation of $\tilde{\mathcal{E}}_n$ by assigning the $i^{th}$ global contact in $\tilde{\mathcal{E}}_n$ to household $\tilde{U}_i^n$. Given that household $\tilde{U}_i^n$ has not previously been infected we assign infectious histories $\{ \bar{\mathbf{X}}_{ij}, \bar{\mathbf{H}}_{ij}; j=1,2,\ldots,h\}$ to the individuals in household $\tilde{U}_i^n$ and assume that the individual contacted globally is individual $(i,1)$. Therefore the number of global contacts emanating from the first household epidemic in household $\tilde{U}_i^n$  is $C_i$.

Let $M_n = \min \left\{ k >1: \tilde{U}_k^n \in \left\{ \tilde{U}_1^n,  \tilde{U}_2^n, \ldots,  \tilde{U}_{k-1}^n \right\} \right\}$, the number of global contacts that occur until the first attempted infection of a previously infected household. This is the well known {\it Birthday Problem}, see for example Ball and Donnelly \cite{BD95}, and
\begin{equation}
\label{equ:birthday}
\P(M_n \le k) \le \frac{k(k-1)}{2n} \qquad (k=2,3,\dots,n).
\end{equation}
Therefore, for any $k=1,2,\ldots$,
\begin{align}
\P (V_{n,h} \leq k) & = \P (V_{n,h} \leq k | M_n > k) \P (M_n >k) +  \P (V_{n,h} \leq k | M_n \leq k) \P (M_n \leq k) \nonumber \\
& =  \P (V \leq k | M_n > k) \P (M_n >k) +  \P (V_{n,h} \leq k | M_n \leq k) \P (M_n \leq k) \nonumber \\
& \to \P (V \leq k) \qquad \mbox{as } n \to \infty,
\end{align}
as required.
\end{proof}

Given that \eqref{equ:birthday} implies  $\P (M_n >\lfloor \log n \rfloor) \to 1$ as $n \to \infty$,  it is straightforward to show that
\[ \left| \P (V_{n,h} \leq \log n \rfloor) - \P (V \leq \log n \rfloor)  \right| \to 0 \qquad \mbox{as } n \to \infty. \]
Since $\P (k < V< \infty) \to 0$ as $k \to \infty$ and $\P (V < \infty) = \rho^m$, it follows by the triangle inequality that
\begin{align}
 \left| \P (V_{n,h} \leq \log n \rfloor) - \rho^m  \right|  & \leq
 \left| \P (V_{n,h} \leq \log n \rfloor) - \P (V \leq \lfloor \log n \rfloor)  \right|  + \P ( \lfloor \log n \rfloor < V < \infty) \nonumber \\
 & \to 0 \qquad \qquad \mbox{as } n \to \infty. \nonumber
\end{align}

 \subsection{Embedding}  \label{sec:clt:Sellke}

In order to obtain a central limit theorem for the final size, we use an embedding argument similar to \cite{ST85}, \cite{BMST} and Ball and Neal \cite{BN03}, utilising a Sellke (\cite{Sellke83}) construction of the epidemic. This involves taking an alternative approach to modelling global infection but we  show that the final size of the epidemic is unchanged. Specifically, we assume that any given individual encounters global infections at the points of a homogeneous unit rate Poisson point process as the amount of global infectious pressure they are exposed to increases. In \cite{BMST} and \cite{BN03}, an infectious individual with infectious period $I$ contributes $\lambda_G I/N$ units of global infectious to each individual in the population with the number of new global infectious encounters arising following a Poisson distribution with mean $\lambda_G I$. In our setting, each infective makes a given number of global contacts distributed according to $X_G$. This means that we cannot directly apply the embedding arguments used in the earlier referenced works but require an additional layer of embedding which links the total number of global contacts in the epidemic process to the independent Poisson point processes of global contacts attached to individuals.

Before defining a sequence of embedded epidemics, $\mathcal{E}_n$, indexed by $n$ the number of households and showing that $\mathcal{E}_n$ and $\tilde{\mathcal{E}}_n$ give the same final size, we require some additional notation.  This includes the formal definition of a susceptibility set whose pgf plays a key role in obtaining, $z$, the mean final proportion infected in a major outbreak given by \eqref{eq:z}.

For $i=1,2,\ldots$ and $j,l =1,2,\ldots, h$, let $(i,j) \leadsto (i,l)$ denote that there is a path of household infection from individual $(i,j)$ to individual $(i,l)$ with the convention that $(i,j) \leadsto (i,j)$. Note that $(i,j) \leadsto (i,l)$ is determined by $\{ (X_{L,(i,k)} ,\mathbf{H}_{ik}) ; k=1,2,\ldots,h \}. $
For $i=1,2,\ldots$ and $j=1,2,\ldots, h$, let $\mathcal{S}^{ij}$ denote the susceptibility set of individual $(i,j)$ which is defined to be
\begin{align} 
\mathcal{S}^{ij} = \{ l \in \{1,2,\ldots,h \} : (i, l) \leadsto (i,j) \}. \nonumber
\end{align}
That is, $\mathcal{S}^{ij}$ is the set of individuals whom if infected by a global infection will infect individual $(i,j)$, if susceptible, via a chain of local infections within the household. Let $S_{ij} = |\mathcal{S}^{ij} |$ denote the size of the susceptibility set of individual $(i,j)$. Note that for all $(i,j)$, $S_{ij} \eqd S_{11}$ and for $k \neq i$, $S_{ij}$ and $S_{kl}$ are independent with the pgf of $S_{11}$ given by $f_S(s)$,  cf.~\eqref{eq:sus:pgf}.

Finally, before introducing the embedded epidemic process we attach to each individual $(i,j)$ an independent, homogeneous Poisson point process, $\eta_{ij}$, with rate 1. For $t \geq 0$, let $\zeta_{ij} (t)$ denote the number of points of $\eta_{ij}$ in  $[0,t]$. Thus $\zeta_{ij} (t) \sim {\rm Po} (t)$.

Suppose that global contacts occur with an individual  at the points of a homogeneous Poisson point process with rate 1. Specifically, we assume that individual $(i,j)$ receives global contacts at the points of $\eta_{ij}$ as the individual is exposed to increasing amounts of global infection. We assume that when an individual is infected globally the local household epidemic from that individual occurs instantaneously. Let $\chi_{ij} (t) = 1 - \prod_{l \in \mathcal{S}^{ij}} 1_{\{ \zeta_{il} (t) =0\}}$.
Then $\chi_{ij} (t)$ is an indicator random variable for whether or not individual $(i,j)$ is infected when all members of the population are exposed to $t$ units of global infectious pressure, since an individual is infected once somebody in their susceptibility set receives a global infectious contact.

For $i=1,2,\ldots$ and $t \geq 0$, let $(R_i (t), G_i (t), Y_i (t))$
be a trivariate random variable determining the state of household when each individual is exposed to $t$ units of global infection.
 Let $R_i (t) = \sum_{j=1}^h \chi_{ij} (t)$ denote the number of individuals infected in the household, let $G_i (t) = \sum_{j=1}^h X_{G,(i,j)} \chi_{ij} (t)$ denote the number of global contacts made by those infected in the household and let $Y_i (t) = \sum_{j=1}^h \zeta_{ij} (t) [ =  \sum_{j=1}^h \chi_{ij} (t) \zeta_{ij} (t) ]$ denote the number of global contacts made into the household. By construction the $\{(R_i (t),G_i (t))\}$s are i.i.d.~copies of $(R (t),G (t))$, defined in Section \ref{sec:main:result}.

For $t \geq 0$, let $\nu_R (t) = \E [R_1 (t)]/h = \E[ \chi_{11} (t)] = 1- f_S (e^{-t})$, cf.~\eqref{eq:nuR}. Since, for all $t \geq 0$, $X_{G,(1,1)}$ and $\chi_{11} (t)$ are independent, we have that
\begin{align} 
\nu_G (t) = \frac{1}{h} \E[G_1 (t)] = \mu_G  [ 1- f_S (e^{-t})]. \nonumber
\end{align}
Finally, $\nu_Y (t) = \E[Y(t)]/h = t$.

We are now in position to describe the construction of the embedded epidemic process $\mathcal{E}_n$ and obtain an expression for the proportion, $\bar{Z}_{n,h}$, of the population infected.

The embedded epidemic process considers each individual, and hence, household being exposed to infection at a constant rate. If each member of the population is exposed to $t$ units of global infection, the total number of global infectious contacts is random and distributed according to ${\rm Po} (n h t)$, the number of points in $[0,t]$ of the Poisson point process $\eta^n$, where $\eta^n$ is defined to be the superposition of the Poisson processes $\{ \eta_{ij}; i=1,2,\ldots, n, j=1,2,\ldots,h\}$.
 To study the original epidemic process using the embedded epidemic process, we reverse this procedure and for a given $x \in \mathbb{R}^+$, we find the random time $S_n (x)$ such that the number of global contacts in the population on the interval $[0,S_n (x)]$ is equal to $\lfloor x n h \rfloor$. More specifically, for $n=1,2,\ldots$ and $x \geq 0$, let 
\begin{align}
\label{eq:Snx}
S_n (x) = \min \left\{ t \geq 0 : \sum_{i=1}^n Y_i (t) = \lfloor x n h \rfloor  \right\}.
\end{align}

  Let $T_0^n$ denote the number of global infections required to generate $m_n$ infectives to initiate the epidemic, and remember that $m_n =m$ for all sufficiently large $n$. Therefore $T_0^n \convp m$ as $n \to \infty$. Let $\bar{T}_0^n = T_0^n/(nh)$. Then
\begin{align}
S_n (\bar{T}_0^n) = \min \left\{ t \geq 0: \sum_{i=1}^n Y_i (t) = T_0^n \right\} \nonumber
\end{align}
is the initial amount of global infection in the epidemic process $\mathcal{E}_n$ to generate $m_n$ infectives ($T_0^n$ global infectious contacts). We say that the set of individuals whose susceptibility set contains an initial infective form generation 0 of $\mathcal{E}_n$. (Therefore generation 0 of $\mathcal{E}_n$ is obtained by running the local epidemics from the initial infectives.) Generation 0 will generate $\sum_{i=1}^n G_i (S_n (\bar{T}_0^n))$ global infectious contacts. Thus
\begin{align}
T_1^n (= n h \bar{T}_1^n) = T_0^n + \sum_{i=1}^n G_i (S_n (\bar{T}_0^n)), \nonumber
\end{align}
is the number of global infections,  including those required for the initial infectives,  after the global infections emanating from generation 0 have been considered.
Following \cite{BMST}, Section 4.2.2, we can define $T_0^n, T_1^n, \ldots$, with $\bar{T}_k^n = T_k^n /(nh)$, to satisify, for $k=0,1,\ldots$,
 \begin{align}
T_{k+1}^n (= n h \bar{T}_{k+1}^n) = T_0^n + \sum_{i=1}^n G_i (S_n (\bar{T}_k^n)). \nonumber
\end{align}
For $k=1,2,\ldots$, we say an individual belongs to the $k^{th}$ generation of infectives if the first time a member of their susceptibility set is infected globally is by a member of generation $k-1$. Using the embedding process an individual $(i,j)$ belongs to generation $k$ if
\[ \chi_{ij} (S_n (\bar{T}_{k-1}^n)) = 0 \qquad \mbox{and} \qquad \chi_{ij} (S_n (\bar{T}_k^n) )= 1, \]
and $T_{k+1}^n$ is the total number of global infections, including  those required for the initial infectives, from the first $k$ generations of infectives.
 The process continues until there are no additional global infections created in a generation. That is, $T_{k+1}^n = T_k^n$, and consequently we can define $T_\infty^n = n h \bar{T}_\infty^n$ to satisfy
\begin{align}
 \label{eq:Tnx}
\bar{T}_\infty^n = \min \left\{ x \geq 0 :  T_0^n + \sum_{i=1}^n G_i (S_n (x)) = \lfloor x nh \rfloor \left( =\sum_{i=1}^n Y_i (S_n (x)) \right) \right\}.
\end{align}
Hence,
 \begin{align}
\bar{T}_\infty^n = \bar{T}_0^n + \frac{1}{nh} \sum_{i=1}^n G_i (S_n (\bar{T}_\infty^n)) \left( = \frac{1}{nh} \sum_{i=1}^n Y_i (S_n (\bar{T}_\infty^n))\right). \nonumber
\end{align}

Therefore, $\bar{Z}_{n,h}$, the proportion of the population infected by the epidemic $\mathcal{E}_n$, satisfies
 \begin{align}
\bar{Z}_{n,h} = \frac{1}{nh} \sum_{i=1}^n R_i (S_n (\bar{T}_\infty^n)) = \frac{1}{nh} \sum_{i=1}^n \sum_{j=1}^h \chi_{ij} (S_n (\bar{T}_\infty^n)). \nonumber
\end{align}

We show how the epidemic processes $\mathcal{E}_n$ and $\tilde{\mathcal{E}}_n$ can be coupled to give the same final size.
We construct $\mathcal{E}_n$ using $\{ \mathbf{X}_{ij} = (X_{G,(i,j)}, X_{L,(i,j)}), \mathbf{H}_{ij}, \eta_{ij}; i=1,2,\ldots, n, j=1,2,\ldots,h\}$. To construct $\tilde{\mathcal{E}}_n$ from  $\mathcal{E}_n$, we use $\{  \mathbf{X}_{ij} = (X_{G,(i,j)}, X_{L,(i,j)}), \mathbf{H}_{ij} ; i=1,2,\ldots, n, j=1,2,\ldots,h\}$ so local epidemics are unchanged and the number of global contacts made by a given individual are the same in both processes. Using $\{ \eta_{ij}; i=1,2,\ldots, n, j=1,2,\ldots,h\}$, we construct $U^n_1, U^n_2, \ldots$. 
For $k=1,2,\ldots$, we set $U_k^n = (i^\prime, j^\prime)$ if the $k^{th}$ point of $\eta^n$ comes from $\eta_{i^\prime j^\prime}$. This construction means that the initial $m_n$ infectives in $\tilde{\mathcal{E}}_n$ are $\mathcal{I}_{K_n}^n$ and that the individual contacted by the $k^{th}$ global contact is the same in both epidemics although the assignment of the infector might be different. Consequently, those individuals whose susceptibility sets have been globally infected, and thus are guaranteed to be infected, by the first $t$ global infections in $\tilde{\mathcal{E}}_n$, is precisely the set of individuals for whom $\chi_{ij} (S_n(t/nh))=1$ $(i=1,2,\ldots, n; j=1,2,\ldots,h)$. Therefore, $T_\infty^n$ is the total number of global infectious contacts in both  $\mathcal{E}_n$ and $\tilde{\mathcal{E}}_n$, with $\bar{Z}_{n,h}$ denoting the proportion of individuals infected.

\subsection{Law of large numbers} \label{sec:clt:lln}

In this section we prove that the proportion of the population infected, $\bar{Z}_{n,h}$, converges to a random variable $Z$ whose probability mass function is defined in \eqref{eq:probZ}.

\begin{lem} \label{lem:lln} Suppose that there exists $m \in \mathbb{N}$ such that $m_n = m$ for all sufficiently large $n$. For $R_\ast >1$, there exists $0 < \tau < \infty$ which solves $\tau = \nu_G (\tau)$ with
 \begin{align} \label{eq:lln}
 \min \left\{ |S_n (\bar{T}_\infty^n)|, |S_n (\bar{T}_\infty^n) - \tau|\right\} \convas 0 \qquad \mbox{as } n \to \infty.
  \end{align}
\end{lem}
\begin{proof} Firstly, $m_n =m$ for all sufficiently large $n$,  implies that $\bar{T}_0^n \convas 0$ as $n \to \infty$.  By the strong law of large numbers, $(nh)^{-1} \sum_{i=1}^n G_i (t) \convas \nu_G (t)$ as $n \to \infty$, for all $t \ge 0$.  Also, $\nu_G (\infty) = \mu_G < \infty$.  A similar, but simpler, argument to the proof of \cite{BN03}, Lemma 3.8, yields
\begin{align}
\label{eq:Gas:ast}
\sup_{t \geq 0} \left| \frac{1}{nh} \sum_{i=1}^n G_i (t) - \nu_G (t) \right| \convas 0 \qquad \mbox{as } n \to \infty.
\end{align}

By a similar argument, for any $T >0$, we have that
\begin{align}
\label{eq:Pas}
\sup_{0 \leq t \leq 2T} \left| \frac{1}{nh} \sum_{i=1}^n Y_i (t) - t\right| \convas 0 \qquad \mbox{as } n \to \infty .
\end{align}
For any $x \geq 0$, using \eqref{eq:Snx}, we have that
\begin{align} 
|S_n(x) - x|  &= \left| S_n (x) - \frac{1}{nh} \sum_{i=1}^n Y_i (S_n (x)) + \frac{\lfloor x nh \rfloor}{nh} - x \right| \nonumber \\
& \leq  \left| S_n (x) - \frac{1}{nh} \sum_{i=1}^n Y_i (S_n (x)) \right| + \left| \frac{\lfloor x nh \rfloor - x nh }{nh} \right|.  \nonumber
\end{align}

Since $S_n (x)$ is increasing in $x$, it follows that for $S_n (T) \leq 2 T$,
\begin{align} \label{eq:Snx:as2}
0 \leq \sup_{0 \leq x \leq T} |S_n(x) - x| & \leq \sup_{0 \leq x \leq T}  \left( \left| S_n (x) - \frac{1}{nh} \sum_{i=1}^n Y_i (S_n (x)) \right| + \left| \frac{\lfloor x nh \rfloor - x nh }{nh} \right|\right)  \nonumber   \\
& \leq \sup_{0 \leq t\leq 2T}  \left| \frac{1}{nh} \sum_{i=1}^n Y_i (t) - t\right| + \frac{1}{nh}.
\end{align}
Given $(nh)^{-1} \sum_{i=1}^n Y_i (2 T) \geq T$ implies that $S_n (T) \leq 2 T$ and $(nh)^{-1} \sum_{i=1}^n Y_i (2 T) \convas 2T$ as $n \to \infty$, it follows from \eqref{eq:Snx:as2} and \eqref{eq:Pas}  that
\begin{align} 
 \sup_{0 \leq x \leq T} |S_n(x) - x|  \convas 0 \qquad \mbox{as } n \to \infty. \nonumber
\end{align}

Let $\mathcal{K}= \{ t \in [0,\infty]: t =  \nu_G (t) \}$. Since $\nu_G (\cdot)$ is a strictly concave function of $t$, it follows that $\mathcal{K}= \{ 0, \tau\}$ for $R_\ast >1$. Also $\nu_G^\prime (\tau) \neq 1$ for all $\tau \in \mathcal{K}$. Let $(\Omega, \mathcal{F}, \mathcal{P})$ denote the probability space on which the random vectors $(R_1 (t), G_1 (t), Y_1(t)), (R_2 (t), G_2 (t), Y_2(t)), \ldots$ are defined. Fix $T > \tau$ and let
\begin{align*}
F_1 &= \left\{ \omega \in \Omega:  \sup_{t \geq 0} \left| \frac{1}{nh} \sum_{i=1}^n G_i (t,\omega) - \nu_G (t) \right| \to 0 \mbox{ as } n \to \infty \right\} \\
F_2& = \left\{ \omega \in \Omega:  \sup_{0 \leq x \leq T} \left| \frac{1}{nh} \sum_{i=1}^n S_n (x,\omega) - x \right| \to 0 \mbox{ as } n \to \infty \right\} 
\end{align*}
and
\[ F_3 = \left\{ \omega \in \Omega:  \lim_{n \to \infty} \bar{T}_0^n (\omega)=0 \right\}. \]
Then
\[ \min \left\{ |S_n (\bar{T}_\infty^n, \omega) - \tau |: \tau \in \mathcal{K}  \right\} \to 0 \qquad \mbox{as } n \to \infty, \]
for all $\omega \in F_1 \cap F_2 \cap F_3$. The lemma follows since $\P (F_1 \cap F_2 \cap F_3) =1$.
\end{proof}

A corollary of Lemma \ref{lem:lln} concerns the proportion infected in the epidemic.  For $R_*>1$, let $z=\tau/\mu_G$.  Note that $z=\nu_R (\mu_G z)$, so $z$ coincides with the definition at \eqref{eq:z}.
\begin{cor} \label{cor:lln} Suppose that there exists $m \in \mathbb{N}$ such that $m_n=m$  for all sufficiently large $n$. For $R_\ast >1$, we have that
 \begin{align} \label{eq:prop_lln}
 \min \left\{ |\bar{Z}_{n,h}|, |\bar{Z}_{n,h} -z|\right\} \convas 0 \qquad \mbox{as } n \to \infty.
  \end{align}
\end{cor}
\begin{proof} An identical line of argument to the derivation of \eqref{eq:Gas:ast} gives
\begin{align}
\label{eq:Ras:ast}
\sup_{t \geq 0} \left| \frac{1}{nh} \sum_{i=1}^n R_i (t) - \nu_R (t) \right| \convas 0 \qquad \mbox{as } n \to \infty,
\end{align}
Then using Lemma \ref{lem:lln}, \eqref{eq:lln} and \eqref{eq:Ras:ast} it is straightforward to prove \eqref{eq:prop_lln} along similar lines to the proof of Lemma  \ref{lem:lln}.
\end{proof}

The final step to prove that  $\bar{Z}_{n,h} \convD Z$, where $Z$ has probability mass function given by \eqref{eq:probZ}, is to show that for any $0 < \epsilon <z$, $\P (\bar{Z}_{n,h} < \epsilon) \rightarrow \rho^m$ as $n \to \infty$. Let $\bar{V}_{n,h} = V_{n,h}/n$. By construction we have that $\bar{V}_{n,h}/h \leq \bar{Z}_{n,h} \leq \bar{V}_{n,h}$ and therefore it suffices to show that there exists $\epsilon^\prime >0$,
\begin{align} \label{eq:barV}
\P (\bar{V}_{n,h} \leq \epsilon^\prime) \rightarrow \rho^m \qquad \mbox{as } n \to \infty.
\end{align}
It is straightforward using  a lower bound branching process, cf.~Whittle \cite{Whittle55}, \cite{BMST}, Ball and Lyne \cite{BL01}, to show that \eqref{eq:barV} holds by following a similar line of argument to the proof of Ball and Neal \cite{BN24}, Theorem 3.2. An outline of the argument is as follows. We couple $\tilde{\mathcal{E}}_n$ and $\mathcal{B}$ until $k_n =\lfloor \log n \rfloor$ households have been infected. The first $k_n$ household epidemics will generate approximately $R_\ast k_n$ global infections. More precisely, we can show that  for any $0 < \delta < R_\ast -1$, the first $\lfloor \log n \rfloor$ household epidemics create at least a further $ \lfloor  \delta \log n \rfloor$ local epidemics in distinct households with probability tending to 1 as $n \to \infty$. For any $0 < \epsilon^\prime < [R_\ast -1]/R_\ast$, we consider a super-critical lower bound branching process approximation to the epidemic starting from $ \lfloor  \delta \log n \rfloor$ individuals where each birth in the branching process is aborted independently with probability $\epsilon^\prime$. Since the lower bound branching process is super-critical and $\E [C^2] < \infty$, we have that the extinction probability, $\rho (\epsilon^\prime)$, from a single ancestor is bounded away from 1 by Ball and Neal \cite{BN17}, Lemma A3, with $\rho (\epsilon^\prime)^{ \lfloor  \delta \log n \rfloor} \to 0$ as $n \to \infty$.  Whilst the proportion of households infected is less than $\epsilon^\prime$, the probability that a global contact is with a previously infected household is at most $\epsilon^\prime$. We can then use the lower bound branching process to show that
\[ \P (\bar{V}_{n,h} \leq \epsilon^\prime | V_{n,h} > \lfloor \log n \rfloor)  \to 0 \qquad \mbox{as } n \to \infty, \]
which combined with \eqref{eq:V:BP} yields \eqref{eq:barV}.

\subsection{Proof of Theorem \ref{thm:clt}} \label{sec:clt:proof}

We are now in position to prove \eqref{eq:thmclt1} in Theorem \ref{thm:clt}. 
By conditioning on the event $\mathcal{G}^{n,h}$, that at least $\log n$ households are infected in the epidemic, it follows from Lemma \ref{lem:lln}, Corollary \ref{cor:lln} and the discussion after Corollary \ref{cor:lln}  that
\[ S_n (\bar{T}_\infty^n)|\mathcal{G}^{n,h}  \convp \tau  \qquad \mbox{and} \qquad \bar{Z}_{n,h} |\mathcal{G}^{n,h}   \convp z \left( = \nu_R (\tau) \right) \qquad \mbox{as } n \to \infty. \]
Throughout the remainder of the proof we implicitly condition on $\mathcal{G}^{n,h}$.

Considering $\bar{Z}_{n,h} | \mathcal{G}^{n,h}$ directly is not straightforward. However, we note that condtional upon $\mathcal{G}^{n,h}$ at least $k_n = \lfloor \log n \rfloor$ households are infected. This allows us to construct lower, $\bar{Z}_{n,h}^L$, and upper, $\bar{Z}_{n,h}^U$, bounds for the proportion infected in the event of a global epidemic by considering who becomes infected in the first $k_n$ households and {\it restarting} the epidemic with $k_n$ infected households and $n-k_n$ initially susceptible households. We show that $\bar{Z}_{n,h}^L \sle \bar{Z}_{n,h} | \mathcal{G}^{n,h} \sle \bar{Z}_{n,h}^U$ with $\sqrt{nh} (\bar{Z}_{n,h}^L - z) \convD N(0,\sigma^2)$ and $\sqrt{nh} (\bar{Z}_{n,h}^U - z) \convD N(0,\sigma^2)$  from which \eqref{eq:thmclt1}, and hence Theorem \ref{thm:clt}, follow.

In order to obtain suitable $\bar{Z}_{n,h}^L$ and $\bar{Z}_{n,h}^U$, we first define a sequence of epidemic processes $\hat{\mathcal{E}}_n (D_n)$, indexed by the number of households $n$, where $D_n \in \mathbb{N}$ denotes the number of global contacts from outside the population to initiate the epidemic. Suppose that in $\hat{\mathcal{E}}_n (D_n)$ there are initially $n-k_n$ totally susceptible households, with the remaining $k_n$ households consisting entirely of removed individuals. We label the initially susceptible households $1,2,\ldots, n-k_n$ and the initially removed households $n+1-k_n, n+2- k_n, \ldots, n$. The epidemic is constructed  in a similar manner to $\mathcal{E}_n$ using $\{ \mathbf{X}_{ij}, \mathbf{H}_{ij}, \eta_{ij}; i=1,2,\ldots,n, j=1,2, \ldots, h\}$ with $D_n$ initial global contacts to determine the initial infectives within the population. However, throughout the epidemic global contacts with households $n+1-k_n, n+2- k_n, \ldots, n$ have no effect as they are with removed individuals. Thus for the initially removed households only the $\eta_{ij}$s are required.  The epidemic only effectively takes place between $n- k_n$ households with the $k_n$ initially removed households included to absorb unsuccessful global infections when we relate $\hat{\mathcal{E}}_n (D_n)$ to $\mathcal{E}_n$. Due to the construction of $\hat{\mathcal{E}}_n (D_n)$, we can use the trivariate random vectors $(R_i (t), G_i (t), Y_i (t))$ for the embedding process. Let $\hat{T}_\infty^n (D_n)$ satisfy
\begin{align} \label{eq:hatT}
\hat{T}_\infty^n (D_n) = \min \left\{ x \geq 0  :   D_n +  \sum_{i=1}^{n-k_n} G_i (S_n (x)) = \lfloor x nh \rfloor \left( = \sum_{i=1}^n Y_i (S_n (x)) \right)  \right\}.
\end{align}
We note that the difference between \eqref{eq:hatT} and \eqref{eq:Tnx} for $\bar{T}_\infty^n$ is the number of global infections to initiate the epidemic and that in  $\hat{\mathcal{E}}_n (D_n)$  only $n-k_n$ households contribute to the generation of new global infections. By construction if $D_n^1 < D_n^2$ then $\hat{T}_\infty^n (D_n^1) \leq \hat{T}_\infty^n (D_n^2)$.

We have the following central limit theorem for the proportion infected in the epidemic  $\hat{\mathcal{E}}_n (D_n)$ with the proof deferred to Section \ref{sec:clt:proofhat}. Theorem \ref{thm:cltbound} is central to proving \eqref{eq:thmclt1}.  (We show in Section~\ref{sec:clt:variance} that the expression for $\sigma^2$ given in~\eqref{eq:sigma3} below is equivalent to that given in~\eqref{eq:sigma} in Section~\ref{sec:SSHEmain:result}.)

\begin{theorem} \label{thm:cltbound}
Let $D_n$ be a sequence of positive integers such that $D_n \to \infty$ and $D_n/\sqrt{n} \to 0$ as $n \to \infty$. Let
\begin{align*}
\hat{Z}_{n,h} (D_n)= \frac{1}{nh} \sum_{i=1}^{n-k_n} R_i (S_n ( \hat{T}_\infty^n(D_n))),
\end{align*}
the proportion of individuals who are infected during the epidemic $\hat{\mathcal{E}}_n (D_n)$. Then
\begin{align}
\sqrt{n h} \left(\hat{Z}_{n,h} (D_n) - z\right) \convD {\rm N} (0, \sigma^2) \qquad \mbox{as } n \to \infty, \nonumber
\end{align}
where
\begin{equation}
\label{eq:sigma3}
\sigma^2 = \var \left( R_1 (\tau) + b (\tau) \left[ G_1(\tau) - Y_1 (\tau)\right] \right).
\end{equation}
\end{theorem}

Conditional on $\mathcal{G}^{n,h}$, we can consider the first $k_n$ households infected. Let $\tilde{D}_n^I$ denote the total number of global infectious contacts emanating from the first $k_n$ local household epidemics plus the initial $T_0^n$ global infectious contacts required to create the $m_n$ initial infectives.  Then $\tilde{D}_n^I/k_n \convp R_\ast $ as $n \to \infty$, where $R_\ast = \E [C]$ is the mean number of global contacts emanating from a local epidemic initiated by a single infective in an otherwise susceptible household.
Let $\tilde{D}_n^B$ denote the number of global infections required to infect $k_n$ distinct households and note that using the birthday problem $\P(\tilde{D}_n^B = k_n) \to 1$ as $n \to \infty$, cf. \ref{equ:birthday}. Let $\tilde{D}_n^A = \tilde{D}_n^I - \tilde{D}_n^B$ the number of excess global contacts between the number of global contacts required to infect the first $k_n$ households and the number of global infectious contacts generated by these first $k_n$ household epidemics. Then $\tilde{D}_n^A/k_n \convp R_\ast -1 >0$ as $n \to \infty$.  Let $\tilde{D}_n^C$ denote the sum of all the global contacts from individuals in the first $k_n$ infected households whether or not they are infected in the initial local epidemic in the household plus the initial $T_0^n$ global infectious contacts. Note that $\tilde{D}_n^C \eqd T_0^n + \sum_{i=1}^{k_n} \sum_{j=1}^h X_{G,(i,j)}$ with $\E[\tilde{D}_n^C] =\E [T_0^n] + h k_n \mu_G$. Let $\tilde{D}_n^D = \tilde{D}_n^C - \tilde{D}_n^B$, the  number of excess global contacts between the number of global contacts required to infect the first $k_n$ households and the total number of potential global infectious contacts generated by these first $k_n$ households should everybody become infected.

We create a lower bounding epidemic process $\bar{\mathcal{E}}_n^L$ by using the same construction as $\mathcal{E}_n$ except that in the first $k_n$ households to be infected only the first global contact is successful. All subsequent global infectious contacts with these $k_n$ households, which we denote by $\mathcal{F}_n$, are unsuccessful. For households in $\mathcal{F}_n^C$, the epidemic progresses as in $\mathcal{E}_n$. Let $\bar{Z}_{n,h}^L$ denote the proportion of the population infected in  $\bar{\mathcal{E}}_n^L$ , then $\bar{Z}_{n,h}^L \leq \bar{Z}_{n,h}$. Similarly we  create an upper bounding epidemic process $\bar{\mathcal{E}}_n^U$ by using the same construction as $\mathcal{E}_n$ except that in the first $k_n$ households all individuals are made infectious. All subsequent global infectious contacts with these $k_n$ households have no effect, as the individual contacted has already been infected. For households in $\mathcal{F}_n^C$, the epidemic again progresses as in $\mathcal{E}_n$. Let $\bar{Z}_{n,h}^U$ denote the proportion of the population infected in  $\bar{\mathcal{E}}_n^U$ , then $\bar{Z}_{n,h}^U \geq \bar{Z}_{n,h}$. Let $\bar{Z}_{n,h}^L = \bar{Z}_{n,h}^{L,0} + \bar{Z}_{n,h}^{L,1}$, where $\bar{Z}_{n,h}^{L,0}$ is the proportion of the population who both  belong to $\mathcal{F}_n$ and are infected in $\bar{\mathcal{E}}_n^L$ and  $\bar{Z}_{n,h}^{L,1}$ is the proportion of the population who both  belong to $\mathcal{F}_n^C$ and are infected in $\bar{\mathcal{E}}_n^L$. Define  $\bar{Z}_{n,h}^{U,0}$ and  $\bar{Z}_{n,h}^{U,1}$ similarly, with $\bar{Z}_{n,h}^U = \bar{Z}_{n,h}^{U,0} + \bar{Z}_{n,h}^{U,1}$.

By construction, the lower bounding and upper bounding epidemic processes behave as if the households in $\mathcal{F}_n$ are removed after considering the first $k_n$ households and $\tilde{D}_n^B$ global infections. Given that Poisson processes have independent increments and $\mathcal{G}^{n,h}$ with $\tilde{D}_n^A = D_n^1$, we can couple the construction of $\bar{\mathcal{E}}_n^L$ to a realisation of $\hat{\mathcal{E}}_n (D_n^1)$ such that $\{\bar{Z}_{n,h}^{L,1} | \mathcal{G}^{n,h}, \tilde{D}_n^A =D_n^1 \} = \hat{Z}_{n,h} (D_n^1)$. Similarly, given that $\tilde{D}_n^D= D_n^2$, we can couple the construction of $\bar{\mathcal{E}}_n^U$ to a realisation of $\hat{\mathcal{E}}_n (D_n^2)$ such that $\{\bar{Z}_{n,h}^{U,1} | \mathcal{G}^{n,h}, \tilde{D}_n^D = D_n^2 \} = \hat{Z}_{n,h} (D_n^2)$.

Let $D_n^L = \lfloor (R_\ast-1) k_n/2 \rfloor$ and $D_n^U = \lfloor 2 h k_n \mu_G \rfloor$. Then $\P (D_n^L \leq \tilde{D}_n^A) \to 1$ and $\P (D_n^U \geq \tilde{D}_n^D) \to 1$ as $n \to \infty$. Also since $D_n^L, D_n^U \to \infty$ and $D_n^L/\sqrt{n}, D_n^U/\sqrt{n} \to 0$ as $n \to \infty$, it follows from Theorem \ref{thm:cltbound} that both $\sqrt{n h} \left(\hat{Z}_{n,h} (D_n^L) - z\right) $ and $\sqrt{n h} \left(\hat{Z}_{n,h} (D_n^U) - z\right) $ converge in distribution to $N(0,\sigma^2)$ as $n \to \infty$.

Let $\bar{Z}_{n,h}^0$ ($\bar{Z}_{n,h}^1$) denote the proportion of the population who both belong to $\mathcal{F}_n$ ($\mathcal{F}_n^C$) and are infected in $\mathcal{E}_n$. We have that if $D_n^L \leq \tilde{D}_n^A$ and $D_n^U \geq \tilde{D}_n^D$,
\[  \hat{Z}_{n,h} (D_n^L) \leq \bar{Z}_{n,h}^{L,1} | \mathcal{G}^{n,h} \leq \bar{Z}_{n,h}^1 | \mathcal{G}^{n,h} \leq \bar{Z}_{n,h}^{U,1} | \mathcal{G}^{n,h} \leq \hat{Z}_{n,h} (D_n^U). \]
Given that $\P (\{ D_n^L \leq \tilde{D}_n^A\} \cup \{ D_n^U \geq \tilde{D}_n^D\}) \to 1$ as $n \to \infty$, it follows that $\sqrt{nh} (\bar{Z}_{n,h}^1 - z) \convD N(0,\sigma^2)$ as $n \to \infty$. Finally, \eqref{eq:thmclt1} follows using Slutsky's theorem (see, for example, Billingsley \cite{Bill99}, Theorem 3.1), since $\sqrt{nh} \bar{Z}_{n,h}^0 \convp 0$ as $n \to \infty$.

\subsection{Proof of Theorem \ref{thm:cltbound}} \label{sec:clt:proofhat}

In order to prove Theorem \ref{thm:cltbound}, we show that  $\sqrt{nh} (\hat{Z}_{nh} (D_n)- z)$ has the same limiting distribution, as $n \to \infty$, as the normalised sum of a certain linear combination of $\{ (R_i (\tau), G_i (\tau), Y_i (\tau)); i=1,2,\ldots,n \}$. This requires first defining for
$T >0$ a sequence of stochastic processes $\mathbf{W}_{[n,T]}$ and showing in Lemma \ref{lem:cltW} that the limiting stochastic process is  a zero-mean trivariate Gaussian process.

For $J = R, G, Y$ and $t \geq 0$, let
\begin{align} \label{eq:W:def}
W_n^J (t) = \frac{1}{\sqrt{nh}} \sum_{i=1}^{n_J} \left[ J_i (t) - h \nu_J (t) \right],
\end{align}
where  $\nu_R (t)$ is defined in \eqref{eq:nuR}, $\nu_G (t) = \mu_G \nu_R (t)$,  $\nu_Y (t) =t$, $n_R = n_G = \hat{n}= n-k_n$ and $n_Y =n$. That is, for $R$ and $G$ we sum over the $n-k_n$ initially susceptible households and for $Y$ we sum over all $n$ households, since global contacts with the initially susceptible households are important. Let $\mathbf{W}_n (t) = (W^R_n (t), W^G_n (t), W^Y_n (t))$ and let, for $T >0$,
\begin{align} \label{eq:WT:def}
\mathbf{W}_{[n,T]} = \{ \mathbf{W}_n (t) : 0 \leq t \leq T \}.
\end{align}
Also for $T >0$, let $\mathbf{W}_{[\ast,T]} =(W^R, W^G, W^Y)$ be a zero-mean trivariate Gaussian process with, for $J, L \in \{R,G,Y \}$ and $0 \leq s, t \leq T$,
\[ \cov (W^J (s), W^L (t)) = \frac{1}{h} \cov (J_1 (s), L_1 (t)). \]

\begin{lem}
\label{lem:cltW}
For any $T \geq 0$,
\[ \mathbf{W}_{[n,T]}\convw \mathbf{W}_{[\ast,T]} \qquad \mbox{as } n \to \infty \]
where $\convw$ denotes weak convergence in the space of bounded functions from $[0,T]$ to $\mathbb{R}^3$ endowed with the supremum metric (see, van der Vaart and Wellner \cite{vW}, page 34).
\end{lem}
\begin{proof}
Fix $T >0$. The lemma follows using \cite{vW}, Theorem 1.5.4, by showing that the finite-dimensional distributions of $ \mathbf{W}_{[n,T]}$ converge to those of $\mathbf{W}_{[\ast,T]}$ and that the sequence $\mathbf{W}_{[n,T]}$ $(n=1,2,\ldots)$ is asymptotically tight.

For any $m \in \mathbb{N}$, $\mathbf{t} \in [0,T]^m$ and $\alpha_{Jk} \in \mathbb{R}$ $(J=R,G,Y;k=1,2, \ldots, m)$,
\begin{align}
\label{eq:Qdef1}
& \sum_{k=1}^m \left\{ \alpha_{Rk} W_{n}^R (t_k) +\alpha_{Gk} W_{n}^G (t_k) + \alpha_{Yk} W_{n}^Y (t_k)   \right\}\\
& \qquad = \frac{1}{\sqrt{nh}}\sum_{i=1}^{\hat{n}} Q_i (\bmalpha, \mathbf{t})+\sum_{k=1}^m  \frac{1}{\sqrt{nh}} \sum_{i=\hat{n}+1}^n  \alpha_{Yk} [Y_i (t_k)  -  h \nu_Y (t_k)],\nonumber
\end{align}
where, for $i=1,2,\dots$,
\[
 Q_i (\bmalpha, \mathbf{t}) =\sum_{k=1}^m \left\{\alpha_{Rk} [R_i(t_k)-h \nu_R (t_k)]+\alpha_{Gk} [G_i (t_k)-h \nu_G (t_k)] +\alpha_{Yk} [Y_i (t_k)  - h \nu_Y (t_k)]  \right\}.
\]
The $\{Q_i  (\bmalpha, \mathbf{t})\}$s are i.i.d.~with $\E [Q_1 (\bmalpha, \mathbf{t})]=0$.
Since, for any $t \geq 0$, $R_1 (t) \leq h$, $G_1 (t) \leq \sum_{j=1}^h C_{1j}^G$ and $Y_1 (t) \sim {\rm Po} (ht)$, it is straightforward to show that
\[ \E \left[ Q_1  (\bmalpha, \mathbf{t})^2\right] < \infty. \]
Therefore, since $k_n/\sqrt{n} \to 0$ as $n \to \infty$, the central limit theorem yields
\[
\frac{1}{\sqrt{nh}} \sum_{i=1}^{\hat{n}} Q_i (\bmalpha, \mathbf{t}) \convD N \left( 0, \frac{1}{h} \var (Q_1  (\bmalpha, \mathbf{t}))\right) \qquad \mbox{as } n \to \infty.
\]
It is straightforward to show that the final term on the right-hand side of  \eqref{eq:Qdef1} converges in probability to 0 as $n \to \infty$, so using Slutsky's theorem,
\[
\sum_{k=1}^m \left\{ \alpha_{Rk} W_{n}^R (t_k) +\alpha_{Gk} W_{n}^G (t_k) + \alpha_{Yk} W_{n}^Y (t_k)   \right\} \convD N \left( 0, \frac{1}{h} \var (Q_1  (\bmalpha, \mathbf{t}))\right) \quad \mbox{as } n \to \infty.
\]
By considering linear combinations of $\mathbf{W}_n (t_k)$ and using the Cram\'er-Wold device, it follows that the finite-dimensional distributions of $ \mathbf{W}_{[n,T]}$ converge to those of $\mathbf{W}_{[\ast,T]}$.

By \cite{vW}, Lemma 1.4.3, the sequence $\mathbf{W}_{[n,T]}$ $(n=1,2,\ldots)$ is asymptotically tight if and only if each of the sequences $W_{[n,T]}^J$ $(J=R,G,Y; n=1,2,\ldots)$ is asymptotically tight.
We start by showing that the sequence $W_{[n,T]}^G$ $(n=1,2,\ldots)$ is asymptotically tight.

For $t \geq 0$, let $\bar{G}_1 (t) = G_1 (t) - h \nu_G (t)$, with $W_n^G (t) = (nh)^{-1/2} \sum_{i=1}^{n_G} \bar{G}_i (t)$.
Since $\bar{G}_1 (\cdot), \bar{G}_2 (\cdot), \ldots$ are i.i.d., the 3 conditions which are given for asymptotic tightness  of  $W_{[n,T]}^G$ $(n=1,2,\ldots)$ in \cite{vW}, Theorem 2.11.9 simplify to showing that as $n \to \infty$:
\begin{itemize}
\item[(i)] For every $\xi >0$,
\[  n_G \E \left[ \left\| \frac{\bar{G}_1}{\sqrt{nh}} \right\|_T 1_{\{ \|[\bar{G}_1 /\sqrt{nh} \|_T > \xi \}} \right] \rightarrow 0, \]
where $\| f \|_T = \sup_{0 \leq t \leq T} |f(t)|$.
\item[(ii)] For every $\delta_n \downarrow 0$,
\[ \sup_{|s-t| <\delta_n} \frac{n_G}{nh} \E \left[ \left( \bar{G}_1(s) - \bar{G}_1(t)\right)^2 \right] \rightarrow 0. \]
\item[(iii)]  For every $\delta_n \downarrow 0$,
\[ \int_0^{\delta_n} \sqrt{ \log N_{[\,]}^n (\epsilon,T)} \, d \epsilon   \rightarrow 0, \]
where for $\epsilon >0$,  the bracketing number $N_{[\,]}^n (\epsilon,T)$ is defined to be the minimum number of sets $N_\epsilon$ in a partition $[0,T] = \cup_{j=1}^{N_\epsilon} \mathcal{A}_{\epsilon j}^n$  such that, for each $ \mathcal{A}_{\epsilon j}^n$, we have
\begin{align} \label{eq:bracketdef}
\frac{n_G}{nh} \E \left[ \sup_{s,t \in \mathcal{A}_{\epsilon j}^n} \left( \bar{G}_1 (t) - \bar{G}_1 (s) \right)^2  \right] \leq \epsilon^2.
\end{align}
\end{itemize}

Given that, for all $t \geq 0$, $|\bar{G}_1 (t)| \leq \sum_{j=1}^h X_{G,(1,j)} + h \mu_G =Q^G$, say, it follows that, for any $\xi >0$,
\begin{align} \label{eq:tight1}
 n_G\E \left[ \left\| \frac{\bar{G}_1}{\sqrt{nh}} \right\|_T 1_{\{ \|[\bar{G}_1 /\sqrt{nh} \|_T > \xi \}} \right] & \leq \sqrt{\frac{n}{h}} \E \left[ Q^G 1_{\{ Q^G > \sqrt{nh} \xi\}}\right].
\end{align}
The same argument as a proof of Markov's inequality yields, for $a >0$,
\begin{align} \label{eq:tight2}
\E \left[ Q^G 1_{\{ Q^G > \sqrt{nh} \xi\}}\right] =(\sqrt{nh} \xi)^{-(1+a)} \E   \left[ Q^G 1_{\{ Q^G > \sqrt{nh} \xi\}} (\sqrt{nh} \xi)^{(1+a)} \right] \leq \frac{ \E   \left[ (Q^G)^{2+a} \right]}{(\xi \sqrt{nh})^{1+a}}.
\end{align}
Since $\E [X_G^{2+a}] < \infty$ implies that $\E  [(Q^G)^{2+a}] < \infty$, it is straightforward to show condition (i) holds using \eqref{eq:tight1} and \eqref{eq:tight2}.

Given that $G_1 (\cdot)$ and $\nu_G (\cdot)$ are non-decreasing in $t$, it is straightforward to show that for any $u \leq s \leq t \leq v$,
\begin{align} \label{eq:tight3}
\left[\bar{G}_1 (t) - \bar{G}_1 (s) \right]^2 & \leq \left[G_1 (t) - G_1 (s) \right]^2 +h^2 \left[ \nu_G (t) - \nu_G (s) \right]^2 \\
&\leq \left[G_1 (v) - G_1 (u) \right]^2 + h^2\left[ \nu_G (v) - \nu_G (u) \right]^2.  \label{eq:tight3a}
\end{align}
Also, jumps in $G_1 (\cdot)$ only occur when a global infectious contacts are made with the household, so, for all $0 \leq s < t$,
 \begin{align} \label{eq:tight4}
|G_1 (t) - G_1 (s)|  \leq \left( \sum_{j=1}^h X_{G,(1,j)} \right) 1_{\{ Y_i(t) \neq Y_i(s)\}}
\end{align}
with $\sum_{j=1}^h  X_{G,(1,j)}$ independent of $1_{\{ Y_i(t) \neq Y_i(s)\}}$.
It then follows from \eqref{eq:tight3}, \eqref{eq:tight4} and
\[ h^2 (\nu_G (t) - \nu_G (s))^2 = h^2 \left( \frac{1}{h} \E \left[ G_1 (t) - G_1(s) \right] \right)^2 \leq \E [( G_1(t) - G_1(s))^2], \]
that for all $0 \leq s < t$,
 \begin{align} \label{eq:tight5}
\E \left[ \left( \bar{G}_1(t) - \bar{G}_1(s)\right)^2  \right] & \leq 2 \E \left[ \left( G_1(t) - G_1(s)\right)^2  \right]  \nonumber \\ &\leq 2  \E \left[  \left( \sum_{j=1}^h X_{G,(1,j)} \right)^2 \right] \E \left[1_{\{ Y_1(t) \neq Y_1(s)\}} \right] \nonumber \\
&\leq 2 \E \left[  \left( \sum_{j=1}^h X_{G,(1,j)} \right)^2 \right] \E \left[Y_1(t) - Y_1(s)\right] \nonumber \\
& \leq  2h^2 \E \left[ X_G^2 \right] h(t-s).
\end{align}
Condition (ii) follows since for all $s \geq 0$, the right hand side of \eqref{eq:tight5} converges to 0 as  $t \downarrow s$.

Fix $\epsilon >0$ and $\mathcal{A} = [u,v]$, where $0 \leq u < v$ such that  $|u-v| \leq \epsilon^2/(4h^2 \E [X_G^2 ])$.  It follows from \eqref{eq:tight3a} and \eqref{eq:tight5} that
\begin{align} 
\frac{n_G}{nh} \E \left[ \sup_{s,t \in \mathcal{A}} \left( \bar{G}_1 (t) - \bar{G}_1 (s) \right)^2  \right] & \leq\frac{1}{h} \E \left[(G_1 (v) - G_1 (u))^2 \right] + h  \left[ \nu_G (v) - \nu_G (u) \right]^2 \nonumber \\
& \leq \frac{2}{h}  \E \left[(G_1 (v) - G_1 (u))^2 \right]  \nonumber \\
& \leq \frac{2}{h} \times 2h^3 \E \left[ X_G^2 \right] \times  \frac{\epsilon^2}{4h^2 \E [ X_G^2 ]} = \epsilon^2. \nonumber
\end{align}
Therefore a partition of $[0,T]$ into intervals $\mathcal{A}_{\epsilon j}^n$ of length $L_\epsilon = \epsilon^2/(4 h^2 \E [X_G^2 ])$ exists such that \eqref{eq:bracketdef} holds. Hence, $N_{[\,]}^n (\epsilon,T) \leq c/\epsilon^2$, where $c =1+ 4 T h^2 \E [ X_G^2 ]$. Then,
\begin{align} 
\int_0^{\delta_n} \sqrt{\log N_{[\,]}^n (\epsilon,T)} \, d \epsilon &\leq \int_0^{\delta_n} \sqrt{\log \left( \frac{c}{\epsilon^2}\right) } \, d \epsilon \nonumber \\
& = \frac{\sqrt{c}}{2} \int_{\log(c/\delta_n^2)}^\infty \sqrt{u} \exp \left(- \frac{u}{2} \right) \, du \to 0 \qquad \mbox{as } n \to \infty. \nonumber
\end{align}
Hence condition (iii) is satisfied, concluding the proof of asymptotic tightness of $W_{[n,T]}^G$ $(n=1,2,\ldots)$.

The asymptotic tightness of $W_{[n,T]}^R$ $(n=1,2,\ldots)$ follows by an identical argument with $X_G \equiv 1$. Finally, using properties of Poisson processes, it is straightforward to show that conditions (i)-(iii) hold with $\bar{G}_1$ replaced by $\bar{Y}_1$, where $\bar{Y}_1 (t) = Y_1 (t) - h \nu_Y (t)$ and $n_G = n-k_n$ replaced by $n_Y = n$. Therefore, the sequence $\mathbf{W}_{[n,T]}$ $(n=1,2,\ldots)$ is asymptotically tight and the lemma follows.
\end{proof}

\begin{proof}[Proof of Theorem \ref{thm:cltbound}]
Using similar arguments to Section \ref{sec:clt:lln}, we have that if $D_n \to \infty$ and $D_n/\sqrt{n} \to 0$ as $n \to \infty$, then
\begin{align}
\label{equ:SnTnCnconvp}
S_n (\hat{T}_\infty^n (D_n) ) \convp \tau, \qquad \mbox{as } n \to \infty.
\end{align}
This is because the probability that the epidemic fails to take-off from $D_n$ initial global contacts, of  which ${\rm Bin} (D_n, (n-k_n)/n)$ are with initially susceptible households, tends to 0 as $n \to \infty$.
Therefore, using similar arguments to Corollary \ref{cor:lln}, we have that $\hat{Z}_{n,h} (D_n) \convp z$ as $n \to \infty$.

Using the mean value theorem, we have that
\begin{align}
\label{eq:finalsize}
\sqrt{nh}& \left( \hat{Z}_{n,h} (D_n) -z \right) = \sqrt{nh} \left[ \frac{1}{nh}  \sum_{i=1}^{\hat{n}} R_i (S_n (\hat{T}^n_\infty (D_n))) -  \nu_R (\tau) \right]\nonumber \\
 & = \sqrt{nh} \left[ \frac{1}{nh}  \sum_{i=1}^{\hat{n}} R_i (S_n (\hat{T}^n_\infty (D_n))) -  \nu_R (S_n (\hat{T}^n_\infty (D_n)))  + \nu_R (S_n (\hat{T}^n_\infty (D_n))) - \nu_R (\tau) \right]  \nonumber \\
 & = W_n^R (S_n (\hat{T}^n_\infty (D_n))) + \tilde{a}_n + \nu_R^\prime (a_{n1}) \sqrt{nh} [ S_n (\hat{T}^n_\infty (D_n)) - \hat{T}^n_\infty (D_n) + \hat{T}^n_\infty (D_n) - \tau ],  \nonumber \\
 \end{align}
 where $a_{n1}$ lies between $S_n (\hat{T}^n_\infty (D_n)) $ and $\tau$ and $\tilde{a}_n = \sqrt{nh} [\hat{n} -n] \nu_R ( S_n (\hat{T}^n_\infty (D_n)))/n \convp 0$
  as $n \to \infty$.
By definition, see  \eqref{eq:hatT}, 
\[ \frac{1}{nh} \sum_{i=1}^n Y_i (S_n (\hat{T}^n_\infty (D_n)) ) =\hat{T}^n_\infty (D_n)
\qquad \mbox{and} \qquad
 \nu_Y (S_n (\hat{T}^n_\infty (D_n)) ) = S_n (\hat{T}^n_\infty (D_n)) . \]
Therefore, we can rewrite \eqref{eq:finalsize} as
\begin{align}
\label{eq:V2}
\sqrt{nh} \left(\hat{Z}_{n,h} (D_n) - z  \right) & =  W_n^R (S_n (\hat{T}^n_\infty (D_n))) + \tilde{a}_n -  \nu_R^\prime (a_{n1}) W_n^Y (S_n (\hat{T}^n_\infty (D_n)))  \nonumber \\ & \quad  \qquad+   \nu_R^\prime (a_{n1}) \sqrt{nh} [  \hat{T}^n_\infty (D_n)) - \tau ].
\end{align}
Hence, we need to consider the distribution of $\sqrt{nh} [  \hat{T}^n_\infty (D_n)) - \tau ]$.

Let $\hat{D}_n = D_n/(nh)$ and note that
\begin{align} 
 \sqrt{nh} [ \hat{T}^n_\infty (D_n)) - \tau ] & = \sqrt{nh} \left[ \hat{D}_n + \frac{1}{nh} \sum_{i=1}^{\hat{n}} G_i (S_n (\hat{T}^n_\infty (D_n))) - \nu_G (\tau) \right] \nonumber \\
&=  \sqrt{nh} \hat{D}_n  + W_n^G (S_n (\hat{T}^n_\infty (D_n))) + \tilde{a}_n \mu_G  \nonumber \\ & \quad + \sqrt{nh} [ \nu_G (S_n (\hat{T}^n_\infty (D_n)) - \nu_G (\hat{T}^n_\infty (D_n)) +  \nu_G (\hat{T}^n_\infty (D_n)) - \nu_G (\tau) ]. \nonumber
\end{align}
By the mean value theorem, there exists $a_{n2}$ lying between $S_n (\hat{T}^n_\infty (D_n))$ and $\tau$, such that
 \begin{align}  
&\sqrt{nh} [ \hat{T}^n_\infty (D_n) - \tau ] \nonumber \\ & \quad =\sqrt{n h} \hat{D}_n + W_n^G (S_n (\hat{T}^n_\infty (D_n))) + \tilde{a}_n \mu_G  \nonumber \\ & \qquad + \sqrt{nh} \nu_G^\prime (a_{n2}) [S_n (\hat{T}^n_\infty (D_n)) -\hat{T}^n_\infty (D_n)+ \hat{T}^n_\infty (D_n)- \tau ] \nonumber \\
& \quad = \sqrt{n h}\hat{D}_n + W_n^G (S_n (\hat{T}^n_\infty (D_n))) + \tilde{a}_n \mu_G - \nu_G^\prime (a_{n2}) W_n^Y(S_n (\hat{T}^n_\infty (D_n)))  \nonumber \\ & \qquad  +\nu_G^\prime (a_{n2}) \sqrt{nh} [\hat{T}^n_\infty (D_n)- \tau], \nonumber
 \end{align}
using~\eqref{eq:hatT}.  Hence,
\begin{align}
 \label{eq:Tn_tau}
& \sqrt{nh} [  \hat{T}^n_\infty (D_n) - \tau ] \nonumber\\ & = \frac{1}{1 -  \nu_G^\prime (a_{n2})} \left\{  \sqrt{n h} \hat{D}_n + W_n^G (S_n (\hat{T}^n_\infty (D_n))) + \tilde{a}_n \mu_G - \nu_G^\prime (a_{n2})  W_n^Y (S_n(\hat{T}^n_\infty (D_n))) \right\}. \nonumber \\
\end{align}
Inserting \eqref{eq:Tn_tau} into \eqref{eq:V2}, we obtain that
\begin{align}
 \label{eq:V3}
  \sqrt{nh} \left( \hat{Z}_{n,h} (D_n)- z \right) & =    W_n^R (S_n (\hat{T}^n_\infty (D_n)) + \tilde{a}_n +  \frac{\sqrt{nh} \hat{D}_n + \tilde{a}_n \mu_G}{ [1 -  \nu_G^\prime (a_{n2})] } \nonumber \\ & \quad + \left[ \frac{\nu_R^\prime (a_{n1})}{1 -  \nu_G^\prime (a_{n2})} \right] \left[ W_n^G (S_n (\hat{T}^n_\infty (D_n)))  -W_n^Y (S_n (\hat{T}^n_\infty (D_n))) \right].
 \end{align}

Using~\eqref{equ:SnTnCnconvp}, we have by the sandwich theorem that $a_{n1}, a_{n2} \convp \tau$ as $n \to \infty$. Therefore, $\nu_G^\prime (a_{n2}) \convp  \nu_G^\prime (\tau) <1$ as $n \to \infty$. Given that $\sqrt{n} \hat{D}_n \to 0$ and $\tilde{a}_n \convp 0$ as $n \to \infty$, the second and third terms on the right-hand side of \eqref{eq:V3} converge in probability to 0 as $n \to \infty$. Also, we have that
 \begin{align}
 \frac{\nu_R^\prime (a_{n1})}{1 -  \nu_G^\prime (a_{n2})}  \convp \frac{\nu_R^\prime (\tau)}{1- \nu_G^\prime (\tau)} \qquad \mbox{as } n  \to \infty. \nonumber
 \end{align}

It follows by Slutsky's theorem that $\sqrt{nh} \left( \hat{Z}_{n,h} (D_n) - z \right)$ and
 \begin{align*}  
   W_n^R (S_n (\hat{T}^n_\infty (D_n))) + b(\tau) \left[ W_n^G (S_n (\hat{T}^n_\infty (D_n)))   -  W_n^Y (S_n (\hat{T}^n_\infty (D_n))) \right] 
 \end{align*}
have the same limiting distribution, should one exist, as $n \to \infty$. 
By Slutsky's lemma and the continuous mapping theorem, \cite{vW}, Example 1.4.7 and Theorem 1.3.6, respectively, it follows from  Lemma \ref{lem:cltW} and~\eqref{equ:SnTnCnconvp} that
\begin{align*} 
\mathbf{W}_n (S_n (\hat{T}^n_\infty (D_n))) \convD \mathbf{W} (\tau) \qquad \mbox{as } n \to \infty .
\end{align*}
Hence,
\begin{align*}
& W_n^R (S_n (\hat{T}^n_\infty (D_n))) + b(\tau) \left[ W_n^G (S_n (\hat{T}^n_\infty (D_n)))   -  W_n^Y (S_n (\hat{T}^n_\infty (D_n))) \right]  \nonumber \\
& \qquad \convD  W^R (\tau) + b(\tau) \left[ W^G  (\tau)  -  W^Y  (\tau) \right] \qquad \mbox{as } n \to \infty,
\end{align*}
and the theorem follows since
\[ \sigma^2 = \var \left( W^R (\tau) + b(\tau) \left[ W^G  (\tau)  -  W^Y  (\tau) \right] \right) = \frac{1}{h} \var \left( R_1 (\tau) + b (\tau) \left[ G_1 (\tau) - Y_1 (\tau)\right]\right).\]
\end{proof}

\subsection{Variance calculations} \label{sec:clt:variance}

In this section we discuss  $\sigma^2$ and present an alternative representation of the variance.
The variance $\sigma^2$ satisfies
\begin{align} \label{eq:sigma1}
\sigma^2 &= \frac{1}{h} \var \left( R_1 (\tau) + b (\tau) [G_1 (\tau) - Y_1 (\tau)] \right) \\
& = (1 + b (\tau) \mu_G )^2 \nu_R (\tau) [1- \nu_R (\tau)] + (h-1) (1 + b (\tau) \mu_G )^2 \cov \left(\chi_{11} (\tau), \chi_{12} (\tau) \right)\nonumber\\
&\qquad + b (\tau)^2 \nu_R (\tau) [\sigma_G^2 - \mu_G] + 2 (h-1) b(\tau) (1 +\mu_G b(\tau) ) \cov (\chi_{11} (\tau), X_{G,(1,2)}),\label{eq:sigma2}
\end{align}
where $b(t) = \nu_R^\prime (t)/[1- \mu_G \nu_R^\prime (t)]$.
Since $\nu^G (\cdot)$ is a concave function, we have that $\nu^G (\tau) <1$ giving $b(\tau) < \infty$.

We make the following observations regarding $\sigma^2$ and defer showing that $\sigma^2$ satisfies \eqref{eq:sigma2} to Appendix \ref{app:sigma2}.
\begin{enumerate}
\item The expression for $\sigma^2$ in \eqref{eq:sigma2} involves simply the relationship by individuals $(1,1)$ and $(1,2)$. Note that if the number of global and local contacts made by individuals are independent then $\cov (\chi_{11} (\tau), X_{G,(1,2)}) =0$.
\item In the case $X_L \equiv 0$ (no local infection) we obtain a homogeneously mixing model with $\P (S_{11} =1)=1$, giving $\nu_R (t) = 1- e^{-t}$ and $\nu_R (t) + \nu_R^\prime (t) - 1 =0$. Therefore for $X_L \equiv 0$, letting $\xi = \exp (-\tau) [ = \nu_R^\prime (\tau)]$, we have that $\nu_R (\tau) =1 -\xi$, $b(\tau) = \xi /[1- \mu_G \xi ]$, $\tau = \mu_G (1-\xi )$ and $1+ b(\tau) \mu_G= [1-\mu_G \xi]^{-1}$, so
\begin{align} \label{eq:sigma:homo}
\sigma^2 &=(1+ b(\tau) \mu_G)^2 \nu_R (\tau) [1 -\nu_R (\tau)] +  b(\tau)^2  \nu_R (\tau) [\sigma_G^2 - \mu_G] \nonumber \\
& = \frac{[1-\xi] \xi}{(1 - \mu_G \xi)^2} + \frac{\xi^2}{(1 - \mu_G \xi)^2} [1-\xi] [\sigma_G^2 - \mu_G] \nonumber \\
& =  \frac{\xi (1-\xi ) + \xi ^2(1-\xi ) [\sigma_G^2 - \mu_G]}{(1-\mu_G \xi )^2}.
\end{align}
The expression for $\sigma^2$ given in \eqref{eq:sigma:homo} agrees with the variance term given in \cite{ML86}, Theorem 1, for a constant number of initial infectives $m$. The model considered in \cite{ML86} is the generalised Reed-Frost model, where infectious individuals make $X_G$ contacts with distinct members of the population. As we note in Section \ref{sec:clt:without} below the difference between sampling global contacts with and without replacement vanishes as $n \to \infty$.
\item The expression for $\sigma^2$ given at~\eqref{eq:sigma} in Section~\ref{sec:SSHEmain:result} is of course equivalent to~\eqref{eq:sigma1} or~\eqref{eq:sigma2} above, as is shown at the end of the proof of \eqref{eq:sigma2} in Appendix \ref{app:sigma2}.
\end{enumerate}

\subsection{Global and local contacts sampled without replacement} \label{sec:clt:without}

In this section, we briefly describe the minor modifications required for the central limit theorem to hold when the global and local contacts made by an infective are {\it without replacement} from the remainder of the population and household, respectively. As $n \to \infty$, the probability an infective makes either a global self-contact or multiple global contacts with a given individual converges to 0 provided that $\mu_G = \E [X_G] <\infty$. Moreover, it is straightforward to show that the total number of global self-contacts and multiple global contacts made by individuals with the same individual, $V_n$ say, satisfies
\begin{align*}
V_n \convD V \sim {\rm Po} \left( \frac{\E[X_G (X_G+1)]}{2} \right) \qquad \mbox{as } n \to \infty,
\end{align*}
provided that $\E[X_G^2] < \infty$, which is the case as under the  assumptions of Theorem \ref{thm:clt},  there exists $a>0$ such that $\E [X_G^{2+a}] <\infty$. The effect of $V_n$ additional global contacts to replace global self- and multiple contacts is negligible and does not affect the law of large numbers and central limit theorem for final proportion infected by a major epidemic. A similar result holds if we preclude the possibility of an individual making global contacts with their own household. 

Turning to local (household) infectious contacts, if $X_L$ denotes the total number of distinct household contacts then $X_L$ has support on $\{0,1,\ldots, h-1\}$. Consequently, $\mathbf{H}_{ij}$, the successive individuals contacted locally by individual $(i,j)$, is  a random vector of length $h-1$, whose entries are a random permutation of $\{1,2,\ldots h\} \backslash j$, with individual $(i,j)$ making a household infectious contact with individual $(i,l)$ if $l \in \{ H_{ij1}, H_{ij2}, \ldots, H_{ijX_{L,(i,j)}}\}$. The susceptibility set of individuals can then be constructed from $\{ (X_{L,(i,k)}, \mathbf{H}_{ik}) ; k=1,2,\ldots,h \}$ in a similar manner to Section \ref{sec:clt:Sellke} with the proof of the central limit theorem continuing unchanged. The only change is the values taken by $z$ and $\sigma^2$, which change owing to the different distribution of $S$, the size of a susceptibility set (see Remark~\ref{rmk:withoutreplacement} in Appendix~\ref{app:Calc}).  Note that the probability of a minor outbreak, $\rho$, also changes, since $C$ has a different distribution.

\section{Proofs of Theorems  \ref{thm:prob:major} and  \ref{thm:finalsize}}\label{Sec-proof}

\subsection{Proof of Theorem \ref{thm:prob:major}} \label{sec:proof:major}

In this section we prove that the probability of a major outbreak, $\pi^{(h,p)}$, is increasing in $h$ and $p$ for any random vector $(X_G, X_L)$. This is proved in two separate lemmas where we vary $h$ (keeping $p$ fixed) in Lemma~\ref{lem:major:h} and vary $p$ (keeping $h$ fixed) in Lemma \ref{lem:major:p}. Lemma \ref{lem:major:p} is proved under weaker assumptions on $\mathbf{X} = (X_G, X_L)$ and the independent replacement of local contacts by global contacts.

We show first that $\pi^{(h,p)}$ is increasing in $h$.  We assume without loss of generality that $p=0$.  Recall the single-household epidemic model from $\mathcal{E}_h^H (X_G,X_L)$ from Section \ref{sec:SSHEmain:result}.
Let $R^{(h)}$ be the size of that epidemic, including the initial infective, and $C^{(h)}$ be the number of global contacts that emanate from infectives in that epidemic.  Thus $C^{(h)}$ is the offspring random variable for the branching process, $\mathcal{B}^{(h)}$, which approximates the of the epidemic $\mathcal{E}_{n,h}  (X_G,X_L)$.  Let $\rho^{(h)}$ denote the extinction probability of $\mathcal{B}^{(h)}$.

\begin{lem}  \label{lem:major:h}
For a given contacts random vector $\mathbf{X} = (X_G, X_L)$, $\rho^{(h)}$ is strictly decreasing in $h$.
\end{lem}
\begin{proof}
It is immediate that $C^{(1)} \sle C^{(2)}$ and hence that $\rho^{(1)} \ge \rho^{(2)}$.  Let $(X_{G,k}, X_{L,k})$ $(k=1,2,\dots)$ be i.i.d.~copies of $(X_G, X_L)$ and $U_k$ $(k=1,2,\dots)$ be an independent sequence of independent ${\rm U}(0,1)$ random variables.  We use these random variables to construct a realisation of $C^{(h)}$ for each $h=2,3,\dots$, as follows.

Fix $h \ge 2$.  We determine $(R^{(h)}, C^{(h)})$ by considering the infectives in $\mathcal{E}_h^H (X_G,X_L)$ one at a time.  We use $X_{L,1}$ to determine the number of distinct local contacts, $Z_1^{(h)}$, made by the initial infective.  Precise details are given below.  If $Z_1^{(h)}=0$ the epidemic stops and $(R^{(h)}, C^{(h)})=(1, X_{G,1})$.  Otherwise, we take one of $Z_1^{(h)}$ newly infected individuals and use $X_{L,2}$ to determine the number of distinct contacts it makes with
the remaining $h-1-Z_1^{(h)}$ susceptibles.  We continue the process in the obvious fashion, stopping when we have run out of infectives to consider.  Let $W_0^{(h)}=1$ and $W_k^{(h)}=1+Z_1^{(h)}+Z_2^{(h)}+\dots+Z_k^{(h)}$ $(k=1,2,\dots)$.
Then $R^{(h)}=\min\{k \ge 1: W_k^{(h)}-k=0\}$ and $C^{(h)}=\sum_{i=1}^{R^{(h)}} X_{G,i}$.  For completeness we define $W_k^{(h)}=W_{R^{(h)}}$ for $k>R^{(h)}$.

To determine whether local contacts are with susceptibles, we treat the local contacts one at a time.  Suppose that just prior to the $l^{\rm th}$ local contact a total of $Y_l^{(h)}$ individuals have been infected (including the initial infective).  Then that local contact is with a susceptible if and only if $U_l \le (h-Y_l^{(h)})/(h-1)$; otherwise the contact is with a non-susceptible individual and does not result in a new infective. Since
$(h-y)/(h-1)\le (h+1-y)/h$ for $y=1,2,\dots$, it follows immediately from the construction that $W_k^{(h+1)}\ge W_k^{(h)}$ for $k=0,1,\dots$, whence $R^{(h+1)} \ge R^{(h)}$ and $C^{(h+1)} \ge C^{(h)}$.
Thus, $C^{(h)} \sle C^{(h+1)}$ and $\rho^{(h)} \ge \rho^{(h+1)}$.  The inequality $(h-y)/(h-1)\le (h+1-y)/h$ is strict for $y>1$, so $\P(C^{(h+1)}>C^{(h)})=1$, whence $\rho^{(h)} > \rho^{(h+1)}$.
\end{proof}

Consider a random vector $\mathbf{X} = (X_G, X_L)$ for the number of global and household contacts made by a typical individual. For $0 \leq p \leq 1$, let $\mathcal{B} (p)$ denote the branching process where a proportion $p$ of household contacts are converted to global contacts.
Throughout, the branching process approximation is based on the assumption that each global infectious contact (birth) is with a previously uninfected household. Therefore, each infected household is infected globally once. The local epidemic (within the household) is determined by the number of local contacts, distributed independently according to $X_L$, and $p$, the proportion of household contacts that are converted to global contacts. 
We allow for a general rule for the transferring of household to global contacts. Let $Y_T^{(p)}$ denote the number of household contacts transferred to global contacts, so that in $\mathcal{B} (p)$ the number of  global and household contacts made by a typical infective are distributed according to $(X_G + Y_T^{(p)}, X_L - Y_T^{(p)})$. Note that $Y_T^{(0)} =0$ and $Y_T^{(1)} = X_L$, and if $Y_T^{(p)} = Y_L^{(p)} \sim {\rm MixBin} (X_L,p)$ we are in the scenario described in Section \ref{sec:main:p}. For $0 \leq p < q \leq 1$, we assume that $Y_T^{(p)}|X_L \sle Y_T^{(q)}|X_L$, that is, a coupling exists such that at least as many household contacts are transferred for an individual in $\mathcal{B} (q)$ as for the corresponding individual in $\mathcal{B} (p)$.

\begin{lem}  \label{lem:major:p}
For a given household size $h$ and contacts random vector $\mathbf{X} = (X_G, X_L)$, the extinction probability, $\rho_p$, of the branching process $\mathcal{B} (p)$  is monotonically decreasing in $p$.
\end{lem}
\begin{proof}
Fix $0 \leq p < q \leq 1$. We prove the lemma by showing that $\rho_q \leq \rho_p$.

Construct $\mathcal{B}(p)$ and $\mathcal{B} (q)$ on a common probability space such that the $i^{th}$ individual in both process makes $X_{G,i} + X_{L,i}$ attempted births and $Y_{T,i}^{(p)} \leq Y_{T,i}^{(q)}$. We construct a lower bound branching process $\hat{\mathcal{B}} (p,q)$ in which the  $i^{th}$ individual makes $X_{G,i}+ Y_{T,i}^{(p)}$ global and $X_{L,i}- Y_{T,i}^{(q)}$ local contacts. Note that in $\hat{\mathcal{B}} (p,q)$ the  $i^{th}$ individual has $Y_{T,i}^{(q)}  -Y_{T,i}^{(p)}$ fewer contacts than its counterparts in $\mathcal{B}(p)$ and $\mathcal{B} (q)$ and we term these missing contacts, {\it ghost} contacts. Let $\rho_p, \rho_q$ and $\hat{\rho}_{p,q}$ denote the extinction probabilities in the branching processes $\mathcal{B}(p), \mathcal{B} (q)$ and $\hat{\mathcal{B}} (p,q)$, respectively. Let $(\hat{V}_{p,q}, \hat{W}_{p,q})$ denote the number of global and ghost contacts emanating from a typical infectious individual in  $\hat{\mathcal{B}} (p,q)$. Then
\begin{align*} 
\hat{\rho}_{p,q} &= \E \left[ \hat{\rho}_{p,q}^{\hat{V}_{p,q}}\right].
\end{align*}
We define
\begin{align*}
\hat{f}_{p,q} (\theta,s) = \E \left[ \theta^{\hat{V}_{p,q}} s^{\hat{W}_{p,q}} \right],
\end{align*}
the joint pgf of $(\hat{V}_{p,q}, \hat{W}_{p,q})$,
so $\hat{\rho}_{p,q}$ solves
\[ \theta = \hat{f}_{p,q} (\theta,1). \]
Also we have that $\rho_q$ solves
\begin{align} 
\theta &= \hat{f}_{p,q} (\theta,\theta) =  \E \left[ \theta^{\hat{V}_{p,q}+\hat{W}_{p,q}} \right], \nonumber
\end{align}
since all ghost contacts in $\hat{\mathcal{B}} (p,q)$ correspond to global contacts in $\mathcal{B} (q)$.

The ghost contacts in  $\hat{\mathcal{B}} (p,q)$ correspond to local contacts within the household in $\mathcal{B} (p)$. The additional $\hat{W}_{p,q}$ local contacts in $\mathcal{B}(p)$ will result in $\tilde{W}_{p,q} \leq \hat{W}_{p,q}$ additional infectives from whom to grow the epidemic. It is likely that $\tilde{W}_{p,q} < \hat{W}_{p,q}$ as some contacts could be with individuals who are already members of the local household epidemic and/or repeat contacts with a new individual. Let $P (\hat{V}_{p,q},\hat{W}_{p,q})$ denote the probability that the branching process goes extinct from those individuals infected by the additional $\hat{W}_{p,q}$ local contacts. Thus, $\rho_p$ solves
\begin{align*}
\theta & = \E \left[\theta^{\hat{V}_{p,q}}   P (\hat{V}_{p,q},\hat{W}_{p,q}) \right].
\end{align*}
The $\tilde{W}_{p,q}$ individuals will initiate a local epidemic in a household with at least one removed individual (the initial infective). The number of global infections emanating from the local epidemic from the $\tilde{W}_{p,q}$ is stochastically smaller than $\sum_{i=1}^{\tilde{W}_{p,q}} \tilde{V}_{p,i}$, where the $\tilde{V}_{p,i}$s are i.i.d.~copies of $\tilde{V}_p$, the number of global contacts emanating from a household where individuals have i.i.d.~contacts according to $(X_{G,i}+ Y_{T,i}^{(p)}, X_{L,i}- Y_{T,i}^{(p)})$ and households initially have 1 infective, 1 removed and $h-2$ susceptibles. Let $\tilde{\rho}_p$ solve
\begin{align*}
\theta &=  \E \left[ \theta^{\tilde{V}_p} \right],
\end{align*}
the extinction probability of a branching process where the offspring distribution is $\tilde{V}_p$. Then $\tilde{V}_p \sle V_p$, where $V_p$ is the number of global contacts emanating from a household where individuals have i.i.d~contacts according to $(X_{G,i}+ Y_{T,i}^{(p)}, X_{L,i}- Y_{T,i}^{(p)})$ and households initially have 1 infective and $h-1$ susceptibles. Thus, $\tilde{\rho}_p \geq \rho_p$ and for $0 \leq \theta \leq1$,
\begin{align} \label{eq:maj:bound}
 \E \left[\theta^{\hat{V}_{p,q}}   P (\hat{V}_{p,q},\hat{W}_{p,q}) \right] \geq \E \left[ \theta^{\hat{V}_{p,q}}   \tilde{\rho}_p^{\hat{W}_{p,q}}  \right] =  \hat{f}_{p,q} (\theta,\tilde{\rho}_p).
\end{align}
Let $\rho_\ast$ solve
\begin{align*}
\rho_\ast =  \hat{f}_{p,q} (\rho_\ast,\tilde{\rho}_p).
\end{align*}
Then by \eqref{eq:maj:bound} it follows that $\rho_p \geq \rho_\ast$.

Since $\tilde{\rho}_p \geq \rho_\ast$, it follows that
\begin{align} \label{eq:maj:zast-bound}
\rho_\ast \geq  \hat{f}_{p,q} (\rho_\ast,\rho_\ast).
\end{align}
Given that $\rho_q$ is the smallest solution in $[0,1]$ of $\theta = \hat{f}_{p,q} (\theta,\theta)$, an immediate consequence of \eqref{eq:maj:zast-bound} is that $\rho_\ast \geq \rho_q$, whence $\rho_p \geq \rho_q$, as required.
\end{proof}
We observe that Lemma \ref{lem:major:p} holds if we assume instead that the, $X_L$, local contacts made by an individual are without replacement, with each local contact made by an individual being equally likely to be with anybody in their household they have not previously contacted.

\subsection{Proof of Theorem \ref{thm:finalsize}} \label{sec:proof:finalsize}

In this section we prove that the final size of a major outbreak, $z^{(h,p)}$, is increasing in $h$ and $p$ for any random vector $(X_G, X_L)$, for which the pgf of $X_L$ is log-convex. As in Section \ref{sec:proof:major} we prove the result in two separate lemmas where we vary $h$ (keeping $p$ fixed) in Lemma \ref{lem:finalsize:h}. and vary $p$ (keeping $h$ fixed) in Lemma \ref{lem:finalsize:p}.

We show first that $z^{(h,p)}$ is increasing in $h$.  We assume without loss of generality that $p=0$ and for ease of notation write $z^{(h,0)}$ as $z^{(h)}$.
\begin{lem}  \label{lem:finalsize:h}
For a given contact random vector $\mathbf{X} = (X_G, X_L)$, with $\log (f_{X_L} (s))$ being convex, the final size of a major outbreak, $z^{(h)}$,  is strictly increasing in $h$.
\end{lem}
\begin{proof}
In the proof, we use the following way of sampling a ${\rm Bin}(n,1-q)$ random variable.  First sample $Z \sim {\rm Po}(\lambda)$, where $\lambda=-n \log q$.  Then place $Z$ balls independently and uniformly at random into $n$ boxes and let $Y$ be the number of boxes that contain at least one ball.  Then $Y \sim {\rm Bin}(n,1-q)$.  (The numbers of balls in the $n$ boxes are independent ${\rm Po}(-\log q)$ random variables.)  Note this implies that $Y \sle Z$.

The susceptibility set $\mathcal{S}^{(h)}$ of a typical individual $a$ in a household of size $h$ can be constructed as follows.  We first look to see which individuals make contact with $a$; there are $X_1^{(h)} \sim {\rm Bin}(h-1, 1-q_0^{(h)})$ such individuals, where
\[
q_0^{(h)}=\E\left[\left(\frac{h-2}{h-1}\right)^{X_L}\right]=f_{X_L}\left(1-\frac{1}{h-1}\right).
\]
If $X_1^{(h)}=0$, the process stops and ${S}^{(h)}=1$. Otherwise, we take one of the $X_1^{(h)}$ individuals that have been added to the susceptibility set, individual $b$ say, and look to see which of the remaining $h-1-X_1^{(h)}$ individuals make contact with $b$.  Each of these individuals have failed to make contact with $a$, so the probability they make contact with $b$ is $1-q_1^{(h)}$, where
\[
q_1^{(h)}=\frac{\E\left[\left(\frac{h-3}{h-1}\right)^{X_L}\right]}{\E\left[
\left(\frac{h-2}{h-1}\right)^{X_L}\right]}=\frac{f_{X_L}\left(1-\frac{2}{h-1}\right)}{f_{X_L}\left(1-\frac{1}{h-1}\right)}.
\]
The process is then continued in the obvious fashion.  Specifically, for $k=2,3,\dots,h-1$,
\[
X_k^{(h)}|X_1^{(h)},X_2^{(h)}, \dots, X_{k-1}^{(h)} \sim {\rm Bin}(h-1-X_1^{(h)}-X_2^{(h)}-\dots-X_{k-1}^{(h)},1-q_{k-1}^{(h)}),
\]
where
\begin{equation}
\label{equ:qkh}
q_k^{(h)}=\frac{\E\left[\left(\frac{h-(k+2)}{h-1}\right)^{X_L}\right]}{\E\left[
\left(\frac{h-(k+1)}{h-1}\right)^{X_L}\right]}=\frac{f_{X_L}\left(1-\frac{k+1}{h-1}\right)}{f_{X_L}\left(1-\frac{k}{h-1}\right)}.
\end{equation}
Let $Y^{(h)}_0=1$ and $Y_k^{(h)}=1+X_1^{(h)}+X_2^{(h)}+\dots+X_k^{(h)}$ $(k=1,2,\dots,h)$, where $X_h^{(h)}=0$.  Then $S^{(h)}\eqd \min\{k \ge 1:Y_k^{(h)}-k=0\}$.

For $k=0,1,\dots, h-2$, let $\lambda_k^{(h)}=-(h-1)\log q_k^{(h)}$.  Note that~\eqref{equ:qkh} holds also for $k=0$, since $f_{X_L}(1)=1$.
Hence, for $k=0,1,\dots, h-2$,
\[
\lambda_k^{(h)}=-(h-1)\log q_k^{(h)}=(h-1)g_{X_L}\left(\frac{k+1}{h-1}\right)-(h-1)g_{X_L}\left(\frac{k}{h-1}\right).
\]
where $g_{X_L}(x)=-\log f_{X_L}(1-x)$ $(0 \le x \le 1)$.  The function $g_{X_L}$ is concave, increasing and differentiable on $[0, 1]$ (recall that $\log f_{X_L}$ is convex).  For $k=0,1,\dots,h-2$,
\begin{equation}
\label{equ:intrep}
\lambda_k^{(h)}=(h-1)\int_{k/(h-1)}^{(k+1)/(h-1)} g_{X_L}'(y) \, {\rm d}y
=\int_{k}^{k+1} g_{X_L}'\left(\frac{u}{h-1}\right) \, {\rm d}u,
\end{equation}
where we have made the substitution $u=(h-1)y$.
Now $g_{X_L}'$ is decreasing on $(0,1)$, since $g_{X_L}$ is concave, so it follows from~\eqref{equ:intrep} that $\lambda_k^{(h)} \ge \lambda_k^{(h')}$ if $h>h'$.

It is immediate that $S^{(h)} \sge S^{(1)}$, for $h>1$, so suppose that $h>h'\ge 2$. We construct coupled realisations of $S^{(h)}$ and $S^{(h')}$, satisfying $S^{(h)} \ge S^{(h')}$, as follows.
Let $\tilde{Z}_k^{(h)}  \sim {\rm Po}(\lambda_k^{(h)})$ $(k=0,1,\dots,h-2)$ be independent random variables and define $\tilde{Z}_k^{(h')}$ $(k=0,1, \dots, h'-2)$ similarly.  Since, $\lambda_k^{(h)} \ge \lambda_k^{(h')}$ $(k=0,1,\dots, h'-2)$,
$\tilde{Z}_k^{(h)}$ and $\tilde{Z}_k^{(h')}$ can be coupled so that $\tilde{Z}_k^{(h)} \ge \tilde{Z}_k^{(h')}$ $(k=0,1,\dots, h'-2)$.  We show by induction that the processes $Y_k^{(h)}$ $(k \ge 0)$ and $Y_k^{(h')}$ $(k \ge 0)$ can be coupled so that $Y_k^{(h)} \ge Y_k^{(h')}$ for all $k=0,1,\dots, h'$, whence $S^{(h)} \sge S^{(h')}$. (Note that $S^{(h')}$ is necessarily $\le h'$.)

Now $Y_0^{(h)}=Y_0^{(h')}=1$.  Suppose that $Y_i^{(h)} \ge Y_i^{(h')}$ for $i=0,1,\dots k$, where $k \le h'-1$.  Let $y=Y_k^{(h)}$ and $y'=Y_k^{(h')}$, so $y \ge y'$.  We use the above balls-in-boxes approach to obtain a realisation of $X_{k+1}^{(h)}$.  We place $\tilde{Z}_k^{(h)}$ balls uniformly at random in $h-1$ boxes, labelled $1,2,\dots, h-1$.  Then $X_{k+1}^{(h)}$ is given by the number of boxes with label $\ge y$ which contain at least one ball.  A realisation of $X_{k+1}^{(h')}$ can be obtained similarly, using $\tilde{Z}_k^{(h')}$.
Let $X_{k+1}'^{(h')}$ be the number of boxes with label $\ge y$ that contain at least one ball in the realisation of $X_{k+1}^{(h')}$.  Now $\tilde{Z}_k^{(h)} \ge \tilde{Z}_k^{(h')}$ and $(h-y)/(h-1)>(h'-y)/(h'-1)$, so using a sequence of independent ${\rm U}(0,1)$ random variables as in the proof of Lemma~\ref{lem:major:h}, it is straightforward to couple $X_{k+1}^{(h)}$ and $X_{k+1}'^{(h')}$ so that $X_{k+1}^{(h)} \ge X_{k+1}'^{(h')}$, whence  $Y_k^{(h)} \ge Y_k^{(h')}$, as required. 

It follows immediately from the above argument that $f_{S^{(h)}}(s) \le f_{S^{(h')}}(s)$ $(0 \le s \le 1)$ if $h > h'$.  Moreover, it easily seen that this inequality is strict for $s \in [0, 1)$.  Hence, $z^{(h)}>z^{(h')}$ if $h > h'$.
\end{proof}

The proof of Lemma \ref{lem:finalsize:p} is similar to that of Lemma  \ref{lem:major:p}. In Lemma \ref{lem:major:p} we select a random typical individual and study the {\it forward} epidemic process of who is infected from the resulting epidemic. We couple this to a {\it forward} branching process and compute the probability of extinction of the branching process.  In Lemma \ref{lem:finalsize:p} we select a random typical individual and study the {\it backward} epidemic process of who, if infected, will infect our selected individual. That is, we identify the susceptibility set of the individual and couple this to a {\it backward} branching process and compute its probability of extinction. Dependencies in the backward branching process mean that we require conditions on  $(X_G, X_L)$, namely, that the pgf of $X_L$ is log-convex and $Y_L^{(p)} \sim {\rm MixBin} (X_L,p)$, each local contact is independently with probability $p$ replaced by a global contact.

For $0 \leq p \leq 1$, let $\mathcal{S}(p)$ denote the susceptibility set of a randomly chosen individual in a household of size $h$, where individuals have household contacts distributed according to $X_L$ and each local contact is replaced by a global contact independently with probability $p$. Let $S (p) = |\mathcal{S}(p)|$,  the size of the susceptibility set. Let $\mathcal{B}^B (p)$ denote the backward branching process where individuals (household susceptibility sets) have sizes independently distributed according to $S(p)$ and an individual, with a susceptibility set of size $S(p)$ has ${\rm Po} ([\mu_G + p \mu_L] S(p))$ offspring. 
The offspring of a household susceptibility set in $\mathcal{B}^B (p)$ correspond to the set of individuals, who if infected, will infect the household susceptibility set via a global infection.  Let $\rho_p^B$ denote the extinction probability of
$\mathcal{B}^B (p)$.  Note that $\rho_p^B$ satisfies
\[
\E\left[\re^{-(\mu_G+p\mu_L)S(p)(1-\rho_p^B)}\right]=\rho_p^B,
\]
so $z^{(h,p)}=1-\rho_p^B$; cf. \eqref{eq:z}.

\begin{lem}  \label{lem:finalsize:p}
For a given household size $h$ and contacts random vector $\mathbf{X} = (X_G, X_L)$, with $\log (f_{X_L} (s))$ being convex, if $Y_L^{(p)} \sim {\rm MixBin} (X_L,p)$, then $\rho_p^B$ is monotonically decreasing in $p$.
\end{lem}
\begin{proof}
Fix $0 \leq p < q \leq 1$. We prove the lemma by showing that $\rho_q^B \leq \rho_p^B$.

Construct $\mathcal{B}_B (p)$ and $\mathcal{B}_B (q)$ on a common probability space as follows. Attach to each individual a local contact random variable $X_L$ to be used to construct susceptibility sets in the household. For each (potential) local contact assign an independent $U \sim U(0,1)$ random variable and if $U \leq p$ ($U \leq q$) convert the local contact to a global contact in $\mathcal{B}_B (p)$ ($\mathcal{B}_B (q)$). Thus in $\mathcal{B}_B (q)$ each individual makes the same number or fewer local contacts than the corresponding individual in $\mathcal{B}_B (p)$. Each individual has {\it backward} global contacts to grow the branching process beyond the current household. Attach to each individual a random variable $X_B \sim {\rm Po} (\mu_G + \mu_L)$ of potential global contacts into the individual. To each (potential) global contact assign an independent $\tilde{U} \sim U(0,1)$ random variable and if $\tilde{U}  \leq  [\mu_G +p \mu_L]/ [\mu_G + \mu_L]$ ($\tilde{U}  \leq  [\mu_G +q \mu_L]/ [\mu_G+ \mu_L]$) the global contact is kept in $\mathcal{B}_B (p)$ ($\mathcal{B}_B (q)$). Thus in $\mathcal{B}_B (q)$ each individual has the same number or more global contacts in than the corresponding individual in $\mathcal{B}_B (p)$.

As in Lemma \ref{lem:major:p}, we construct a lower bound branching process $\hat{\mathcal{B}}_B (p,q)$ in which the  $i^{th}$ individual has the same number of global contacts in as the $i^{th}$ individual in $\mathcal{B}_B (p)$ and the same number of local contacts out
 as the $i^{th}$ individual in $\mathcal{B}_B (q)$. Let $\hat{\mathcal{S}}(p,q)$ denote the susceptibility set of a randomly chosen individual in the branching process $\hat{\mathcal{B}}_B (p,q)$ with $\hat{S}(p,q) = |\hat{\mathcal{S}}(p,q)|$. Then $\hat{S}(p,q) \eqd S(q)$.
More explicitly, by selecting a typical individual in a typical household we can construct realistions of $\mathcal{S}(p)$, $\mathcal{S}(q)$ and $\hat{\mathcal{S}}(p,q)$ such that $\hat{\mathcal{S}}(p,q) = \mathcal{S}(q) \subseteq \mathcal{S}(p)$, with $\hat{S}(p,q) = S(q) \leq S(p)$.

 Let $\hat{W}^{(p,q)}$ be the number of potential global infectious contacts made with individuals in $\hat{S}(p,q)$ in $\hat{\mathcal{B}}_B (p,q)$. Then
\begin{align*}
\hat{W}^{(p,q)} | \hat{S}(p,q) \sim {\rm Po} \left(\lambda_p \hat{S}(p,q) \right),
\end{align*}
where $\lambda_p = \mu_G + p \mu_L$.
Let $\hat{\rho}^B_{p,q}$ denote the extinction probability of the branching process  $\hat{\mathcal{B}}_B (p,q)$. Then
\begin{align*}
\hat{\rho}^B_{p,q} = \E \left[ \exp \left( - \lambda_p \left[ 1 -  \hat{\rho}^B_{p,q}  \right] \hat{S}(p,q) \right) \right].
\end{align*}
For $0 \leq \theta \leq 1$ and $s=1,2,\ldots$, let
\begin{align*}
f_p (\theta; s) = \exp \left( - \lambda_p [1-\theta] s\right).
\end{align*}
Hence, $\hat{\rho}^B_{p,q}$ solves
\begin{align*}
\hat{\rho}^B_{p,q} =  \E \left[   f_p (\hat{\rho}^B_{p,q}; \hat{S}(p,q)) \right] =  \E \left[ f_p (\hat{\rho}^B_{p,q}; S(q))\right].
\end{align*}
Similarly,
\begin{align} \label{eq:fs:rhoBq}
\rho^B_q= \E \left[ f_q (\rho^B_q; S(q))\right] = \E \left[ f_p (\rho^B_q; S(q))  \exp \left( - (q-p) \mu_L [1- \rho^B_q]  S(q)\right)\right] .
\end{align}

Let $\mathcal{V}^C$ denote the individuals in $\mathcal{S}^c(q)$ that make contact with members of $\mathcal{S}(q)$ in the construction of $\mathcal{S}(p)$. Let $V^C = |\mathcal{V}^C|$. For these $V^C$ individuals we can construct the restricted susceptibility set, $\bar{\mathcal{S}}_R$, from members of $\mathcal{S}^c(q)$.  In other words, the restricted susceptibility set precludes individuals in $\mathcal{S}(q)$. Let $\bar{S}_R = |\bar{\mathcal{S}}_R|$. (Note that if $V^C=0$ then  $\bar{\mathcal{S}}_R = \emptyset$.) Then
\begin{align} \label{eq:fs:barS}
\bar{S}_R | V^C, S_0 \sle \sum_{i=1}^{V^C} S_i(p),
\end{align}
where $S_1(p), S_2(p), \ldots$ are i.i.d.~according to $S(p)$. The justification for \eqref{eq:fs:barS} is as follows.
Recall the definition of $q_k^{(h)}$ at~\eqref{equ:qkh} and note that~\eqref{equ:intrep} implies $q_k^{(h)}$ is nondecreasing in $k$.  Let $\bar{X}_L$ denote the local infectious contact distribution of a member of $\mathcal{S}^c(q) \backslash \mathcal{V}^C$ who does not make any household contacts with $\mathcal{S}(q)$. Then $\bar{X}_L \sle X_L$, since $q_k^{(h)}$ is nondecreasing in $k$. Further, since $|\mathcal{S}^c(q)| \leq h-1$, we can couple the construction of the susceptibility set of one member of $\mathcal{V}^C$ with the construction of the susceptibility set of an individual in a new household of size $h$ where all individuals have local contact distributions according to $X_L$, so that the size of the susceptibility set in the latter case is no smaller than the  former case. If $V^C >1$ we can repeat the process in turn for each member of  $\mathcal{V}^C$ considering only those individuals in  $\mathcal{S}^c(q)$ who have not previously been added to $\bar{\mathcal{S}}_R$.

Let $P(V^C, S(q))$ denote the probability of extinction of a branching process with an atypical initial individual, whose suspectibility set is formed of  $\bar{\mathcal{S}}_R$, and subsequent individuals have susceptibility sets of size i.i.d.~according to $S(p)$, and each member of the susceptibility set has ${\rm Po} (\lambda_p)$ offspring.
Then it follows from \eqref{eq:fs:barS} that
\begin{align} \label{eq:fs:barSorder}
P(V^C, S(q)) \geq [\rho_p^B]^{V^C}.
\end{align}
Also we have that $\rho^B_p$ solves
\begin{align} \label{eq:fs:rhoBp}
\rho^B_p= \E \left[ f_p (\rho^B_p; S(q)) P(V^C, S(q))\right].
\end{align}
However, from \eqref{eq:fs:barSorder}, we have that
\begin{align}  \label{eq:fs:rhoBp2}
 \E \left[ f_p (\rho^B_p; S(q)) P(V^C, S(q))\right] & \geq \E \left[f_p (\rho^B_p; S(q)) (\rho^B_p)^{V^C} \right] \nonumber \\ &= \E \left[f_p (\rho^B_p; S(q)) \E \left[  (\rho^B_p)^{V^C} |S(q) \right] \right].
\end{align}
Therefore if $\rho_\ast$ is the smallest solution in $[0,1]$ of
\begin{align*}
\theta =  \E \left[f_p (\theta_p; S(q)) \theta^{V^C} \right] =  \E \left[f_p (\theta; S(q)) \E \left[  \theta^{V^C} |S(q)\right] \right],
\end{align*}
it follows from \eqref{eq:fs:rhoBp} and \eqref{eq:fs:rhoBp2} that $\rho_p^B \geq \rho_\ast$.

We complete the proof of the lemma by showing, for $0 \leq \theta \leq 1$, that
\begin{align} \label{eq:fs:thetaVC}
\E [   \theta^{V^C} |S(q)] \geq  \exp \left( - [1- \theta] (q-p) \mu_L S(q)\right).
\end{align}
Since then it follows that
\begin{align} \label{eq:fs:bound2}
 \E \left[f_p (\theta_p; S(q)) \theta^{V^C} \right] &=  \E \left[f_p (\theta; S(q)) \E \left[  \theta^{V^C} |S(q) \right] \right] \nonumber \\
& \geq \E \left[ f_p (\theta; S(q))  \exp \left( - [1- \theta] (q-p) \mu_L S(q)\right) \right] = \E \left[ f_q (\theta; S(q))\right],
\end{align}
and together with \eqref{eq:fs:rhoBq}, \eqref{eq:fs:bound2} implies that $\rho_\ast \geq \rho^B_q$, whence $\rho^B_p \geq \rho^B_q$, as required.

For $0 \leq s \leq 1$, let $f_{X_L,p} (s)$ be the pgf of ${\rm MixBin} (X_L,1-p)$, so
\begin{align}
f_{X_L,p} (s) &= \E \left[ \E \left[ s^{X_L,p} | X_L \right] \right] \nonumber \\
& = \E \left[ \left(p + (1-p) s \right)^{X_L} \right] = f_{X_L} \left( p + [1-p]s\right). \nonumber
\end{align}
Note that
\begin{align*}
f_{X_L,p} (s) = f_{X_L} \left( p + [1-p][1-s]\right)  = f_{X_L} \left(1 - [1-p]s \right)
\end{align*}
Now
\begin{align*}
V^C | S(q) \sim {\rm Bin} \left( h - S(q), 1 - r_{S(q)} \right)
\end{align*}
where for $k=1,2,\ldots, h-1$,
\begin{equation*}
r_k  = \frac{f_{X_L,p} \left(1 - k/[h-1] \right)}{f_{X_L,q} \left(1 - k/[h-1] \right)}
 = \frac{f_{X_L} \left(1 - (1-p) k/[h-1] \right)}{f_{X_L} \left(1 - (1-q) k/[h-1] \right)}
\end{equation*}
is the probability that  an individual fails to infect locally a given set of $k$ individuals, when the probability of a local contact being transferred to a global contact is $p$, given the individual fails to infect locally a given set of $k$ individuals, when the probability of a local contact being transferred to a global contact is $q$. Hence, using ${\rm Bin} (n,1-r) \sle {\rm Po} (-n \log r)$ and $g_{X_L}(s)=-\log f_{X_L}(1-s)$ $(0 \le s \le 1)$, we have that
\begin{align}
\E \left[ \theta^{V^C} | S(q) \right] & \geq \exp \left( \left[h-  S(q) \right] \log \left\{ \frac{f_{X_L} \left(1 - (1-p)  S(q)/[h-1] \right)}{f_{X_L} \left(1 - (1-q)  S(q)/[h-1] \right)} \right\} [ 1-\theta]\right) \nonumber \\
& = \exp \left( -  \left[h- S(q) \right] [1-\theta] \left\{g_{X_L}\left(\frac{(1-p) S(q)}{h-1}\right)-g_{X_L}\left(\frac{(1-q) S(q)}{h-1}\right) \right\} \right) \nonumber \\
& = \exp \left( -  \frac{h-  S(q)}{h-1}  [1-\theta] \{ q -p\}  S(q) g_{X_L}^\prime (\xi) \right), \nonumber
\end{align}
where $(1-q) S(q)/(h-1)< \xi < (1-p) S(q)/(h-1)$.
Now $g_{X_L}^\prime(\theta)$ is decreasing as $f_{X_L}$ is log-convex, so for $0 \leq \theta \leq 1$,
\begin{align*}
g_{X_L}^\prime (\theta) = \frac{f_{X_L}^\prime (1-\theta)}{f_{X_L} (1-\theta)} \leq  \frac{f_{X_L}^\prime (1)}{f_{X_L} (1)} = \mu_L.
\end{align*}
Hence,
\begin{align*}
\E \left[ \theta^{V^C} |  S(q)\right]  &\geq \exp \left( - \frac{h- S(q)}{h-1} [1- \theta] (q-p) \mu_L  S(q)\right) \\
 &\geq \exp \left( - [1- \theta] (q-p) \mu_L  S(q)\right),
\end{align*} proving \eqref{eq:fs:thetaVC} and completing the proof of the lemma.
\end{proof}

\section{Discussion}\label{Sec-discussion}

In the paper we analysed a stochastic household epidemic model characterized by the random vector $(X_G, X_L)$ describing the number of global and local (=household) contacts individuals have, all global contacts being uniform in the entire community and all local contacts uniform in the household. Large population properties of the epidemic model were derived for the probability and size of a major outbreak. Then it was shown that the outbreak probability increases the larger household are considered, and the more of the local contacts are transferred to global contacts. The corresponding monotonicity results for the limiting relative final size $z$ were shown to require conditions on the distribution of $X_L$ with counter examples provided when these conditions were not satisfied.

For ease and clarity of presentation we have assumed that all households are of the same size. It is trivial to extend the central limit theorem to the case of unequal sized households provided that there exists $h_{max}<\infty$ such that all households are of size at most $h_{max}$. Additional conditions on the household size distribution will be required to extend the central limit theorem to the case where there is no maximum household size, see for example \cite{BL01} Section 5. The monotonicity results with increasing household size are conjectured to hold if we replace increasing household size by a stochastically increasing household distribution. That is, if we have epidemics in two populations with the same $(X_G,X_L)$ and household size distributions $H_1$ and $H_2$, in populations 1 and 2, respectively, such that $H_1$ is stochastically smaller than $H_2$ then $\pi_1 < \pi_2$ and, provided that $X_L$ has a log-convex pgf, $z_1 < z_2$, where $\pi_k$ and $z_k$ $(k=1,2)$ are the probability of, and proportion infected in, a major outbreak in populations 1 and 2, respectively.

The somewhat surprising counter examples to the monotonicity result: bigger epidemics with larger households or when swapping local to global contacts, occurred when the number of local contacts $X_L$ had low or no randomness. For example, in a household of size 3 and $X_L\equiv 1$ this would mean that an individual who gets infected would certainly infect one but not both of its household members. From an applied point of view this seems like an exceptional case, so we believe the monotonicity results are valid in most real world situations.

In Ball {\it et al.}~\cite{BBN24} we analysed an epidemic model with two types of subgroups where each individual belongs to precisely one subgroup of each type. Therefore each type of subgroup forms a partition of the population and it was assumed all subgroups of a given type have a common size. We allowed for the possibility of overlap between subgroups, that is, the possibility of two or more individuals belonging to the intersection of a subgroup of type 1 and a subgroup of type 2. The model was defined by contact rates during the  infectious periods (rather than arbitrary random vector as in the current paper), leading to mixed-Poisson distributed contacts of different categories. A branching process approximation for the initial stages of the epidemic and a law of large number approximation for the final proportion infected  were derived for the model. Numerical investigations suggested that the final size is increasing in the size of both group structures, and also increases as the amount of overlap between the two group structures decreases. These results served as inspiration for the current paper, but here we simplified to having only one group structure. A relevant question is of course if the monotonicity can be proven also when there are two (possibly overlapping) group structures. It follows immediately that in the case that there is no overlap between the two subgroups, our results in the present paper carry over to the case with two (and more) group structures: the final size increases if either (or both) of the two subgroup sizes increases. This follows, since as noted in Ball and Neal \cite{BN02}, the construction of the susceptibility set, which now extends across multiple, and in the limit possibly infinitely many groups, alternates between the two types of subgroups, so the distribution of the size of a susceptibility set of a typical individual  is stochastically increasing as the size of subgroups increases. For the situation where the two group structures are partly overlapping it remains an open problem, as is the numerically motivated conjecture that the final size increases as the amount of overlap between the two group structures decreases.

The embedding argument employed in the proof of the central limit theorem in Section \ref{sec:clt} can be utilised to study a wide range of epidemic models. As has been noted above, the central limit theorem can be applied to extensions of the Reed-Frost epidemic model where individuals are assumed to make infectious contacts with members of the population without replacement, see \cite{ML86} and \cite{PL90}. In a household context the key elements of the proof evolve around deriving the joint distribution of the number infected in a household, $R(t)$, and the number of global infections out of the household, $G(t)$, given that there has been a specific number of infectious contacts into the household, $Y(t)$. Therefore the approach is applicable to a wider class of models including, for example, assuming that not every individual is infected the first time they are contacted by an infectious individual but instead assuming there is a distribution on the number of infectious contacts required to infect an individual. Beyond the household model, the embedding argument could be employed to central limit theorems for the final size of epidemics in other two-level mixing population structures such as the great circle epidemic model, \cite{BN03}, and  network epidemic models, Ball and Neal \cite{BN08}, allowing progress beyond the mixed-Poisson distributions of global and local contacts in these earlier works.

An important assumption in the current model is that the random number of (uniformly chosen) local contacts $X_L$ is independent of household size. Some household epidemic models are defind by assuming the contact rate to \emph{each} household member equals some constant $\beta_H$. The overall rate to infect household members if in a household of size $h$ then equals $(h-1)\beta_H$ (in our model the contact rate, or equivalently total number of contacts, is assumed to be the same irrespective of household size). In such a situation, the epidemic is easily shown to increase the larger the household size is. Most network epidemic models makes a similar assumption: the rate or probability of contacting a given neighbour is fixed and independent of the number of neighbours. For a network epidemic model to more closely mimic the current model a fixed overall rate of infecting neighbours would be required, which is then distributed uniformly among the neighbours. The effect would be that highly connected individuals are no longer necessarily super-spreaders to the same extent. Would such an epidemic increase if the mean degree increased? Would the final size increase if the degree distribution has a heavier tail? These are some interesting open questions.

\begin{appendix}
\section{Calculation of asymptotic properties of $\mathcal{E}_{n,h} (X_G, X_L)$}
\label{app:Calc}
In this appendix we outline calculation of the major outbreak probability, $1-\rho$, and the asymptotic mean, $z$, and scaled variance, $\sigma^2$ of the fraction infected by a major outbreak; see Theorem~\ref{thm:clt}.  We make extensive use of Gontcharoff polynomials (see Lef\`evre and Picard~\cite{LP90}).  Let $U=u_0, u_1, \dots$ be a sequence of real numbers.  The Gontcharoff polynomials associated with $U$, i.e.~$G_i(x|U)$ $(i=0,1,\dots)$ are defined by
\begin{equation}
\label{equ:gontdef}
\sum_{i=0}^n n_{[i]} u_i^{n-i} G_i(x|U)=x^n \qquad (n=0,1,\dots),
\end{equation}
where $n_{[i]}=n(n-1)\dots(n-i+1)$ denotes a falling factorial, with the convention that $n_{[0]}=1$.
Note that $G_0(x)=1$ $(x \in \mathbb{R})$ and that $G_i(x|U)$ $(i=1,2,\dots)$ can be computed recursively for fixed $x$.  Further,
$G_i(x|U)$ is a polynomial of degree $i$ and (Lef\`evre and Picard~\cite{LP90}, Property 2.4) for $0 \le j \le i$,
\begin{equation}
\label{equ:gontdiff}
G_i^{(j)}(x|U)=G_{i-j}(x|E^jU),
\end{equation}
where $G_i^{(j)}(x|U)$ denotes the $j^{\rm th}$ derivative of $G_i(x|U)$ and $E^jU$ is the sequence $u_j, u_{j+1}, \dots$.

For $h=1,2,\dots$ and $\pi \in [0,1]$, recall the epidemic model $\mcEt^H_h(X_G, X_L, \pi)$ defined in Section~\ref{sec:SSHEmain:result}.
Let $\tilde{R}^{(h)}$ be the number of individuals infected in $\mcEt^H_h(X_G, X_L, \pi)$ and $\tilde{G}^{(h)}$ be the total number of global contacts that emanate from those infectives.  Further, let $\tilde{S}^{(h)}=h-\tilde{R}^{(h)}$ be the number of susceptibles at the end of the epidemic. Note that if $\pi=\re^{-t}$ then $(\tilde{R}^{(h)}, \tilde{G}^{(h)}) \eqd (R(t), G(t))$, defined in Section~\ref{sec:SSHEmain:result}.  We derive expressions for $\E[\tilde{S}^{(h)}_{[i]}]$ $(i=1,2)$ and $\E[\tilde{S}^{(h)}\tilde{G}^{(h)}]$, from which $\nu_R(t)=h^{-1}\E[R(t)], \var(R(t))$ and $\cov(R(t), G(t))$ follow easily.

For $k=1,2, \dots,h-1$, let $A_k^{(h)}$ be the event that an infective in $\mcEt^H_h(X_G, X_L, \pi)$ fails to contact anyone in a given set of $k$ susceptibles in the household.  For $s \in [0,1]$, let $q_0(s)=\E[s^{X_G}]$ 
and
$q_k(s)=\E[s^{X_G} 1_{A_k^{(h)}}]$ $(k=1,2\dots, h-1)$.  Then,
\begin{equation}
\label{equ:qksdef}
q_k(s)=f_{X_G, X_L}\left(s, 1-\frac{k}{h-1}\right) \qquad (k=0,1,\dots, h-1),
\end{equation}
where $f_{X_G, X_L}$ is the joint pgf of $(X_G, X_L)$.  Let $\tilde{f}_h(s_1, s_2)=\E[s_1^{\tilde{S}^{(h)}} s_2^{\tilde{G}^{(h)}}]$ $(s_1, s_2, \in [0,1])$.  Then it follows using Ball~\cite{Ball19}, Theorem 3.3, that
\begin{equation}
\label{equ:ftilde}
\tilde{f}_h(s_1, s_2)=\sum_{i=0}^h h_{[i]} (q_i(s_2))^{h-i} \pi^i G_i(s_1|U(s_2)) \qquad (s_1, s_2, \in [0,1]),
\end{equation}
where the sequence $U(s_2)$ satisfies $u_i(s_2)=q_i(s_2)$ $(i=0,1,\dots, h-1)$.  (Note from~\eqref{equ:gontdef} that, for $i=1,2,\dots$, $G_i(x|U)$ is determined by $u_0, u_1, \dots, u_{i-1}$.)

\begin{remark}
Observe from~\eqref{equ:qksdef} that $q_k(s)$, and hence also $U(s_2)$ and $G_i(s_1|U(s_2))$, depends on $h$.  We have suppressed this dependence for ease of presentation but note that in the sequel all Gontcharoff polynomials need to be recalculated if the household size is changed.
\end{remark}

\begin{remark}
\label{rmk:withoutreplacement}
Note that if the local contacts made by an infective are without replacement then, for $k=0,1,\dots, h-1$,
\[
\P(A_k^{(h)}|X_L)=\frac{\binom{h-1-k}{X_L}}{\binom{h-1}{X_L}}=\frac{(h-1-X_L)_{[k]}}{(h-1)_{[k]}},
\]
so~\eqref{equ:qksdef} is replaced by
\[
q_k(s)=\E\left[s^{X_G}\frac{(h-1-X_L)_{[k]}}{(h-1)_{[k]}}\right]  \qquad (k=0,1,\dots, h-1).
\]
\end{remark}

Differentiating~\eqref{equ:ftilde} partially $k$ times with respect to $s_1$, using~\eqref{equ:gontdiff}, yields (see 
\cite{Ball19}, Proposition 3.1)
\[
\E[\tilde{S}^{(h)}_{[k]}]=\sum_{i=k}^h h_{[i]} q_i^{h-i} \pi^{i} G_{i-k}(1|E^k U) \qquad (k=1,2,\dots,h),
\]
where $q_i=q_i(1)=f_{X_L}\left(1-\frac{i}{h-1}\right)$ $(i=0,1,\dots, h-1)$ and $U$ satisfies $u_i=q_i$ $(i=0,1,\dots, h-1)$.

\begin{remark}
\label{rmk:susset}
Note from 
\cite{Ball19}, Lemma 3.1, that $\P(S^{(h)}=i)=(h-1)_{[i-1]}G_{i-1}(1|EU)q_i^{h-i}$ $(i=1,2,\dots,h)$, where $S^{(h)}$ is the size of the susceptibility set of a typical individual in a household of size $h$.  Setting $\pi=\re^{-t}$, so $R(t) \eqd h-\tilde{S}^{(h)}$, yields~\eqref{eq:nuR}, and also enables $\nu_R^\prime(t)$ to be computed easily.
\end{remark}

For a function $f:\mathbb{R}^2 \to \mathbb{R}$ and $(k_1, k_2) \in \mathbb{Z}_+^2$, let $f^{(k_1, k_2)}(s_1, s_2)$ denote the partial derivative of $f$ of order $k_1$ in $s_1$ and $k_2$ in $s_2$.  Then,
$\E[\tilde{S}^{(h)}\tilde{G}^{(h)}]=\tilde{f}_h^{(1,1)}(1, 1)$.  For $i=0,1,\dots, h-1$ and $(s_1, s_2) \in [0,1]^2$, let $\alpha_i(s_1,s_2)=G_i(s_1|EU(s_2))$.  Differentiating~\eqref{equ:ftilde} with respect to $s_1$, using~\eqref{equ:gontdiff} and noting that
$G_0^{(1)}(s_1|U(s_2))=0$, yields
\begin{equation}
\label{equ:tildefdiff}
\tilde{f}_h^{(1,0)}(s_1,s_2)=\sum_{i=1}^h h_{[i]} (q_i(s_2))^{h-i} \pi^i \alpha_{i-1}(s_1, s_2).
\end{equation}
Recalling that $q_i=q_i(1)$, differentiating~\eqref{equ:tildefdiff} with respect to $s_2$ yields
\[
\E[\tilde{S}^{(h)}\tilde{G}^{(h)}]=\sum_{i=1}^h h_{[i]} q_i^{h-i} \pi^{i}\alpha_{i-1}^{(0,1)}(1,1)+\sum_{i=1}^{h-1}h_{[i+1]}q_i^{(1)}(1)q_i^{h-i-1}\pi^i \alpha_{i-1}(1,1).
\]

Using~\eqref{equ:gontdef},
\begin{equation}
\label{equ:alphagont}
\sum_{i=0}^n n_{[i]} (q_{i+1}(s_2))^{n-i}\alpha_i(s_1, s_2)=s_1^n \qquad (n=0,1,\dots, h-1),
\end{equation}
whence, differentiating \eqref{equ:alphagont} partially with respect to $s_2$,
\begin{equation}
\label{equ:alphadiff}
\sum_{i=0}^{n-1} n_{[i+1]} q_{i+1}^{(1)}(1)q_{i+1}^{n-i-1}\alpha_i(1, 1)+\sum_{i=1}^n n_{[i]}q_{i+1}^{n-i}\alpha_i^{(0,1)}(1,1)=0\qquad(n=0,1,\dots,h-1).
\end{equation}
Now $\alpha_{i}(1,1)$ $(i=0,1,\dots,h-1)$ can be computed using~\eqref{equ:alphagont}, $\alpha_0^{(0,1)}=0$ and $\alpha_i^{(0,1)}$ $(i=1,2,\dots,h-1)$ can be computed using~\eqref{equ:alphadiff},
thus enabling $\E[\tilde{S}^{(h)}\tilde{G}^{(h)}]$ to be computed.

For $h=1,2,\dots$, let $C^{(h)}$ be the total number of global contacts that emanate from $\mcE_h^H(X_G, X_L)$ defined in Section~\ref{sec:SSHEmain:result}.
Then, using 
\cite{Ball19}, Theorem 3.3,
\[
f_{C^{(h)}}(s)=\sum_{i=0}^{h-1} (h-1)_{[i]} (q_i(s))^{h-i} G_i(1|U(s)) \qquad(s \in [0,1]),
\]
thus enabling $f_{C^{(h)}}(s)$, and hence $\rho$, to be computed.

The above enables the asymptotic properties of $\mathcal{E}_{n,h}  (X_G,X_L)$ to be computed.  To compute the asymptotic properties of $\mathcal{E}_{n,h}  (X_G,X_L,p)$, for $p \neq 0$, note that elementary calculation yields
\begin{align*}
f_{X_G^{(p)}, X_L^{(p)}}(s_1, s_2)&=f_{X_G, X_L}(s_1, ps_1+(1-p)s_2)\qquad ((s_1, s_2) \in [0,1]^2),\\
\E[X_G^{(p)}]&=\E[X_G]+p\E[X_L],\\
\var(X_G^{(p)})&=\var(X_G)+2p\cov(X_G, X_L)+p^2 \var(X_L)+p(1-p)\E[X_L].
\end{align*}

\section{Proof of Theorem~\ref{thm:pnear1}}
\label{app:pnear1}
For $h=1,2,\dots$ and $p \in [0,1]$, let $S^{(h,p)}$ denote the size of the (household) susceptibility set of a typical individual in $\mathcal{E}_{n,h}  (X_G,X_L,p)$.  The mean number of global contacts made by a typical individuals in
$\mathcal{E}_{n,h}  (X_G,X_L,p)$ is $\mu_G+p\mu_L$, so it follows from~\eqref{eq:z} that $z^{(h,p)}$ is given by the largest solution in $[0,1]$ of
\[
1-z=f_{S^{(h,p)}}(\re^{-z(\mu_G+p\mu_L)}).
\]
Suppose that $\alpha=\mu_G+\mu_L>1$, so the (homogeneously mixing) epidemic when $p=1$ is supercritical, and let $z_1$ be the unique solution of $1-z=\re^{-\alpha z}$ in $(0,1)$.  (Note that $z_1=z_{\rm hom}(\alpha)$, defined just before the statement of Theorem~\ref{thm:pnear1} in Section~\ref{sec:SSHEmain:result}.)
For $h=1,2,\dots$ let
\begin{equation}
\label{equ:ghp}
g_h(p)=f_{S^{(h,p)}}(\re^{-z_1(\mu_G+p\mu_L)}) \qquad (p \in [0,1]).
\end{equation}
The behaviour of $z^{(h,p)}$ near $p=1$ is determined by the derivatives of $g_h$ at $p=1$.

It follows from Remark~\ref{rmk:susset}, \eqref{equ:gontdef} and a little algebra that
\begin{equation}
\label{equ:sussetT}
\sum_{i=0}^n \binom{h-1-i}{n-i}v_i^{n+1-h} \P(S^{(h)}=i+1)=\binom{h-1}{n}\qquad (n=0,1,\dots, h-1),
\end{equation}
where
\[
v_i=q_{i+1}=f_{X_L}\left(1-\frac{i+1}{h-1}\right) \qquad (i=0,1,\dots,h-2).
\]
Now $f_{X_L^{(p)}}(s)=f_{X_L}(p+(1-p)s)$, so let
\[
v_i(p)=f_{X_L^{(p)}}\left(1-\frac{i+1}{h-1}\right)=f_{X_L}\left(1-\frac{(1-p)(i+1)}{h-1}\right) \qquad (p \in [0,1]).
\]
(Note that $v_i(p)$ depends also on $h$ but we suppress that dependence for ease of notation.)

Fix $h \ge 2$ and let $f_i(p)=\P(S^{(h,p)}=i+1)$ $(i=0,1,\dots h-1)$.  Then, using~\eqref{equ:sussetT},
\begin{equation}
\label{equ:fip}
\sum_{i=0}^n \binom{h-1-i}{n-i}(v_i(p))^{n+1-h} f_i(p)=\binom{h-1}{n}\qquad (n=0,1,\dots, h-1).
\end{equation}
Noting that $v_i(1)=1$ $(i=0,1,\dots, h-2)$, it follows from~\eqref{equ:fip} that
\begin{equation}
\label{equ:fip1}
f_0(1)=1\qquad \text{and} \qquad  f_i(1)=0 \quad (i=1,2,\dots,h-1).
\end{equation}
Now,
\[
\frac{\mathrm{d}}{\mathrm{d}p}\left.\left(v_i(p)^{n+1-h}\right)\right|_{p=1}=\frac{(n+1-h)(i+1)\mu_L}{h-1},
\]
so differentiating~\eqref{equ:fip} and using~\eqref{equ:fip1} yields
\begin{equation}
\label{equ:fipdiff1}
\sum_{i=0}^n \binom{h-1-i}{n-i} f_i^{(1)}(1)=\binom{h-1}{n}\left(\frac{h-1-n}{h-1}\right)\mu_L \qquad (n=0,1,\dots, h-1).
\end{equation}
Successively setting $n=0,1,2, \dots, h-1$ in~\eqref{equ:fipdiff1} yields, after a little algebra,
\begin{equation}
\label{equ:fipdiff1a}
f_0^{(1)}(1)=\mu_L,\quad f_1^{(1)}(1)=-\mu_L \qquad \text{and} \qquad  f_i^{(1)}(1)=0 \quad (i=2,3,\dots,h-1).
\end{equation}
Now
\[
\frac{\mathrm{d}^2}{\mathrm{d}p^2}\left.\left(v_0(p)^{n+1-h}\right)\right|_{p=1}=\frac{(h-1-n)(h-n)\mu_L^2}{(h-1)^2}-\frac{(h-1-n)\mu_{L,[2]}}{(h-1)^2},
\]
where $\mu_{L,[2]}=\E[X_L(X_L-1)]$. Differentiating~\eqref{equ:fip} twice yields, after using~\eqref{equ:fip1} and~\eqref{equ:fipdiff1a},
\begin{align}
\label{equ:fipdiff2}
\sum_{i=0}^n &\binom{h-1-i}{n-i} f_i^{(2)}(1)=2\frac{(h-1-n)}{(h-1)}\left[\binom{h-1}{n}-2\binom{h-2}{n-1}1_{\{n \ge 1\}}\right]\mu_L^2\\
&+\binom{h-1}{n}\frac{(h-1-n)}{(h-1)^2}[\mu_{L,[2]}-(h-n)\mu_L^2] \qquad (n=0,1,\dots, h-1)\nonumber.
\end{align}
Successively setting $n=0,1,2, \dots, h-1$ in~\eqref{equ:fipdiff2} yields, after a little algebra,
\begin{align}
\label{equ:fipdiff2a}
&f_0^{(2)}(1)=\frac{1}{h-1}[\mu_{L,[2]}+(h-2)\mu_L^2],\quad f_1^{(2)}(1)=-\frac{1}{h-1}[\mu_{L,[2]}+4(h-2)\mu_L^2],\\
& f_2^{(2)}(1)= 3\frac{h-2}{h-1}\mu_L^2 \qquad \text{and} \qquad  f_i^{(2)}(1)=0 \quad (i=3,4,\dots,h-1). \nonumber
\end{align}

Returning to $g_h$, note from~\eqref{equ:ghp} that
\begin{equation}
\label{equ:ghp1}
g_h(p)=\sum_{k=1}^{h} f_{k-1}(p) \re^{-kz_1(\mu_G+p\mu_L)}   \qquad (p \in [0,1]).
\end{equation}
Differentiating~\eqref{equ:ghp1} yields, after using~\eqref{equ:fip1} and~\eqref{equ:fipdiff1a},
\begin{equation}
\label{equ:ghpdiff1}
g_h^{(1)}(1)=\mu_L\re^{-\alpha z_1} (1-z_1-\re^{-\alpha z_1})=0,
\end{equation}
since $1-z_1=\re^{-\alpha z_1}$.
Differentiating~\eqref{equ:ghp1} twice, and using ~\eqref{equ:fip1}, \eqref{equ:fipdiff1a} and~\eqref{equ:fipdiff2a}, yields after a little algebra
\begin{equation}
\label{equ:ghpdiff2}
g_h^{(2)}(1)=\frac{z_1(1-z_1)}{h-1}[\mu_{L,[2]}+(2-3z_1)\mu_L^2]
=\frac{z_1(1-z_1)}{h-1}[\sigma_L^2-\mu_L+3(1-z_1)\mu_L^2].
\end{equation}

Recall that $h\ge 2$ is fixed and let $z(p)=z^{(h,p)}$ $(p \in [0,1])$.  Then~\eqref{equ:ghpdiff1} and~\eqref{equ:ghpdiff2} imply that $z^{(1)}(1)=0$ and
$z^{(2)}(1) >0 (<0)$ if $\sigma_L^2-\mu_L+3(1-z_1)\mu_L^2 >0 (<0)$, from which Theorem~\ref{thm:pnear1} follows easily.

\section{Proof of Theorem \ref{lem:hinf}} \label{app:hinf}

We prove the theorem in the case $p=0$, with a similar proof holding for $0 < p <1$. Note that the case $p=1$ is trivial as the epidemic is a homogeneously mixing epidemic with mean number of contacts made by each individual being $\mu_G + \mu_L$.

For $h=1,2,\ldots$, let $S^{(h)}$ denote the size of the susceptibility set of a typical individual in a household of size $h$. The probability that an individual with $X_L=x_L$ contacts the same individual twice in the household converges to 0 as the household size $h \to \infty$. Therefore for large $h$, the probability an individual contacts a given individual in their household via a local infection is approximately $\mu_L/(h-1)$. It is then straightforward to couple the construction of $S^{(h)}$ to a branching process with offspring distribution $V_h \sim {\rm Bin} (h-1, \mu_L/(h-1))$, with $V_h \convD \tilde{V} \sim {\rm Po} (\E[X_L])$ as $h \to \infty$.

Let $\tilde{\mathcal{B}}$ denote the branching process with offspring distribution $\tilde{V}$ and let $\tilde{S}$ denote the total size of the branching process $\tilde{\mathcal{B}}$. Then for $0 \leq s \leq 1$, the probability generating function of $\tilde{S}$ satisfies
\begin{align} \label{eq:hinf:1} \E \left[ s^{\tilde{S}} \right] = f_{\tilde{S}} (s) &= s \E \left[ f_{\tilde{S}} (s)^{\tilde{V}} \right] \nonumber \\
& = s \exp \left( - \mu_L \left[ 1 -  f_{\tilde{S}} (s) \right] \right).
\end{align}
 It follows that $z^{(h,0)} \to \tilde{z}$ as $h \to \infty$, where $\tilde{z}$ satisfies
\begin{align} \label{eq:hinf:2} 1 - \tilde{z} = f_{\tilde{S}} (\exp(-\mu_G \tilde{z})), \end{align}
and it remains to show that $\tilde{z} = z_{\rm hom}(\alpha)$, where $\alpha=\mu_G+\mu_L$.

We set $s = \exp (- \mu_G \tilde{z})$ in \eqref{eq:hinf:1}, and then using \eqref{eq:hinf:2}, we have that
\begin{align*}
 f_{\tilde{S}} \left( {\rm e}^{-\mu_G \tilde{z}} \right) & = {\rm e}^{-\mu_G \tilde{z}} \exp \left( -\mu_L \left[ 1 -  f_{\tilde{S}} \left( {\rm e}^{-\mu_G \tilde{z}} \right)  \right] \right)  \\
 &= {\rm e}^{-\mu_G \tilde{z}} \exp \left( -\mu_L \tilde{z} \right) =  \exp \left( -[ \mu_G + \mu_L] \tilde{z} \right).
\end{align*}
Therefore $\tilde{z}$ solves
\[ \tilde{z} = 1 -  f_{\tilde{S}} \left( {\rm e}^{-\mu_G \tilde{z}} \right)  = 1 - \exp \left( -[ \mu_G + \mu_L] \tilde{z} \right)=1-\exp(-\alpha \tilde{z}), \]
which is the defining equation for $z_{\rm hom}(\alpha)$.  Therefore, $\tilde{z} = z_{\rm hom}(\alpha)$, as required.

\section{Comparison of variance expressions} \label{app:sigma2}

In this appendix we prove that the expressions for $\sigma^2$ in \eqref{eq:sigma1} and \eqref{eq:sigma2} are equivalent.
The first step, using \eqref{eq:sigma1},   is to note that,
\begin{align} \label{eq:clt:var:C1}
\sigma^2 &= \frac{1}{h} \left[\var (R_1 (\tau)) + b(\tau)^2 \var (G_1 (\tau)) +   b(\tau)^2 \var (Y_1 (\tau)) \right.  \\
& \qquad \left. +2 b (\tau) \cov (R_1 (\tau), G_1 (\tau)) -2 b (\tau) \cov (R_1 (\tau), Y_1 (\tau)) - 2 b (\tau)^2 \cov (G_1 (\tau), Y_1 (\tau))  \right]\nonumber.
\end{align}
We consider separately each of the variance and covariance terms on the right-hand side of \eqref{eq:clt:var:C1}.

Using exchangeability of individuals,
\begin{align} \label{eq:clt:var:R1}
\var (R_1 (\tau)) & = h \var (\chi_{11} (\tau)) + h (h-1) \cov (\chi_{11} (\tau), \chi_{12} (\tau)) \nonumber \\
& = h \nu_R (\tau) [ 1-\nu_R (\tau)] + h (h-1) \cov (\chi_{11} (\tau), \chi_{12} (\tau)) .
\end{align}
Similarly,
\begin{align} \label{eq:clt:var:G1}
\var (G_1 (\tau)) & = h \var (\chi_{11} (\tau) X_{G,(1,1)}) + h (h-1) \cov (\chi_{11} (\tau)X_{G,(1,1)}, \chi_{12} (\tau)X_{G,(1,2)}).
\end{align}
Since $\chi_{11} (\tau)^2=\chi_{11} (\tau)$, 
\begin{align} \label{eq:clt:var:G2}
\var (\chi_{11} (\tau) X_{G,(1,1)})  &= \E \left[\chi_{11} (\tau) X_{G,(1,1)}^2 \right] - \nu_R (\tau)^2 \mu_G^2 \nonumber\\
&= \nu_R (\tau) [ \sigma_G^2 + \mu_G^2]  -  \nu_R (\tau)^2 \mu_G^2 \nonumber \\
&= \nu_R (\tau) [1-\nu_R (\tau)] \mu_G^2 + \nu_R (\tau) \sigma_G^2.
\end{align}

An observation similar to that made in \cite{BN03}, Section 4, is that, conditional upon $\chi_{11} (\tau)=0$ (individual $(1,1)$'s susceptibility set is not contacted when each members of the population is exposed to $\tau$ units of global infectious pressure), $(X_{G,(1,1)}, X_{L,(1,1)})$ and $\chi_{12} (\tau)$ are independent, so
\[
\E\left[(1-\chi_{11} (\tau))(1-\chi_{12} (\tau))X_{G,(1,1)}X_{G,(1,2)}\right]=\mu_G^2\E\left[(1-\chi_{11} (\tau))(1-\chi_{12} (\tau))\right].
\]
Since also
\begin{align*}
&\E\left[(1-\chi_{11} (\tau))(1-\chi_{12} (\tau))X_{G,(1,1)}X_{G,(1,2)}\right]\\
&\; =\E\left[X_{G,(1,1)}X_{G,(1,2)}\right]-2\E\left[\chi_{11} (\tau) X_{G,(1,1)}X_{G,(1,2)}\right]+\E\left[\chi_{11} (\tau) X_{G,(1,1)} \chi_{12} (\tau) X_{G,(1,2)}\right]\\
&\; =\mu_G^2-2\mu_G\E\left[\chi_{11} (\tau) X_{G,(1,2)}\right]+\E\left[\chi_{11} (\tau) X_{G,(1,1)} \chi_{12} (\tau) X_{G,(1,2)}\right],
\end{align*}
it follows that
\begin{align}
\label{equ:EchiX1}
&\E\left[\chi_{11} (\tau) X_{G,(1,1)} \chi_{12} (\tau) X_{G,(1,2)}\right]\nonumber\\
&\qquad=\mu_G^2\left\{\E\left[(1-\chi_{11} (\tau))(1-\chi_{12} (\tau))\right] -1\right\}+2\mu_G \E\left[\chi_{11} (\tau) X_{G,(1,2)}\right]\nonumber\\
&\qquad=\mu_G^2\E\left[\chi_{11} (\tau)  \chi_{12} (\tau)   \right]+2\mu_G \cov \left(\chi_{11} (\tau), X_{G,(1,2)} \right).
\end{align}
Thus,
\begin{align} \label{eq:clt:var:G3}
&\cov (\chi_{11} (\tau)X_{G,(1,1)}, \chi_{12} (\tau)X_{G,(1,2)})\nonumber\\
&\qquad\qquad  = \E \left[\chi_{11} (\tau) X_{G,(1,1)} \chi_{12} (\tau) X_{G,(1,2)}  \right] - \nu_R (\tau)^2 \mu_G^2 \nonumber \\
&\qquad\qquad= \mu_G^2 \E\left[\chi_{11} (\tau)  \chi_{12} (\tau)   \right] + 2 \mu_G \cov \left(\chi_{11} (\tau),X_{G,(1,2)} \right)- \nu_R (\tau)^2 \mu_G^2   \nonumber \\
&\qquad\qquad= \mu_G^2 \cov (\chi_{11} (\tau), \chi_{12} (\tau)) + 2 \mu_G  \cov \left(\chi_{11} (\tau),X_{G,(1,2)} \right).
\end{align}
Hence, substituting~\eqref{eq:clt:var:G2} and~\eqref{eq:clt:var:G3} into~\eqref{eq:clt:var:G1},
\begin{align} \label{eq:clt:var:G4}
\var (G_1 (\tau)) & = h \nu_R (\tau) [ 1-\nu_R (\tau)] \mu_G^2 +h \nu_R (\tau) \sigma_G^2 \nonumber \\
& \qquad + h (h-1)  \mu_G^2 \cov (\chi_{11} (\tau), \chi_{12} (\tau)) + 2 h (h-1) \mu_G  \cov \left(\chi_{11} (\tau),X_{G,(1,2)} \right) \nonumber \\
& = \mu_G^2 \var (R_1 (\tau)) + h \nu_R (\tau) \sigma_G^2  +  2 h (h-1) \mu_G  \cov \left(\chi_{11} (\tau),X_{G,(1,2)} \right).
\end{align}
Since $Y_1 (\tau) \sim {\rm Po} (h \tau)$, we have that $\var (Y_1 (\tau)) = h \tau$.

Turning to the covariance terms, since $\chi_{11} (\tau)$ and $X_{G,(1,1)}$ are independent,
\begin{align*}
\cov (R_1 (\tau), G_1 (\tau)) &= h \cov \left( \chi_{11} (\tau), \chi_{11} (\tau) X_{G,(1,1)} \right) + h (h-1) \cov \left( \chi_{11} (\tau), \chi_{12} (\tau) X_{G,(1,2)} \right) \\
&= h \mu_G \nu_R (\tau) [ 1- \nu_R (\tau)] + h (h-1)  \cov \left( \chi_{11} (\tau), \chi_{12} (\tau) X_{G,(1,2)} \right).
\end{align*}

A similar argument to the derivation of~\eqref{equ:EchiX1} yields
\begin{equation*}
\E\left[\chi_{11} (\tau)  \chi_{12} (\tau) X_{G,(1,2)}\right]=\mu_G \E\left[\chi_{11} (\tau)  \chi_{12} (\tau)   \right]+\cov \left(\chi_{11} (\tau), X_{G,(1,2)} \right),
\end{equation*}
so
\begin{align*}
\cov \left( \chi_{11} (\tau), \chi_{12} (\tau) X_{G,(1,2)} \right) & = \E \left[ \chi_{11} (\tau) \chi_{12} (\tau) X_{G,(1,2)} \right]  - \mu_G \nu_R (\tau)^2 \\
& = \mu_G \E \left[ \chi_{11} (\tau) \chi_{12} (\tau)\right] + \cov (\chi_{11} (\tau),X_{G,(1,2)})  - \mu_G \nu_R (\tau)^2  \\
& = \mu_G  \cov (\chi_{11} (\tau), \chi_{12} (\tau)) + \cov (\chi_{11} (\tau),X_{G,(1,2)}).
\end{align*}
Hence,
\begin{align} \label{eq:clt:covRG}
\cov (R_1 (\tau), G_1 (\tau)) &= h \mu_G \nu_R (\tau) [1-\nu_R (\tau)] + h (h-1)  \mu_G  \cov (\chi_{11} (\tau), \chi_{12} (\tau)) \nonumber\\
&\qquad + h (h-1) \cov (\chi_{11} (\tau),X_{G,(1,2)}) \nonumber \\
&= \mu_G \var (R_1 (\tau)) + h (h-1)  \cov (\chi_{11} (\tau),X_{G,(1,2)}).
\end{align}

Next, we have that
\begin{align}
\label{eq:clt:covRYa}
\cov (R_1 (\tau), Y_1 (\tau) & = h \cov (\chi_{11} (\tau), \zeta_{11} (\tau)) + h (h-1)  \cov (\chi_{11} (\tau), \zeta_{12} (\tau)).
\end{align}
Since $\chi_{11} (\tau) =1$ if $ \zeta_{11} (\tau) >0$ and $\zeta_{11}(\tau) \sim {\rm Po}(\tau)$, we have that
\begin{align}
\label{eq:clt:covRYb}
\cov (\chi_{11} (\tau), \zeta_{11} (\tau)) = \E [\chi_{11} (\tau) \zeta_{11} (\tau)] - \tau \nu_R (\tau)
&= \E [ \zeta_{11} (\tau)] - \tau \nu_R (\tau) \nonumber \\
&= \tau [1- \nu_R (\tau)].
\end{align}
Note that $\P(\chi_{11} (\tau)=1\,|\,\zeta_{12} (\tau)=k)= \P(\chi_{11} (\tau)=1\,|\,\zeta_{12} (\tau)=1)$ $(k=1,2,\dots)$, so
\begin{align*}
\E\left[\chi_{11} (\tau)\zeta_{12} (\tau)\right]&=\sum_{k=0}^{\infty} k \P(\chi_{11} (\tau)=1\,|\,\zeta_{12} (\tau)=k)\P(\zeta_{12} (\tau)=k)\\
&=\tau\P(\chi_{11} (\tau)=1\,|\,\zeta_{12} (\tau)=1),
\end{align*}
since $\zeta_{12}(\tau) \sim {\rm Po}(\tau)$.  Also, $\P((1,2) \notin \mathcal{S}_{1,1}\,|\,S_{11}=i)=\frac{h-i}{h-1}$ $(i=1,2,\dots, h)$, so
\[
\P(\chi_{11} (\tau)=0\,|\,\zeta_{12} (\tau)=1)=\sum_{i=1}^h \P(S_{11}=i)\left[\frac{h-i}{h-1}\right] \re^{-i \tau}.
\]
Noting that $\nu_R(\tau)=1-f_S(\re^{-\tau})=1-\sum_{i=1}^h \P(S_{11}=i)\re^{-i\tau}$, a short calculation yields
\[
\E\left[\chi_{11} (\tau)\zeta_{12} (\tau)\right]=\tau\left[\nu_R(\tau)+\frac{\nu_R (\tau) + \nu_R^\prime (\tau)-1}{h-1}\right],
\]
whence
\begin{align}
\label{eq:clt:covRYc}
\cov (\chi_{11} (\tau), \zeta_{12} (\tau)) &= \E [\chi_{11} (\tau) \zeta_{12} (\tau)] - \tau \nu_R (\tau) \nonumber \\
&= \frac{\tau [\nu_R (\tau) + \nu_R^\prime (\tau)-1]}{h-1}.
\end{align}
Substituting~\eqref{eq:clt:covRYb} and~\eqref{eq:clt:covRYc} into~\eqref{eq:clt:covRYa} yields
\begin{align} \label{eq:clt:covRY}
\cov (R_1 (\tau), Y_1 (\tau)) & = h\tau [1- \nu_R (\tau)] + h\tau [\nu_R (\tau) + \nu_R^\prime (\tau)-1] \nonumber \\
&= h \tau \nu_R^\prime (\tau).
\end{align}
Also
\begin{align} \label{eq:clt:covGY}
\cov (G_1 (\tau), Y_1 (\tau) )& = \mu_G \cov (R_1 (\tau), Y_1 (\tau)  )= h \mu_G \tau \nu_R^\prime (\tau).
\end{align}

Substituting \eqref{eq:clt:var:R1},  \eqref{eq:clt:var:G4}, $\var (Y_1 (\tau))=h \tau$, \eqref{eq:clt:covRG},  \eqref{eq:clt:covRY} and \eqref{eq:clt:covGY} into \eqref{eq:clt:var:C1} yields
\begin{align}
\sigma^2 &= \frac{1}{h}\left[\var (R_1 (\tau)) +b (\tau)^2 \mu_G^2 \var (R_1 (\tau)) + h b(\tau)^2 \nu_R (\tau) \sigma_G^2 \right.  \nonumber \\
& \qquad + 2 h (h-1) b(\tau)^2 \mu_G \cov (\chi_{11} (\tau), X_{G,(1,2)}) + b (\tau)^2 h \tau  + 2 b (\tau) \mu_G \var (R_1 (\tau)) \nonumber \\
& \qquad  \left. + 2 b(\tau) h (h-1)  \cov (\chi_{11} (\tau), X_{G,(1,2)})  - 2 b(\tau) [1 + \mu_G b(\tau) ] h \tau \nu_R^\prime (\tau) \right]. \nonumber
\end{align}
Now,
\begin{align*}
1 + \mu_G b (\tau) = 1 + \frac{\mu_G \nu_R^\prime (\tau)}{1- \mu_G \nu_R^\prime (\tau)} = \frac{1}{1- \mu_G \nu_R^\prime (\tau)},
\end{align*}
so
\begin{align*}
\left\{1 + \mu_G b (\tau) \right\} \nu_R^\prime (\tau) =\frac{\nu_R^\prime (\tau)}{1- \mu_G \nu_R^\prime (\tau)} = b (\tau).
\end{align*}
Hence, using $\tau = \nu_G (\tau) = \mu_G \nu_R (\tau)$,
\begin{align}  \label{eq:clt:sigmaF}
\sigma^2 &= \frac{1}{h} \left[ \{ 1+ \mu_G b (\tau)\}^2 \var (R_1 (\tau)) + b (\tau)^2 \{ h \nu_R (\tau) \sigma_G^2 + h \tau - 2 h \tau\} \right. \nonumber  \\
& \qquad \left.  +2 h (h-1) b (\tau) [1 + \mu_G b (\tau) ] \cov (\chi_{11} (\tau), X_{G,(1,2)}) \right] \nonumber \\
& = \left( 1 + b (\tau) \mu_G \right)^2 \nu_R (\tau) [1 - \nu_R (\tau)] + b (\tau)^2 \nu_R (\tau) [\sigma_G^2 - \mu_G] \nonumber \\
& \qquad +(h-1) \left[  \left( 1 + b (\tau) \mu_G \right)^2 \cov (\chi_{11} (\tau), \chi_{12} (\tau) + 2 b(\tau) \left( 1 +  \mu_G b (\tau)  \right) \cov \left(\chi_{11} (\tau),  X_{G,(1,2)}\right)   \right]. \nonumber \\
\end{align}
The right-hand side of \eqref{eq:clt:sigmaF} agrees with \eqref{eq:sigma2} completing the proof.

The expression for $\sigma^2$ given by~\eqref{eq:sigma} follows after a little algebra by substituting~\eqref{eq:clt:covRY}, \eqref{eq:clt:covGY} and $\var (Y_1 (\tau))=h \tau$ into~\eqref{eq:clt:var:C1} and noting that~\eqref{eq:clt:var:G4} and~\eqref{eq:clt:covRG} imply
\begin{equation*}
\var (G_1 (\tau))=2\mu_G\cov(R_1(\tau), G_1(\tau))-\mu_G^2\var(R_1(\tau))+h\nu_R(\tau)\sigma_G^2.
\end{equation*}

\end{appendix}

{\bf Funding:}
T.B. was supported in part by the Swedish Research Council (grant 2020-0474).

\end{document}